\documentclass[reqno]{amsart}

\usepackage{geometry}
\usepackage{amsmath,amssymb,amsthm,amsfonts}
\usepackage{physics}
\allowdisplaybreaks[4]
\usepackage{hyperref}
\usepackage{mathrsfs}
\usepackage{appendix}
\usepackage{graphicx}
\usepackage{setspace}
\usepackage{mathtools}
\usepackage{latexsym}
\usepackage{mleftright}
\usepackage{cite}
\usepackage{color}

\newtheorem{lemma}{Lemma}[section]
\newtheorem{theorem}{Theorem}[section]

\newtheorem{proposition}{Proposition}[section]
\newtheorem{remark}{Remark}[section]

\numberwithin{equation}{section}
\begin{document}
\title[NSCH equations modeling the dynamics of contact line]{Well-posedness of Navier-Stokes/Cahn-Hilliard equations modeling the dynamics of contact line in a channel}
\thanks{{\it Keywords}: Navier-Stokes/Cahn-Hilliard equations; generalized Navier boundary condition;  relaxation  boundary condition; existence; uniqueness}
\thanks{{\it AMS Subject Classification}: 35A01, 35Q35, 76D03, 76T05}%
\author[Shijin Ding]{Shijin Ding}
\address[S. Ding]{School of Mathematical Sciences, South China Normal University,
Guangzhou, 510631, China}
\email{dingsj@scnu.edu.cn}
\author[Yinghua Li]{Yinghua Li}
\address[Y. Li]{School of Mathematical Sciences, South China Normal University,
Guangzhou, 510631, China}
\email{yinghua@scnu.edu.cn}
\author[Zhilin Lin]{Zhilin Lin}
\address[Z. Lin]{School of Mathematical Sciences, South China Normal University,
Guangzhou, 510631, China}
\email{zllin@m.scnu.edu.cn}
\author[Yuanxiang Yan]{Yuanxiang Yan}
\address[Y. Yan]{School of Mathematical Sciences, South China Normal University,
Guangzhou, 510631, China}
\email{yuanxiangyan@m.scnu.edu.cn}


\begin{abstract}
In this paper, we study the contact line problem in a channel. Precisely, we consider the incompressible Navier-Stokes/Cahn-Hilliard system with generalized Navier boundary condition and relaxation boundary condition in a channel, which is the phase field model for the moving contact line problem in fluid mechanics. We establish the existence and uniqueness of local-in-time strong solution to this initial boundary value problem in 2D. To our knowledge, this is the first result to give the local-in-time well-posedness of Navier-Stokes/Cahn-Hilliard system with generalized Navier boundary condition and relaxation boundary condition.  This result provides a rigorous mathematical analysis to confirm that the physical and numerical results by Qian-Wang-Sheng [Phys. Rev. E 68 (2003), 016306, 1-15; J. Fluid Mech. 564 (2006), 333-360] are well-posed and reasonable.
\end{abstract}

\maketitle
\tableofcontents
\vspace{-5mm}


\newpage
\section{Introduction}
\subsection{Formulation of the problem}
The so called phase field model for the 2D binary mixture of incompressible fluid reads as  the following Navier-Stokes/Cahn-Hilliard (NSCH) equations:
\begin{align}\label{111111}
\begin{cases}
{\bf u}_t+({\bf u}\cdot\nabla){\bf u}+ \nabla p = {\rm div}\mathbb{S}({\bf u})+\mu\nabla \phi, &~~\text{in}~~\Omega,
\\
{\rm div}{\bf u}=0, &~~\text{in}~~\Omega,
\\
\phi_t+({\bf u}\cdot\nabla)\phi=M\Delta\mu, &~~\text{in}~~\Omega,
\\
\mu=-\Delta\phi+f(\phi), &~~\text{in}~~\Omega,
\end{cases}
\end{align}
where $\Omega \subset \mathbb{R}^2$ is a bounded smooth domain, the strain tensor is defined by $\mathbb{S}({\bf u})=\nu\left(\nabla{\bf u} + \nabla^\top{\bf u}\right)$, and ${\bf u}=(u_1,u_2),\phi,\mu,M,p$ denote the mean velocity of the fluid mixture, the concentration difference of the two components, chemical potential, mobility coefficient and the pressure, respectively. In our setting, the free energy density is taken as $f(\phi)=\phi^3-\phi$.
 It is noted that system \eqref{111111} is also called ``model H", see \cite{H-H} for details.

To understand the behavior of fluid near the boundary, from the physical view \cite{Q-W-S-03,Q-W-S-03-02,Q-W-S-04,Q-W-S-06}, the generalized Navier boundary condition (GNBC) and the relaxation boundary condition may be more suitable to formulate the motion of fluid on boundaries, which is stated as follows
\begin{align}\label{BC-1}
\begin{cases}
{\bf u}\cdot {\bf n} = 0, &~~\text{on}~~\Gamma ,
\\
\partial_{\bf n}\mu =0, &~~\text{on}~~\Gamma ,
\\
\beta {\bf u}_{\boldsymbol \tau} + \left(\mathbb{S}({\bf u}) \cdot {\bf n}\right)_{\boldsymbol \tau}= \mathcal{L}(\phi) \nabla_{\boldsymbol \tau} \phi, &~~\text{on}~~\Gamma ,
\\
\phi_t+ {\bf u}_{\boldsymbol \tau} \cdot \nabla_{\boldsymbol \tau} \phi = -\hat{\Gamma} \mathcal{L}(\phi), &~~\text{on}~~\Gamma,
\end{cases}
\end{align}
 where $\Gamma$ is the boundary of $\Omega$, $\mathcal{L}(\phi)=\partial_{\bf n} \phi + \frac{\partial\gamma_{fs}(\phi)}{\partial \phi} $ in which $\gamma_{fs}(\phi)$ is the fluid-solid interfacial free energy per unit area, $\beta > 0 $ is the slip coefficient, $\hat{\Gamma}$ is a positive phenomenological parameter.
${\boldsymbol \tau}$ and ${\bf n}$ represent the unit tangent vector and the unit outward normal vector on the solid boundary $\Gamma$, and
\begin{equation*}
{\bf u}_{\boldsymbol \tau} = {\bf u}-({\bf n}\cdot {\bf u}){\bf n},~~~\partial_{\bf n}={\bf n}\cdot \nabla,~~~\nabla_{\boldsymbol \tau} =  \nabla-( {\bf n}\cdot\nabla){\bf n} .
\end{equation*}

In particular, if we consider the system \eqref{111111} with \eqref{BC-1} in a channel, the case which has been studied in lots of physical and numerical literatures, see \cite{Q-W-S-03,Q-W-S-03-02,Q-W-S-04,Q-W-S-06} ,  $\Omega=\mathbb{T}\times (-1,1)$,  i.e., the horizontal direction is periodic, then
$$
\Gamma = \mathbb{T} \times(\{ -1\} \cup \{ 1\}).
$$
Moreover, there are ${\bf u}_{\boldsymbol \tau} = (u_{1},0) $  and $\nabla_{\boldsymbol \tau}= (\partial_{x},0)$ on $\Gamma$, $\partial_{\bf n}=\pm\partial_{y}$  on $\{y=\pm 1\}$. The outward unit normal vector of $\Gamma$ becomes:
$$
{\bf n} = (0,\pm 1)~ \text{on}~\{y=\pm 1\}.
$$

In this case, boundary condition \eqref{BC-1} can be rewritten as
\begin{align}\label{BC}
\begin{cases}
u_2=0,& \ \mathrm{on} \ \Gamma,
\\
\partial_{\bf n} \mu =0, &~~\text{on}~~\Gamma,
\\
\beta u_1 + \partial_{\bf n} u_1= \mathcal{L}(\phi)\partial_x  \phi, &~~\text{on}~~\Gamma,
\\
\phi_t+u_1\partial_x \phi = - \mathcal{L}(\phi), &~~\text{on}~~\Gamma,
\end{cases}
\end{align}
where  $\mathcal{L}(\phi)= \partial_{\bf n} \phi + \frac{\partial\gamma_{fs}(\phi)}{\partial \phi}$ on $\Gamma$, and we have set $M=\hat{\Gamma}=\nu=1$ without loss of generality. Inspired by \cite{Q-W-S-03,Q-W-S-04,Q-W-S-06}, we take the fluid-solid interfacial free energy density $\gamma_{fs}(\phi)$ as
\begin{equation}\label{gfs}
\gamma_{fs}(\phi)=-\dfrac{\gamma}{2}\cos \theta_s \sin\left(\frac{\pi \phi}{2}\right),
\end{equation}
in which $\gamma>0$ is the fluid-fluid interface tension constant and $\theta_s$ is the static contact angle which is determined by the Young equation $\gamma_{fs}(\phi_{-}) = \gamma_{fs}(\phi_{+}) + \gamma \cos \theta_{s}$, see Figure 1.

\begin{center}\label{fig1}
\begin{figure}[h]
\includegraphics[width=0.6\linewidth]{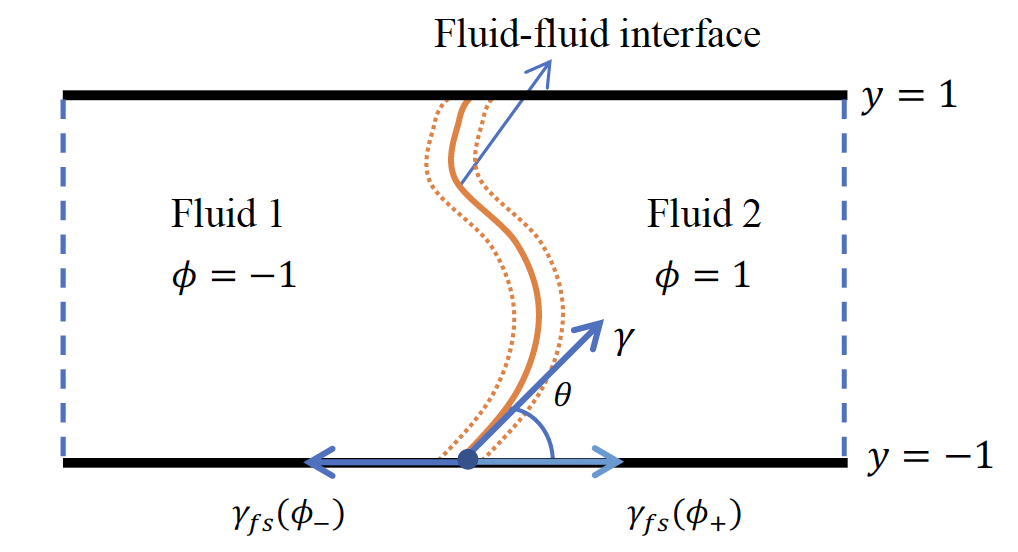}
\caption{Contact angle formed by fluid-fluid interface with the solid boundary}
\end{figure}
\end{center}
\ \\[-3em]

For simplicity, we set $ \frac{\gamma}{2}\cos \theta_s = 1 $. We impose the initial conditions as follows
\begin{equation}\label{initial}
u(x, y, 0)=u_0(x, y),\ \phi(x, y, 0)=\phi_0(x, y)
\end{equation}
and define
\begin{equation}\label{I-mu-0}
\mu_0= -\Delta\phi_0+f(\phi_0).
\end{equation}

Our goal of this paper is to establish the local-in-time well-posedness for problem \eqref{111111} , \eqref{BC} and \eqref{initial}.

\subsection{Literatures review}
There are lots of works about the well-posedness for system \eqref{111111} with different type of boundary conditions, in particular the no-slip boundary condition ${\bf u}\big|_{\Gamma}=0$ for velocity  and no-flux boundary condition for $\phi$ and $\mu$. The well-posedness for the system \eqref{111111} with the boundary condition
\begin{align}\label{222222}
\begin{cases}
{\bf u} = 0, &~~\text{on}~~\Gamma,
\\
\partial_{\bf n}\phi = \partial_{\bf n}\mu =0, &~~\text{on}~~\Gamma
\end{cases}
\end{align}
has been widely studied, see \cite{A,B,M,G-G,G-M-T} for instance.
However, as stated in \cite{D-D}, for immiscible two-phase flows with a moving contact line which is defined as the intersection of the fluid-fluid interface with the solid wall, the no-slip boundary condition does not apply. In fact, the authors showed that for incompressible Newtonian fluids, if one assigns no-slip boundary condition on smooth rigid walls on which the interface of fluid-fluid slips, then there is a velocity discontinuity at the moving contact line, and the tangential force exerted by the fluids on the solid surface in the vicinity of the contact line becomes infinite. This type of singularity is well known as the moving contact line singularity.

To treat with the singularity, Qian-Wang-Sheng \cite{Q-W-S-03, Q-W-S-04, Q-W-S-06} proposed the GNBC (i.e., the boundary condition \eqref{BC-1}$_{1,3}$) and relaxation boundary condition (i.e., the boundary condition \eqref{BC-1}$_{4}$) which can be obtained from the laws of thermodynamics and the principle of minimum dissipation Onsager energy. They found that GNBC can formulate the relative slipping between fluids and solid wall well.
In particular, for single phase fluid, or away from the contact line, \eqref{BC-1}$_{1,3}$  are reduced to
\begin{align}
\begin{cases}
{\bf u}\cdot {\bf n} = 0, &~~\text{\rm on}~~\Gamma ,
\\
\beta {\bf u}_{\boldsymbol \tau}+ \left(\mathbb{S}({\bf u}) \cdot {\bf n}\right)_{\boldsymbol \tau}= 0, &~~\text{\rm on}~~\Gamma ,
\end{cases}
\end{align}
which are called as the Navier boundary condition (NBC), more about the Navier-Stokes equations with NBC can be found in  \cite{A-R, B-Y, C-M-R, K-Z, L-D, M-R, S-M, S-S, X-X-07, X-X-13} and the references therein. Indeed, the NBC is valid for case of away from the contact line region \cite{T-T}. However, in the interfacial region, due to the singularity near the contact line (or contact angle), the NBC fails to describe the motion of fluids near contact line. Therefore, the formulation for the motion of fluids in the interfacial region is nontrivial. Let us point out that the item $\mathcal{L}(\phi)\nabla_{\boldsymbol\tau} \phi$ on the right hand of $(\ref{BC-1})_{3}$ denotes the uncompensated Young stress, which is resulted from the deviation of the fluid-fluid interface from its static configuration, see \cite{Q-W-S-03-02, Q-W-S-06} and references therein for details. In addition, a special case of \eqref{BC-1}$_{4}$:
$
\partial_{\bf n}\phi = 0
$
leads to a static contact line with contact angle $\frac{\pi}{2}$.




In the special case of slow dynamics, one can neglect the effect of the flow, i.e., ${\bf u}=0$, the problem is reduced to
the Cahn-Hilliard equations $(\ref{111111})_{3,4}$ with the homogeneous Neumann boundary condition $(\ref{BC-1})_{2}$ and
the relaxation boundary condition $(\ref{BC-1})_{4}$.
To our knowledge, the only theoretical analysis was given by Chen-Wang-Xu \cite{C-W-X}, in which the existence of global solutions was established for the case without fluid effect. 
Precisely, by adding a surface diffusion term $\delta \Delta_{\boldsymbol \tau}\phi$, \cite{C-W-X} considered the following dynamic boundary condition to replace the relaxation boundary condition \eqref{BC-1}$_{4}$:
$$
\phi_t= -\partial_{\bf n}\phi + \delta\Delta_{\boldsymbol \tau}\phi- k(\phi)\ \ \  \text{on}\ \Gamma.
$$
Compared with the original boundary condition \eqref{BC-1}$_{4}$, a surface diffusion term $ \delta\Delta_{\boldsymbol \tau}\phi$ has been added which may provide certain regular effect near boundaries, see also Gal \cite{G-1} for more details about this point.

The well-posedness theory of Cahn-Hilliard equations with dynamic boundary condition has been studied in several works. Racke-Zheng \cite{R-Z} obtained the global existence and uniqueness of strong solution to this initial boundary value problem in $\mathbb{R}^N(N=1,2,3)$. In \cite{P-R-Z}, Pr${\rm\ddot{u}}$ss et al. established the maximal $L^p$-regularity and the existence of a global attractor. Later, convergence of the solution to Cahn-Hilliard equation with dynamic boundary condition  to equilibrium was proved by Wu-Zheng \cite{W-Z}, as $t \to \infty$. For the free energy density $f(\phi)$ with arbitrary growths, Miranville-Zelik \cite{M-Z} obtained the existence and uniqueness of solution and constructed a robust family of exponential attractors in 3D. Recently, Liu-Wu\cite{L-W} proved the existence and uniqueness of global weak/strong solution to this initial boundary value problem with or without surface diffusion. It is worth mentioning that this is a well-posedness result for a fourth-order parabolic equation subject to a fourth-order dynamic boundary condition\cite{L-W}. Other progresses can be found in \cite{C-F-W,M-W,W} and so on.

There are lots of works devoted to the numerical analysis for NSCH with GNBC and relaxation boundary condition, see \cite{S-Y-Y,W-H,X-D-Y,Y-Y-1,Y-Y-2,Z-K-Y-L-S-1,Z-K-Y-L-S-2} and the references therein. However, to our knowledge, there are no results about the existence and uniqueness of strong/weak solution to 
the NSCH system \eqref{111111} with GNBC and relaxation boundary condition \eqref{BC-1} due to the difficulties in the analysis involving the relaxation boundary conditions. So the known works up to now are devoted to consider the following dynamic boundary condition
\begin{equation}\label{dynamic}
\phi_t+ {\bf u}_{\boldsymbol \tau} \cdot \nabla_{\boldsymbol \tau} \phi = -\partial_{\bf n}\phi + \delta\Delta_{\boldsymbol \tau}\phi - \alpha\phi - k(\phi), ~\text{on}~\Gamma,
\end{equation}
which can be viewed as a regularization of the relaxation boundary condition $(\ref{BC-1})_{4}$ by adding surface diffusion $\delta\Delta_{\boldsymbol \tau}\phi$ and a damping term $-\alpha\phi$. For NSCH system \eqref{111111} with GNBC and \eqref{dynamic},
Gal-Grasselli-Miranville \cite{G-G-M} obtained the existence of a suitable global energy solution and established the convergence of such solution to a single equilibrium as $t\to\infty$. Later, Gal-Grasselli-Wu \cite{G-G-W} considered this initial boundary value problem with unmatched densities and proved the existence of a global weak solution for arbitrary initial data in both two and three dimensions. For compressible fluid, Cherfils et al. \cite{C-F-M-M-P-P} proved the existence of global-in-time weak solutions with finite initial energy. Recently, Chen et al. \cite{C-H-H-S} got the existence and uniqueness of the local strong solution to the Navier-Stokes/Allen-Cahn equation  with GNBC and relaxation boundary condition in 3D. In addition, Gal-Grasselli-Poiatti \cite{G-G-P} showed the existence of global weak solutions to the Allen-Cahn-Navier-Stokes-Voigt Systems with GNBC and dynamic boundary condition.

\subsection{Main result}
To state our main result, we will introduce some notations used throughout the paper. For a domain $\Omega \subset \mathbb{R}^2$, $ 1 \le p \le + \infty $, and $ s >0 $, we denote
\begin{align*}
\begin{cases}
\langle f \rangle = \frac{1}{|\Omega|}\int_\Omega f \mathrm{d}x {\rm d}y,~~~~~\langle f , g\rangle = \int_\Omega f g\mathrm{d}x {\rm d}y, \\
L^p=L^p(\Omega), ~~~~~W^{s,p}=W^{s,p}(\Omega),~~~~~H^s= W^{s,2},\\
\mathbb{D} = \left\{ ({\bf u},\phi)\in H^2\times H^3\,\big|\, {\bf u},\phi \,{\rm satisfy}\, (\ref{BC})_{1,2,3}~ {\rm and}~ (\ref{111111})_2 \right\},
\end{cases}
\end{align*}
where $W^{s,p}$ is fractional Sobolev space if $ s $ is non-integer, one can see the  references \cite{D-P-V, B-M} for details.
\vskip 3mm
Our main result can be stated as follows:
\begin{theorem}
\label{TH}
Suppose that $\Omega=\mathbb{T}\times (-1,1)$, $ \mu_0\in H^3$, $ \partial_{x}^{2}\phi_{0} \in H^{2} $ and $({\bf u}_0,\phi_{0})\in \mathbb{D}$, then there is a positive constant $T_{0}$, which may depend on the initial data such that the problem \eqref{111111}, \eqref{BC}, \eqref{initial} admits a unique local-in-time strong solution $({\bf u},\phi ) $ in $[0,T_0]$ with
\begin{align}
&{\bf u} \in L^{\infty}(0,T_{0};H^2) \cap  L^{2}(0,T_{0};H^3), ~\,\,{\bf u}_t\in L^{\infty}(0,T_{0};L^2) \cap  L^{2}(0,T_{0};H^1),  \nonumber \\
&\phi \in L^{\infty}(0,T_{0};H^3),~~\Delta \phi \in L^{2}(0,T_{0};H^2),~\,\,\mu \in  L^{\infty}(0,T_{0};H^3) \cap  L^{2}(0,T_{0};H^4) , \nonumber \\
&\phi_t \in L^{\infty}(0,T_{0};H^1) \cap  L^{2}(0,T_{0};H^2), ~\,\, \mu_t \in  L^{2}(0,T_{0};H^1).\nonumber
\end{align}
\end{theorem}

\begin{remark}
To our knowledge, Theorem \ref{TH} is the first result (even in a channel) to establish the well-posedness  of the strong solutions to the
NSCH system \eqref{111111}
with the original GNBC \eqref{BC-1}, which is proposed in \cite{Q-W-S-03, Q-W-S-03-02}.
\end{remark}

\begin{remark}
Compared with the previous works \cite{C-F-M-M-P-P, G-G-M, G-G-W}, Theorem \ref{TH} treated with the case without surface diffusion and damping in $\mathcal{L}(\phi)$. Therefore, in additional to the difficulties inherited from the Model H, the main challenges in the mathematical analysis
come from the relaxation boundary condition in (\ref{BC-1}). See section \ref{sec:2}.
\end{remark}

\begin{remark}
Up to now, even for the Cahn-Hilliard equation, except for \cite{C-W-X}, most works are focused
on the dynamic boundary condition
$$
\partial_t\phi=-\mathcal{L}(\phi), \quad
\mathcal{L}(\phi)=-\delta\Delta_{\boldsymbol \tau}\phi+\partial_{\bf n}\phi+\alpha\phi+k(\phi),\quad
\text{on}~~\Gamma,
$$
where $\delta, \alpha$ are fixed positive constants, see \cite{ G-1, M-Z, P-R-Z, R-Z} for example. In addition, the effect of the flow was neglected in \cite{C-W-X},  
therefore there is not strong nonlinear term ${\bf u}_{\boldsymbol \tau}\cdot \nabla_{\boldsymbol \tau}\phi $ on relaxation boundary condition.
\end{remark}

\begin{remark}  Our result implies that there is a damping effect coming from the relaxation boundary condition {\rm(}see \eqref{damping}{\rm)}, which confirms that the physical and numerical results in \cite{Q-W-S-03,Q-W-S-03-02,Q-W-S-04,Q-W-S-06} are well-posed and reasonable.
\end{remark}

The structure of this paper is as follows. In Section \ref{sec:2}, we state the key ideas of proof and recall some auxiliary results that will be used later. In Section \ref{sec:3},  we do some a priori estimates for the original problem \eqref{111111}, \eqref{BC}, \eqref{initial}, which will be used for deriving the estimates uniform in $\delta$ to the $\delta$-approximate problem \eqref{A-NSCH}-\eqref{A-IC}. In Section \ref{sec:5}, we construct a $\delta$-approximate problem and obtain the estimates uniform in $\delta$. Finally, we finish the proof of Theorem 1.1 by letting $\delta \to 0$ in Section \ref{sec:55}.


\section{Key ideas of proof and some preliminaries}\label{sec:2}

\subsection{Key ideas of proof}
In this subsection, let us state the key ideas in our arguments.

To establish the well-posedness for problem \eqref{111111}, \eqref{BC}, \eqref{initial}, the key points are the {\it a priori estimates} and the construction of the suitable approximate solution.

In order to derive the {\it a priori estimates}, compared with the previous works \cite{C-F-M-M-P-P, C-P,  G-G-M, G-G-W, Y, Y-L-Z}, the main trouble comes from the boundary condition \eqref{BC}, i.e., there is no diffusion in the relaxation boundary condition in \eqref{BC}$_4$, namely,
$$\partial_t \phi+u_1\partial_x\phi = \partial_{\bf n} \phi + \cdots,\ \text{on}\ \Gamma$$
from which one can not get smoothing effect near the boundaries explicitly. In addition, the coupling effect in the system \eqref{111111} does not provide the positive structure. Indeed, the above relaxation boundary condition can be viewed formally as a transport equation on boundaries, then one can not desire the regularity of $\phi$ can be achieved via the transport effect.

To overcome this difficulties, we have the following two observations.

The first key observation is that it holds from  \eqref{111111}$_{3}$ and \eqref{BC}$_{{4}}$ that
\begin{align}\label{Key-B}
\Delta\mu = -\partial_{\bf n}\phi - \gamma_{fs}'(\phi), ~~~~\text{on}~~\Gamma.
\end{align}
With this, we have the following equality which plays a key role to achieve the desired (especially higher order) estimates
\begin{align}\label{damping}
-\langle \partial_t^i\partial_x^j\Delta \mu,\partial_t^i\partial_x^j \Delta \phi\rangle=&\langle  \partial_t^i\partial_x^j \nabla \Delta \mu,\partial_t^i\partial_x^j  \nabla \phi\rangle-\int_{\Gamma}\partial_t^i\partial_x^j \Delta \mu \, \partial_t^i\partial_x^j  \partial_{\bf n} \phi \mathrm{d}x
\nonumber\\
=&\langle \partial_t^i\partial_x^j  \nabla \Delta \mu,\partial_t^i\partial_x^j \nabla \phi\rangle+\int_{\Gamma} |\partial_t^i\partial_x^j \partial_{\bf n} \phi|^2 \mathrm{d}x+ \int_{\Gamma} \partial_t^i\partial_x^j\gamma_{fs}'(\phi) \partial_t^i\partial_x^j \partial_{\bf n} \phi {\rm d}x,
\end{align}
where the boundary terms $\|\partial_t^i\partial_x^j \partial_{\bf n}\phi \|_{L^2(\Gamma)}$ provide the positive structure in the estimates, which implies that there is a damping effect in the relaxation boundary condition.

Our second key observation is that, in order to obtain the estimates uniform in $\delta$,  the term $\delta \partial_{x}^{2}\phi$ in \eqref{A-NSCH}$_{3}$ is very important to match the boundary condition \eqref{A-BC}$_{4}$ when constructing the approximate solution, see Remark \ref{R-add} for more details.

With the {\it a priori estimates} at hand, it turns to construct the suitable approximate solution to problem \eqref{111111}, \eqref{BC}, \eqref{initial}. As mentioned before, the main challenges are
the relaxation boundary condition and the coupling effect in system. The key point in construction is to look for a suitable approximate scheme, which is suitable for {\it a priori estimates} with fixed point argument.

Formally, there are two main parts to be considered in the approximation: the fluid part (Navier-Stokes part) and the Cahn-Hilliard part. The trouble coming from the relaxation boundary condition mainly affects the Cahn-Hilliard part.
To search for a suitable approximate scheme, inspired by \cite{C-W-X, M-Z}, we should consider the following toy models corresponding to \eqref{111111}:
\begin{align}\label{toy-2}
\begin{cases}
 \psi_t-\delta \partial_{x}^{2} \psi =G_2, &\ \text{on}\ \Gamma,\\
\psi|_{t=0}=\psi_0=:\phi_0|_{\Gamma}, &\ \text{on}\ \Gamma,
\end{cases}
\end{align}
and
\begin{align}\label{toy-3}
\begin{cases}
\delta \phi_t-\Delta \phi=G_{3}, &\ \ \text{in} \ \Omega,\\
\phi|_{\Gamma}=\psi, &  \ \ \text{on} \ \Gamma, \\
\phi|_{t=0}=\phi_0, &\ \ \text{in} \ \Omega,
\end{cases}
\end{align}
where $\delta>0$ is a constant. Indeed, $G_{3} = \mu+ \cdots$, thus we need to consider additionally the following equation for $\mu$:
\begin{align}\label{toy-1}
\begin{cases}
\mu-\delta\Delta \mu =G_1,\ &\text{in}\ \Omega,\\
\partial_{\bf n} \mu \big|_{\Gamma}=0,\ &\text{on}\ \Gamma.
\end{cases}
\end{align}

Since we are working in the channel domain,  which ensures that the above toy models are well-defined and well-posed, see section 4 for details.
The above toy models modify the Cahn-Hilliard part of the system.

With the toy models above, the toy model for the fluid part  can be constructed as
\begin{align}\label{toy-4}
\begin{cases}
{\bf u}_t-\mathrm{div}\mathbb{S}({\bf u})+\nabla p=G_4, &\text{in}\ \Omega,\\
\mathrm{div} {\bf u}=0,& \text{in}\ \Omega,\\
u_2=0, & \text{on}\ \Gamma,\\
\beta u_1+\partial_{\bf n} u_1=g_4,& \text{on}\ \Gamma,\\
{\bf u}|_{t=0}={\bf u}_0, & \text{in}\ \Omega.
\end{cases}
\end{align}

The above toy models come from the linearization processes to solve the $\delta$-approximate problem  by fixed point theorem. In addition, the above toy models allow us to derive the uniform-in-$\delta$ estimates for the approximate solution, which ensure that the fixed point argument can be processed and the damping effect of the relaxation boundary condition can be preserved in the approximation, see Section \ref{sec:5} for details.

Each toy model in \eqref{toy-2}-\eqref{toy-4} is well-posed. With this, the iteration to construct the approximate solution can be done  by solving $\eqref{toy-1}$ and $\eqref{toy-2}$, then solving $\eqref{toy-3}$ and $\eqref{toy-4}$.

\subsection{Some preliminaries}

This subsection will state some preliminaries, which will be used in this paper.
\begin{lemma}
\label{uuu}
Let $ \Omega \subset \mathbb{R}^N(N=2,3)$ be a bounded domain with smooth boundary $\Gamma$. Then there exists a constant $ C > 0 $ independent of $ {\bf u}$, such that
\begin{align}
\| {\bf u}\|_{H^s} \le C\left( \| \nabla \times {\bf u} \|_{H^{s-1}} + \| {\rm div} {\bf u} \|_{H^{s-1}} + \| {\bf u}\cdot{\bf n} \|_{H^{s-\frac{1}{2}}(\Gamma)} +  \| {\bf u} \|_{L^{2}} \right),
\end{align}
for any  $ {\bf u} \in H^s(\Omega)$, $s \ge 1$.
\end{lemma}
\begin{proof}[\bf Proof]
See \cite{B-B,F-T}.
\end{proof}
The following lemma is a Gagliardo-Nirenberg type inequality, which can be found in \cite{B-M}.
\begin{lemma}\label{GN}
Let $ \Omega $ be either $ \mathbb{R}^N $ or a half space or a Lipschitz bounded domain in $ \mathbb{R}^N $. Assume that the real number $0\le s_1<s<s_2$, $1\le p_1,p_2,p\le\infty$ and $0<\theta<1$ satisfy the relations
\begin{align}
s= \theta s_1 + (1-\theta)s_2 ~\,\,\text{and} ~\,\,\dfrac{1}{p}=\dfrac{\theta}{p_1}+\dfrac{1-\theta}{p_2}.
\end{align}
Then there exists a constant $C$ depending on $s_1,s_2,p_1,p_2,\theta$ and $\Omega$ such that
\begin{align}
\|f\|_{W^{s,p}} \le C \|f\|^{\theta}_{W^{s_1,p_1}} \|f\|^{1-\theta}_{W^{s_2,p_2}}, ~\forall f \in W^{s_1,p_1}(\Omega)\cap W^{s_2,p_2}(\Omega)
\end{align}
holds if and only if
\begin{align}
s_2 ~\text{is an integer}~ \ge 1, p_2=1 ~\text{and}~ s_2-s_1\le1-\dfrac{1}{p_1}
\end{align}
fail.

In particular, for $k\in \mathbb{N}$, it holds that 
\begin{align}\label{1-2}
\|f\|_{H^{{k+1/2}}}^{2} \le C \|f\|_{H^{k}} \|f\|_{H^{k+1}}.
\end{align}

\end{lemma}

\begin{lemma}[Trace Imbedding Theorem]\label{Trace}
Let $ \Omega $ be a domain in $\mathbb{R}^N$ satisfying the $C^{k-1,1}$-regularity condition, and suppose $u\in W^{k,p}(\Omega)$, $p>1$, $k$ an integer. If $l\le k-1$, then there is trace operator ${\bf tr}$ such that
\begin{align}
{\bf tr}: ~  W^{k,p}(\Omega) \hookrightarrow  W^{k-l-\frac{1}{p},p}(\Gamma)
\end{align}
holds. If $u|_{\Gamma} = \varphi$, we denote ${\bf tr}(u)=\varphi$.

\end{lemma}
\begin{proof}[\bf Proof]
 See Ch.2, Theorem 5.5 in \cite{N}.
\end{proof}

\begin{lemma}[\cite{B-C,C-M-R,K}]\label{lemma-curl-u}
Let $ \Omega \subset \mathbb{R}^N(N=2,3)$ and $ \omega = \nabla \times {\bf u} $. Suppose ${\bf u}\in H^2$ and ${\bf u} \cdot {\bf n}=0$ on $\Gamma$. Then we have
\begin{align}
\left[ \mathbb{S}({\bf u})\cdot {\bf n} \right] \cdot {\boldsymbol \tau} = (\omega\times {\bf n})\cdot {\boldsymbol \tau} - 2{\bf u} \cdot \dfrac{\partial {\bf n}}{\partial {\boldsymbol \tau}}, ~~~~\text{for}~~ N=3,
\end{align}
and
\begin{align}
\left[ \mathbb{S}({\bf u})\cdot {\bf n} \right] \cdot {\boldsymbol \tau} = \omega - 2{\bf u} \cdot \dfrac{\partial {\bf n}}{\partial {\boldsymbol \tau}}, ~~~~\text{for}~~ N=2,
\end{align}
where $\left | \dfrac{\partial {\bf n}}{\partial {\boldsymbol \tau}}\right|$ is the normal curvature in the ${\boldsymbol \tau}$ direction when $N=3$, $\dfrac{\partial {\bf n}}{\partial {\boldsymbol \tau}}= \kappa {\boldsymbol \tau} $,  while $\kappa$ is the curvature of $\Gamma$ when $N=2$.
\end{lemma}

\begin{lemma}[Korn's inequality,\cite{L-D}] \label{L-Korn}
 There exists a constant $C>0$, there holds
\begin{align}
\| \mathbb{S}({\bf u}) \|_{L^2(\Omega)}\ge C\| \nabla{\bf u} \|_{L^2(\Omega)},
\end{align}
for any  $ {\bf u} \in V=\left\{ {\bf u} \in H^{1}(\Omega) \,\big| \,{\rm div}{\bf u} = 0,  {\bf u}\cdot{\bf n} = 0 ~~{\rm on}~\Gamma \right\} $.
\end{lemma}
We also need the following product law, which can be found in \cite{G-T}.
\begin{lemma}\label{L-Product}
The following estimates hold on sufficiently smooth subsets of  $\mathbb{R}^{N}$:

{\rm 1.} Let $0\le r \le s_{1} \le s_{2}$ be such that $s_{1}>\frac{N}{2}$. Let $f \in H^{s_{1}}$, $g \in H^{s_{2}}$. Then $fg\in H^{r}$ and
\begin{align}\label{Product-1}
\| fg \|_{H^r}\le C\| f \|_{H^{s_{1}}} \| g\|_{H^{s_{2}}}.
\end{align}

{\rm 2.} Let $0\le r \le s_{1} \le s_{2}$ be such that $s_{2}>r+\frac{N}{2}$. Let $f \in H^{s_{1}}$, $g \in H^{s_{2}}$. Then $fg\in H^{r}$ and
\begin{align}\label{Product-2}
\| fg \|_{H^r}\le C\| f \|_{H^{s_{1}}} \| g\|_{H^{s_{2}}}.
\end{align}
\end{lemma}


\section{A priori estimates}\label{sec:3}

For any $\delta>0$, we  will construct the $\delta$-approximate problem \eqref{A-NSCH}-\eqref{A-IC} to approximate the original problem \eqref{111111}, \eqref{BC}, \eqref{initial}. The purpose of this section is to obtain {\it a priori estimates} of the original problem \eqref{111111}, \eqref{BC}, \eqref{initial}, which will be used for deriving the estimates uniform in $\delta$ to the $\delta$-approximate problem \eqref{A-NSCH}-\eqref{A-IC} in subsection \ref{subs-d}.

To derive the {\it a priori estimates}, it is important to obtain the damping effect from relaxation boundary condition. Before giving the details of estimates, let us give an observation about relaxation boundary condition.

It follows from $(\ref{111111})_{3}$ that
\begin{align*}
\phi_{t} + u_1 \partial_{x}\phi  +  u_{2} \partial_{y}\phi = \Delta\mu, ~\,\,\,\,\text{in}~~\Omega .
\end{align*}
Taking the trace for the above equation to yield that
\begin{align}\label{ttrace}
\phi_{t} + u_1 \partial_{x}\phi = \Delta\mu, ~\,\,\,\,\text{on}~~\Gamma ,
\end{align}
since $ u_{2}= 0$ on $\Gamma$. Thus, we find from \eqref{ttrace} and $(\ref{BC})_{4}$ that
\begin{align}\label{N-BC}
\Delta \mu =-\partial_{\bf n}\phi-\gamma'_{fs}(\phi),~\,\,\,\, \text{on}\ \Gamma.
\end{align}

Now it is ready to derive the {\it a priori estimates}. To this end, let us define the following energy functional
\begin{align}\label{tE}
\mathcal{E}(t) = 1 + \| {\bf u}\|_{H^{2}}^{2} + \| {\bf u}_{t}\|_{L^{2}}^{2} + \| \phi\|_{H^{3}}^{2} + \| \phi_{t}\|_{H^{1}}^{2} +  \| \mu \|_{H^{3}}^{2},
\end{align}
and dissipative functional
\begin{align}
\mathcal{D}(t) = \| {\bf u}\|_{H^{3}}^{2} + \| {\bf u}_{t}\|_{H^{1}}^{2} + \| \Delta\phi\|_{H^{2}}^{2} + \| \phi_{t}\|_{H^{2}}^{2} +  \| \mu \|_{H^{4}}^{2} +   \| \mu_{t} \|_{H^{1}}^{2} +   \| \partial_{\bf n}\phi \|_{H^{2}(\Gamma)}^{2}.
\end{align}
The following proposition is about {\it a priori estimates} for \eqref{111111}, \eqref{BC}, \eqref{initial}.
\begin{proposition}\label{apriori}
Under the assumptions of Theorem $\ref{TH}$, then for $0<T<1$,  it holds that
\begin{align}
\mathcal{E}(t)   + \int_{0}^{T} \mathcal{D}(t) {\rm d}t\le C   \int_{0}^{T} \mathcal{E}(t)^{8} {\rm d}t + C,
\end{align}
where $C$ is a constant that depends on $\Omega$, $\beta$, $ \hat{\mathcal{E}}(0)$ and the initial value  $\left(\| {\bf u}_{0}\|_{H^{2}}, \| \phi_{0}\|_{H^{3}}, \| \mu_{0}\|_{H^{3}}\right)$.
\end{proposition}

Proposition \ref{apriori} is the direct consequence of the following lemmas (Lemma $\ref{energy}$-$\ref{phi-H33}$). The rest work of this section is to establish Lemma $\ref{energy}$-$\ref{phi-H33}$.

We start with the basic energy.

\begin{lemma} \label{energy}
Under the assumptions of Theorem $\ref{TH}$, then for $0<T<1$, it holds that
\begin{align}\label{basic-E}
\sup\limits_{0\le t\le T}\hat{\mathcal{E}}(t) +\int_{0}^{T} \hat{\mathcal{D}}(t) {\rm d}t \le C\hat{\mathcal{E}}(0),
\end{align}
where
\begin{align*}
\hat{\mathcal{E}}(t) & := \| {\bf u}\|_{L^2}^{2}+\|  \phi \|_{H^1}^{2}+ \frac{1}{2}\| \phi^2-1 \|_{L^2}^2 + \int_{\Gamma}\gamma_{fs}(\phi) {\rm d}x,
\\
\hat{\mathcal{D}}(t) & := \| \mathbb{S}({\bf u}) \|_{L^2}^{2}+\beta \| u_1\|_{L^2(\Gamma)}^{2}+\| \nabla \mu \|_{L^2}^{2}+\| \mathcal{L}(\phi)\|_{L^2(\Gamma)}^{2},
\end{align*}
and  $C$ is a positive constant.
\end{lemma}

\begin{proof}[\bf{Proof}]  Multiplying $(\ref{111111})_1$ by ${\bf u}$ and integrating over $\Omega$ by parts, we obtain
\begin{align}\label{E-1}
\frac{1}{2}\frac{{\rm d}}{{\rm d}t}\| {\bf u}\|_{L^2}^{2}+\| \mathbb{S}({\bf u}) \|_{L^2}^{2}+\beta \| u_1\|_{L^2(\Gamma)}^{2} = \int_{\Omega} \mu {\bf u} \cdot \nabla \phi {\rm d}x{\rm d}y+\int_{\Gamma} \mathcal{L}(\phi)  u_{1} \partial_{x} \phi  {\rm d}x.
\end{align}
Multiplying $(\ref{111111})_3$ by $\mu$ and integrating by parts, we get
\begin{align}\label{E-2}
&\frac{1}{2}\frac{{\rm d}}{{\rm d}t} \left[ \| \nabla \phi \|_{L^2}^{2}+ \frac{1}{2}\| \phi^2-1 \|_{L^2}^2 + \int_{\Gamma}\gamma_{fs}(\phi) {\rm d}x\right]+\| \nabla \mu \|_{L^2}^{2}+\| \mathcal{L}(\phi)\|_{L^2(\Gamma)}^{2}
\nonumber \\
 &= -\int_{\Omega} \mu {\bf u} \cdot \nabla \phi {\rm d}x{\rm d}y-\int_{\Gamma} \mathcal{L}(\phi) u_{1}\partial_{x} \phi {\rm d}x.
\end{align}
Combining $(\ref{E-1})$ with $(\ref{E-2})$, we have
\begin{align*}
 \sup \limits_{0\le t\le T}\hat{\mathcal{E}}_{1} (t)+ \int_{0}^{T} \left(\| \mathbb{S}({\bf u}) \|_{L^2}^{2}+\beta \| u_1\|_{L^2(\Gamma)}^{2}+\| \nabla \mu \|_{L^2}^{2}+\| \mathcal{L}(\phi)\|_{L^2(\Gamma)}^{2} \right) {\rm d}t= \hat{\mathcal{E}}_{1}(0),
\end{align*}
where
\begin{align}\label{h-E-1}
\hat{\mathcal{E}}_{1} (t) : =  \| {\bf u} \|_{L^{2}}^{2} + \| \nabla \phi \|_{L^2}^{2}+ \frac{1}{2}\| \phi^2-1 \|_{L^2}^2 + \int_{\Gamma}\gamma_{fs}(\phi) {\rm d}x.
\end{align}
Note that
\begin{align}\label{00}
 \frac{\rm d}{{\rm d}t} \langle \phi   \rangle & = \frac{1}{|\Omega|}\int_\Omega \phi_t   \mathrm{d}x{\rm d}y  =\frac{1}{|\Omega|} \int_\Omega \Delta\mu   \mathrm{d}x{\rm d}y  -\frac{1}{|\Omega|} \int_\Omega {\bf u} \cdot \nabla \phi   \mathrm{d}x {\rm d}y
\nonumber \\
& = - \frac{1}{|\Omega|}\int_{\Gamma} \partial_{\bf n} \mu  {\rm d}x  + \frac{1}{|\Omega|}\int_\Omega {\rm div} {\bf u} \,  \phi   \mathrm{d}x{\rm d}y -\frac{1}{|\Omega|}\int_{\Gamma}  {\bf u}\cdot {\bf n} \, \phi   \mathrm{d}x = 0.
\end{align}
We use Poincar${\rm \acute{e}}$ inequality  and obtain that
\begin{align}\label{dijie-1}
 \| \phi \|_{L^{2}}^{2}  \le    \left\| \phi -  \langle \phi   \rangle  \right\|_{L^{2}}^{2} + \left\|  \langle \phi   \rangle  \right\|_{L^{2}}^{2} \le    C \| \nabla\phi\|_{L^{2}}^{2} + C\left|\langle \phi_{0} \rangle \right|^{2}  \le C\hat{\mathcal{E}}(0) .
 \end{align}
This completes the proof of Lemma $\ref{energy}$.
\end{proof}

\begin{lemma}\label{lemma-phi-1}
Under the assumptions of Theorem $\ref{TH}$, then for $0<T<1$,  it holds that
\begin{align}\label{D-1}
 \|  \phi_{t}\|_{H^{1}}^{2} +  \|{\bf u}_t  \|^2_{L^2} + \| {\bf u}   \|^2_{H^2}  + \int_{0}^{T} \mathcal{D}_{1}(t) {\rm d}t \le C  \int_{0}^{T} \mathcal{E}(t)^{4} {\rm d}t + C,
\end{align}
where
\begin{align}
\mathcal{D}_{1}(t) := \| \nabla \mu_{t} \|_{L^{2}}^{2} + \| \partial_{\bf n}\phi_{t}\|_{L^{2}(\Gamma)}^{2} + \|\partial_{x}\phi_{t}\|_{L^{2}(\Gamma)}^{2}+  \|\nabla{\bf u}_t \|^2_{L^2} + \| {\bf u}   \|^2_{H^3} + \beta \| \partial_{t}u_{1}\|^2_{L^2(\Gamma)},
\end{align}
and  $C$ is a constant that depends on $\Omega$, $\beta$, $ \hat{\mathcal{E}}(0)$ and the initial value  $\left(\| {\bf u}_{0}\|_{H^{2}}, \| \phi_{0}\|_{H^{3}}, \| \mu_{0}\|_{H^{3}}\right)$.
\end{lemma}

\begin{proof}[\bf{Proof}]
The proof will be completed by several steps.

{\bf Step 1. Estimate for $\nabla \phi_t$.}

Differentiating $(\ref{111111})_{3,4}$ and $(\ref{N-BC})$ with respect to $t$ respectively, we get
\begin{align}\label{1-t}
\begin{cases}
\Delta\mu_{t}  = \phi_{tt} + {\bf u}_{t} \cdot \nabla\phi  + {\bf u} \cdot \nabla\phi_{t}  , & {\rm in}~~ \Omega,
\\
 \Delta\phi_{t} =  -\mu_{t}  +f_{t}, & {\rm in}~~ \Omega,
\\
\Delta\mu_{t} = - \partial_{\bf n}\phi_{t} - \gamma_{fs}^{(2)}(\phi)\phi_{t},  & {\rm on}~~ \Gamma.
\end{cases}
\end{align}
We now consider
\begin{align}\label{E-1-11}
-\langle \Delta\mu_{t} , \Delta\phi_{t}  \rangle=  \langle  \nabla\Delta\mu_{t},  \nabla\phi_{t} \rangle-   \int_{\Gamma}  \Delta\mu_{t}   \partial_{\bf n}\phi_{t}  {\rm d}x =: I_{1} + I_{2}.
\end{align}
In order to obtain the information of $I_{1}$, we put $(\ref{1-t})_{1}$ into $I_{1}$ and find
\begin{align}\label{E-1-0}
I_{1} & =   \int_{\Omega}  \nabla \left(  \phi_{tt} + {\bf u}_{t} \cdot \nabla\phi  + {\bf u} \cdot \nabla \phi_{t}  \right)   \cdot  \nabla\phi_{t}  {\rm d}x{\rm d}y
\nonumber \\
& =  \frac{1}{2} \frac{\rm d}{{\rm d}t} \| \nabla \phi_{t}\|_{L^{2}}^{2} +  \int_{\Omega}  \left(  \nabla {\bf u}_{t} \cdot \nabla\phi +    {\bf u}_{t} \cdot \nabla^{2}\phi + \nabla{\bf u} \cdot \nabla \phi_{t} +{\bf u} \cdot \nabla^{2} \phi_{t} \right)   \cdot  \nabla\phi_{t}  {\rm d}x{\rm d}y
\nonumber \\
& =  \frac{1}{2} \frac{\rm d}{{\rm d}t} \| \nabla \phi_{t}\|_{L^{2}}^{2} + \int_{\Omega}  \left(  \nabla {\bf u}_{t} \cdot \nabla\phi +    {\bf u}_{t} \cdot \nabla^{2}\phi + \nabla{\bf u} \cdot \nabla \phi_{t} \right)   \cdot  \nabla\phi_{t}  {\rm d}x{\rm d}y,
\end{align}
where we have used the following fact
\begin{align}
 \int_{\Omega}{\bf u} \cdot \nabla^{2} \phi_{t}   \cdot  \nabla\phi_{t}  {\rm d}x{\rm d}y  & =   \frac{1}{2}\int_{\Omega}{\bf u} \cdot \nabla \left( |\nabla\phi_{t}|^{2} \right) {\rm d}x{\rm d}y
\nonumber \\
& = -\frac{1}{2}\int_{\Omega} {\rm div}{\bf u} \left( |\nabla\phi_{t}|^{2} \right) {\rm d}x{\rm d}y  + \frac{1}{2}\int_{\Gamma} {\bf u}\cdot {\bf n} \left( |\nabla\phi_{t}|^{2} \right) {\rm d}x = 0.
\end{align}
Put $(\ref{1-t})_{3}$ into $I_{2}$ to obtain
\begin{align}\label{E-1-1}
I_{2}  =  \int_{\Gamma} \left(  \partial_{\bf n}\phi_{t} + \gamma_{fs}^{(2)}(\phi)\phi_{t} \right) \partial_{\bf n}\phi_{t}  {\rm d}x = \| \partial_{\bf n}\phi_{t}\|_{L^{2}(\Gamma)}^{2} +   \int_{\Gamma} \gamma_{fs}^{(2)}(\phi)\phi_{t}  \partial_{\bf n}\phi_{t}  {\rm d}x.
\end{align}
On the other hand, it follows from integration by parts that
\begin{align}\label{E-1-2}
-\langle \Delta\mu_{t},  \Delta\phi_{t} \rangle = \int_{\Omega} \Delta\mu_{t}  \left(  \mu_{t}  -f_{t} \right)  {\rm d}x{\rm d}y = - \| \nabla \mu_{t} \|_{L^{2}}^{2} + \int_{\Omega} \nabla\mu_{t} \cdot \nabla f_{t}  {\rm d}x{\rm d}y,
\end{align}
where we have used $(\ref{1-t})_{2}$ and $(\ref{BC})_{2}$.

Substituting $(\ref{E-1-0})$, $(\ref{E-1-1})$ and $(\ref{E-1-2})$ into $(\ref{E-1-11})$, we arrive at
\begin{align}\label{C-1}
& \frac{1}{2} \frac{\rm d}{{\rm d}t} \| \nabla \phi_{t}\|_{L^{2}}^{2} +  \| \nabla \mu_{t} \|_{L^{2}}^{2} + \| \partial_{\bf n}\phi_{t}\|_{L^{2}(\Gamma)}^{2}
\nonumber\\
& = -  \int_{\Omega}  \left(  \nabla {\bf u}_{t} \cdot \nabla\phi +    {\bf u}_{t} \cdot \nabla^{2}\phi + \nabla{\bf u} \cdot \nabla \phi_{t} \right)   \cdot  \nabla\phi_{t}  {\rm d}x{\rm d}y
\nonumber\\
&~\,\,\, +\int_{\Omega} \nabla\mu_{t} \cdot \nabla f_{t}  {\rm d}x{\rm d}y -\int_{\Gamma} \gamma_{fs}^{(2)}(\phi)\phi_{t}  \partial_{\bf n}\phi_{t}  {\rm d}x
\nonumber\\
& \le \epsilon \left(  \| \nabla {\bf u}_{t} \|_{L^{2}}^{2} + \| \nabla {\bf u} \|_{L^{\infty}}^{2} \right)  +  \frac{1}{4} \| \nabla \mu_{t} \|_{L^{2}}^{2} +  \frac{1}{6} \| \partial_{\bf n}\phi_{t}\|_{L^{2}(\Gamma)}^{2}  +    \| \nabla\phi \|_{L^{\infty}}^{2}\| \nabla \phi_{t}\|_{L^{2}}^{2}
\nonumber\\
& ~\,\,\,  + \| \nabla \phi_{t} \|_{L^{2}}^{4} +   \| {\bf u}_{t} \|_{L^{4}}   \| \nabla^{2}\phi \|_{L^{4}}\| \nabla \phi_{t}\|_{L^{2}} + C\| \nabla f_{t} \|_{L^{2}}^{2} + C\| \phi_{t} \|_{L^{2}(\Gamma)}^{2}
\nonumber\\
& \le \epsilon \left(  \| \nabla {\bf u}_{t} \|_{L^{2}}^{2} + \| \nabla {\bf u} \|_{H^{2}}^{2} \right)  +  \frac{1}{4} \| \nabla \mu_{t} \|_{L^{2}}^{2} +  \frac{1}{6} \| \partial_{\bf n}\phi_{t}\|_{L^{2}(\Gamma)}^{2}  +    \| \nabla\phi \|_{H^{2}}^{2}\| \nabla \phi_{t}\|_{L^{2}}^{2}
\nonumber\\
&  ~\,\,\, + \| \nabla \phi_{t} \|_{L^{2}}^{4} +   \| {\bf u}_{t} \|_{L^{2}}^{\frac12}  \| {\bf u}_{t} \|_{H^{1}}^{\frac12}  \| \nabla^{2}\phi \|_{H^{1}}\| \nabla \phi_{t}\|_{L^{2}}  + C\| \nabla f_{t} \|_{L^{2}}^{2} + C\| \phi_{t} \|_{H^{1}}^{2}
\nonumber\\
& \le \epsilon \left(  2 \| \nabla {\bf u}_{t} \|_{L^{2}}^{2} + \| \nabla {\bf u} \|_{H^{2}}^{2} \right)  +  \frac{1}{4} \| \nabla \mu_{t} \|_{L^{2}}^{2} +  \frac{1}{6} \| \partial_{\bf n}\phi_{t}\|_{L^{2}(\Gamma)}^{2}  + C\mathcal{E}(t)^{2}
\nonumber\\
&  ~\,\,\, +  C \| {\bf u}_{t} \|_{L^{2}}^{2} +  C\| \nabla^{2}\phi \|_{H^{1}}^{2} \| \nabla \phi_{t}\|_{L^{2}}^{2}  + C\| \nabla f_{t} \|_{L^{2}}^{2} + C\| \phi_{t} \|_{H^{1}}^{2}
\nonumber\\
& \le \epsilon \left(  2 \| \nabla {\bf u}_{t} \|_{L^{2}}^{2} + \| \nabla {\bf u} \|_{H^{2}}^{2} \right)  +  \frac{1}{4} \| \nabla \mu_{t} \|_{L^{2}}^{2} +  \frac{1}{6} \| \partial_{\bf n}\phi_{t}\|_{L^{2}(\Gamma)}^{2}  + C\mathcal{E}(t)^{3},
\end{align}
where we have used Lemma $\ref{GN}$, Lemma $\ref{Trace}$ and the following fact
\begin{align*}
\| \nabla f_{t} \|_{L^{2}}^{2} & = \left\| (3\phi^{2}-1)\nabla\phi_{t} + 6\phi\nabla\phi\phi_{t} \right\|_{L^{2}}^{2}
\\
& \le C \left(   \| \phi\|_{L^{\infty}}^{4}  + 1 \right)  \| \nabla \phi_{t} \|_{L^{2}}^{2} + C \| \phi\|_{L^{\infty}}^{2} \| \nabla \phi\|_{L^{\infty}}^{2} \| \phi_{t}\|_{L^{2}}^{2}
 \\
& \le C  \left(   \| \phi\|_{H^{3}}^{4}  + 1 \right)  \| \phi_{t}\|_{H^{1}}^{2} \le C\mathcal{E}(t)^{3}.
\end{align*}

Meanwhile, since $\Gamma$ is one-dimensional manifold, it follows from $(\ref{BC})_{4}$, $(\ref{gfs})$, $H^{1}(\Gamma)\hookrightarrow L^{\infty}(\Gamma) $ and Lemma $\ref{Trace}$ that
\begin{align}\label{C-yb}
\|\partial_{x}\phi_{t}\|_{L^{2}(\Gamma)}^{2} &\le  \left\| -\partial_{x} u_{1} \partial_{x} \phi - u_{1} \partial_{xx}\phi - \partial_{x}\partial_{\bf n}\phi -  \gamma_{fs}^{(2)}(\phi) \partial_{x}\phi \right\|_{L^{2}(\Gamma)}^{2}
\nonumber \\
& \le  \| \partial_{x}u_{1} \|_{L^{2}(\Gamma)}^{2}  \| \partial_{x} \phi \|_{L^{\infty}(\Gamma)}^{2} +   \|u_{1} \|_{L^{\infty}(\Gamma)}^{2}  \| \partial_{xx} \phi \|_{L^{2}(\Gamma)}^{2}
\nonumber \\
& ~\,\,\, +   \| \partial_{x}\partial_{\bf n}\phi  \|_{L^{2}(\Gamma)}^{2} + C\| \partial_{x} \phi \|_{L^{2}(\Gamma)}^{2}
\nonumber \\
& \le  \|  {\bf u} \|_{H^{2} }^{2}  \|\partial_{x} \phi \|_{H^{1}(\Gamma)}^{2} +   \|  u_{1} \|_{H^{1}(\Gamma)}^{2}  \|  \phi \|_{H^{3}}^{2} +  C \| \phi  \|_{H^{3}}^{2}
\nonumber \\
&\le C  \| {\bf u} \|_{H^{2} }^{2} \|  \phi \|_{H^{3}}^{2} + C \|  \phi \|_{H^{3}}^{2}  \le  C\mathcal{E}(t)^{2}.
\end{align}

\ \\[-2em]

{\bf Step 2. Estimate for $ \nabla {\bf u}_{t}$. }

Differentiating $(\ref{111111})_1$ with respect to $t$, one has
\begin{align*}
{\bf u}  _{tt} + {\bf u}_{t} \cdot \nabla {\bf u} +  {\bf u}  \cdot \nabla {\bf u}_{t} + \nabla p_t   = {\rm div}\mathbb{S}({\bf u}_t  ) + \mu_{t} \nabla \phi + \mu \nabla\phi_{t} .
\end{align*}
Multiplying it by $ {\bf u}_t   $, and integrating over $\Omega$ by parts, combining with the boundary condition $(\ref{BC})_{1,3}$, we derive
\begin{align}\label{C-2}
& \dfrac{1}{2} \dfrac{\rm d}{{\rm d}t} \|{\bf u}_t  \|^2_{L^2} + \dfrac{1}{2}\|\mathbb{S}({\bf u}_t  ) \|^2_{L^2} + \beta \| \partial_{t}u_{1}\|^2_{L^2(\Gamma)}
\nonumber \\
& = -\int_{\Omega}   {\bf u}_{t} \cdot \nabla {\bf u} \cdot {\bf u}_{t} {\rm d}x{\rm d}y - \int_{\Omega}   {\bf u}_{t} \cdot \nabla\mu_{t} \phi {\rm d}x{\rm d}y + \int_{\Omega}    {\bf u}_{t} \cdot\nabla \mu  \phi_{t} {\rm d}x{\rm d}y
\nonumber \\
& ~\,\,\,+ \int_{\Gamma}  \left( \partial_{\bf n}\phi_{t} +\gamma_{fs}^{(2)}(\phi)\phi_{t} \right) \partial_{x} \phi  \partial_{t}u_{1} {\rm d}x + \int_{\Gamma}  \left( \partial_{\bf n}\phi + \gamma_{fs}'(\phi) \right) \partial_{x}\phi_{t}  \partial_{t}u_{1} {\rm d}x
\nonumber \\
& \le   \frac{1}{4} \| \nabla \mu_{t} \|_{L^{2}}^{2} +  \frac{1}{6} \| \partial_{\bf n}\phi_{t}\|_{L^{2}(\Gamma)}^{2} + \frac{1}{4}\|\partial_{x}\phi_{t}\|_{L^{2}(\Gamma)}^{2} + \| {\bf u}_{t} \|_{L^{2}} \| \nabla{\bf u}\|_{L^{4}} \| {\bf u}_{t} \|_{L^{4}} 
\nonumber \\
&~\,\,\, +\| \phi\|_{L^{\infty}}^{2} \| {\bf u}_{t} \|_{L^{2}}^{2} + \| {\bf u}_{t} \|_{L^{2}} \| \nabla\mu\|_{L^{4}} \| \phi_{t} \|_{L^{4}}+ C\| \phi_{t} \|_{L^{2}(\Gamma)}^{2}
\nonumber \\
& ~\,\,\, +  C \left( \|\partial_{x}\phi\|^2_{L^\infty(\Gamma)} +  \|\partial_{\bf n}\phi\|^2_{L^\infty(\Gamma)}  + 1 \right) \| \partial_{t}u_{1}\|^2_{L^2(\Gamma)} 
\nonumber \\
& \le   \frac{1}{4} \| \nabla \mu_{t} \|_{L^{2}}^{2} +  \frac{1}{6} \| \partial_{\bf n}\phi_{t}\|_{L^{2}(\Gamma)}^{2} + \frac{1}{4}\| \partial_{x}\phi_{t}\|_{L^{2}(\Gamma)}^{2} + \| {\bf u}_{t} \|_{L^{2}} \| \nabla{\bf u}\|_{L^{4}} \| {\bf u}_{t} \|_{L^{4}}
\nonumber \\
& ~\,\,\,   +C\| \phi\|_{H^{2}}^{2} \| {\bf u}_{t} \|_{L^{2}}^{2} + \| {\bf u}_{t} \|_{L^{2}} \| \nabla\mu\|_{H^{1}} \| \phi_{t} \|_{H^{1}}+ C\| \phi_{t} \|_{H^{1}}^{2}\nonumber \\
& ~\,\,\,+ C  \left( \|\partial_{x}\phi\|^2_{H^1(\Gamma)} +  \|\partial_{\bf n}\phi\|^2_{H^1(\Gamma)}  + 1 \right)\| \partial_{t}u_{1}\|^2_{L^2(\Gamma)}  
\nonumber \\
& \le   \frac{1}{4} \| \nabla \mu_{t} \|_{L^{2}}^{2} +  \frac{1}{6} \| \partial_{\bf n}\phi_{t}\|_{L^{2}(\Gamma)}^{2} + \frac{1}{4}\|\partial_{x}\phi_{t}\|_{L^{2}(\Gamma)}^{2} + \| {\bf u}_{t} \|_{L^{2}} \| \nabla{\bf u}\|_{L^{4}} \| {\bf u}_{t} \|_{L^{4}}  
 \nonumber \\
&~\,\,\,+ C\mathcal{E}(t)^{2}
 +C\mathcal{E}(t)  \| \partial_{t} u_{1}\|^2_{L^2(\Gamma)} 
\nonumber \\
& \le   \frac{1}{4} \| \nabla \mu_{t} \|_{L^{2}}^{2} +  \frac{1}{6} \| \partial_{\bf n}\phi_{t}\|_{L^{2}(\Gamma)}^{2}    + \frac{1}{4}\| \partial_{x}\phi_{t}\|_{L^{2}(\Gamma)}^{2}  + \frac{1}{6} \|\mathbb{S}({\bf u}_t  ) \|^2_{L^2} + C\mathcal{E}(t)^{3},
\end{align}
where we have used Lemma $\ref{Trace}$, Korn's inequality and the following fact
\begin{align*}
 &\| {\bf u}_{t} \|_{L^{2}} \| \nabla{\bf u}\|_{L^{4}} \| {\bf u}_{t} \|_{L^{4}}   \le  \| {\bf u}_{t} \|_{L^{2}} \| {\bf u}\|_{H^{2}} \| {\bf u}_{t} \|_{L^{2}}^{\frac{1}{2}} \| {\bf u}_{t} \|_{H^{1}}^{\frac{1}{2}} \le C\mathcal{E}(t)^{2}
+   \frac{1}{12} \|\mathbb{S}({\bf u}_t  ) \|^2_{L^2},
\\
&\mathcal{E}(t) \|\partial_{t}u_{1}\|^2_{L^2(\Gamma)}  \le C\mathcal{E}(t) \| {\bf u}_{t} \|_{H^{\frac{1}{2}}}^{2} \le C\mathcal{E}(t) \| {\bf u}_{t} \|_{L^{2}}  \| {\bf u}_{t} \|_{H^{1}} \le  C\mathcal{E}(t)^{3} +  \frac{1}{12} \|\mathbb{S}({\bf u}_t  ) \|^2_{L^2}.
 \end{align*}

\ \\[-2em]

{\bf Step 3. Estimates for ${\bf u} $ .}

 Let $\omega = \nabla \times {\bf u}$. Operating ${\nabla \times}$ to $(\ref{111111})_{1}$ and using  Lemma $\ref{lemma-curl-u}$, we have
\begin{align}\label{curl-u}
\begin{cases}
 - \Delta\omega   = -\omega  _t - \nabla \times \left( {\bf u}\cdot \nabla{\bf u}  - \mu \nabla\phi\right)  , & {\rm in}~~ \Omega,
\\
\omega = \pm \left[\beta u_{1} - \big(\partial_{\bf n}\phi +\gamma_{fs}'(\phi)\big)\partial_{x}\phi\right], & {\rm on}~~ \{y=\pm1\}.
\end{cases}
\end{align}
Applying the classical elliptic $H^k$-regularity theory  to $(\ref{curl-u})_1$ with Dirichlet boundary condition $(\ref{curl-u})_2$ and noting $\|\omega_{t}\|_{L^{2}} \le C\|\nabla {\bf u}_{t}\|_{L^{2}}$, we obtain
\begin{align*}
\| \omega   \|_{H^k}
& \le C(\beta) \left( \|\nabla {\bf u}  _t\|_{H^{k-2}} + \|\nabla \left( {\bf u}\cdot \nabla{\bf u}  + \mu \nabla\phi\right)   \|_{H^{k-2}} + \|u_{1}\|_{H^{k-\frac{1}{2}}(\Gamma)}   \right)
\nonumber \\
& ~\,\,\,+ C  \left\| \left(\partial_{\bf n}\phi +\gamma_{fs}'(\phi)\right)\partial_{x}\phi \right\|_{H^{k-\frac{1}{2}}(\Gamma)} .
\end{align*}
Thus, using Lemma $\ref{uuu}$, Lemma $\ref{Trace}$, the above inequality and $(\ref{Product-1})$ for $r=s_{1}=s_{2}=\frac32 > \frac 12$, we infer
\begin{align}\label{u-H33}
\| {\bf u}   \|^2_{H^3} & \le  C \left( \| \omega   \|^2_{H^2} + \| {\bf u}   \|^2_{L^2} \right)
\nonumber \\
& \le C \Big( \|\nabla {\bf u}  _t\|^2_{L^{2}} + \|\nabla \left( {\bf u}\cdot \nabla{\bf u}  + \mu \nabla\phi\right)   \|_{L^{2}}^{2}  + \| {\bf u}   \|^2_{L^2} + \| u_{1}\|_{H^{\frac{3}{2}}(\Gamma)}^{2} \Big)
\nonumber \\
& ~\,\,\,+C\left\|\left(\partial_{\bf n}\phi +\gamma_{fs}'(\phi)\right)\right\|^{2}_{H^{\frac{3}{2}}(\Gamma)} \left\|\partial_{x}\phi \right\|^{2}_{H^{\frac{3}{2}}(\Gamma)}
 \nonumber \\
 & \le C \Big( \|\nabla {\bf u}  _t\|^2_{L^{2}} + \| \nabla {\bf u} \|^2_{L^4} \| \nabla   {\bf u} \|^2_{L^4}  + \| {\bf u} \|^2_{L^\infty} \| \nabla^{2}{\bf u} \|^2_{L^2} + \| \nabla\mu   \|_{L^{2}}^{2} \| \nabla\phi   \|_{L^{\infty}}^{2} \Big)
\nonumber \\
& ~\,\,\, +  C \| \mu   \|_{L^{\infty}}^{2} \| \nabla^{2}\phi   \|_{L^{2}}^{2}   + C\| {\bf u}   \|^2_{H^2}  
\nonumber \\
& ~\,\,\, +C \left( \left\|\partial_{\bf n}\phi\right\|^{2}_{H^{\frac{3}{2}}(\Gamma)} +\left\|\gamma_{fs}'(\phi)\right\|^{2}_{H^{\frac{3}{2}}(\Gamma)}\right) \left\|\partial_{x}\phi \right\|^{2}_{H^{\frac{3}{2}}(\Gamma)}
 \nonumber \\
 & \le C \Big( \|\nabla {\bf u}  _t\|^2_{L^{2}} + \|  {\bf u} \|^4_{H^2} + \| \mu   \|_{H^{2}}^{2} \| \phi   \|_{H^{3}}^{2} + \| {\bf u}   \|^2_{H^2} +\| \phi \|_{H^{3}}^{4} \Big)
\nonumber \\
& \le  C \|\nabla {\bf u}  _t\|^2_{L^{2}} + C\mathcal{E}(t)^{2}.
\end{align}
Testing $(\ref{curl-u})_1$ by $-\Delta\omega$ and $ \omega $ respectively, integrating over $\Omega$ by parts, taking the sum of the result, and noting that $\nabla \times ({\bf u} \cdot \nabla{\bf u}) = {\bf u} \cdot \nabla \omega$,
\begin{align*}
&   \| \omega \|_{L^{2}} \le C  \| \nabla{\bf u}\|_{L^{2}}, ~~~~~~~~~~~~  \| \omega \|_{L^{2}(\Gamma)} \le  C\| \nabla{\bf u} \|_{L^{2}(\Gamma)} \le C\| {\bf u} \|_{H^{\frac32}},
\\
& \| \partial_{\bf n} \omega \|_{L^{2}(\Gamma)}  \le  C\| \partial_{\bf n} \nabla{\bf u} \|_{L^{2}(\Gamma)}  \le  C\| \partial_{\bf n} \nabla{\bf u} \|_{H^{\frac12}} \le C\| {\bf u} \|_{H^{\frac52}},
\end{align*}
and
\begin{align*}
& - \int_{\Omega}   {\bf u} \cdot \nabla \omega  \omega  {\rm d}x{\rm d}y = - \frac{1}{2} \int_{\Omega}   {\bf u} \cdot \nabla |\omega|^{2}  {\rm d}x{\rm d}y = \frac{1}{2} \int_{\Omega}  {\rm div} {\bf u}  |\omega|^{2}  {\rm d}x{\rm d}y  - \frac{1}{2} \int_{\Gamma}   {\bf u} \cdot {\bf n} |\omega|^{2}  {\rm d}x = 0 ,
\end{align*}
we have
\begin{align}\label{C-3}
& \dfrac{1}{2} \dfrac{\rm d}{{\rm d}t} \| \omega  \|^2_{H^1} + \| \nabla\omega  \|^2_{L^2} + \| \Delta\omega  \|^2_{L^2}
\nonumber \\
& = \int_{\Omega}   \left[ {\bf u} \cdot \nabla \omega + \nabla \times ( \mu \nabla\phi )  \right](-\omega+ \Delta\omega)  {\rm d}x{\rm d}y + \int_{\Gamma} \partial_{\bf n} \omega \left( \omega + \omega_{t} \right) {\rm d}x \nonumber \\
& = \int_{\Omega}   {\bf u} \cdot \nabla \omega  \Delta\omega  {\rm d}x{\rm d}y + \int_{\Omega}    \nabla \times ( \mu \nabla\phi )  (-\omega+ \Delta\omega)  {\rm d}x{\rm d}y + \int_{\Gamma} \partial_{\bf n} \omega \left( \omega + \omega_{t} \right) {\rm d}x
\nonumber \\
& \le \frac{1}{2} \| \Delta\omega  \|^2_{L^2}  +  \| {\bf u}  \|^2_{L^\infty} \| \nabla\omega  \|^2_{L^2}  + \|\nabla( \mu \nabla\phi )  \|^2_{L^2}  + \|\nabla( \mu \nabla\phi )  \|_{L^2} \| \omega  \|_{L^2}
\nonumber \\
&~\,\,\,+ \| \partial_{\bf n} \omega \|_{L^{2}(\Gamma)}  \left( \| \omega \|_{L^{2}(\Gamma)} + \| \omega_{t} \|_{L^{2}(\Gamma)}\right)
\nonumber \\
& \le \frac{1}{2} \| \Delta\omega  \|^2_{L^2}  +\| {\bf u}  \|^2_{H^2} \| \nabla^{2}  {\bf u} \|^2_{L^2}  + 2\left( \|\nabla \mu \|^2_{L^2} \| \nabla\phi  \|^2_{L^\infty} + \|\mu \|^2_{L^\infty} \| \nabla^{2}\phi  \|^2_{L^2}  \right) +  \| \nabla{\bf u}  \|_{L^2}^{2}
\nonumber \\
&~\,\,\,+ \| \partial_{\bf n} \omega \|_{L^{2}(\Gamma)}  \left( \| \omega \|_{L^{2}(\Gamma)} + \| \omega_{t} \|_{L^{2}(\Gamma)}\right)
\nonumber \\
& \le \frac{1}{2} \| \Delta\omega  \|^2_{L^2}   + C\| {\bf u}  \|^2_{H^2} \| \nabla^{2}  {\bf u} \|^2_{L^2}  + C\left( \|\nabla \mu \|^2_{L^2} \| \nabla\phi  \|^2_{H^2} + \|\mu \|^2_{H^2} \| \nabla^{2}\phi  \|^2_{L^2}  \right) + C \| \nabla{\bf u}  \|_{L^2}^{2}
\nonumber \\
&~\,\,\,+  C\| {\bf u} \|_{H^{\frac52}}\left( \| {\bf u} \|_{H^{\frac32}} + \left\|\partial_{t}u_{1} + \left(\partial_{\bf n}\phi_{t} + \phi_{t}\right)\partial_{x}\phi +  \left(\partial_{\bf n}\phi + \gamma_{fs}'(\phi) \right)\partial_{x}\phi_{t}\right\|_{L^{2}(\Gamma)}\right)
\nonumber \\
& \le \frac{1}{2} \| \Delta\omega  \|^2_{L^2}  +  C\mathcal{E}(t)^{2}  
\nonumber \\
& ~\,\,\,+  C\mathcal{E}(t)^{\frac12}
\| {\bf u} \|_{H^{2}}^{\frac12}\| {\bf u} \|_{H^{3}}^{\frac12}\left( 1  + \|  \partial_{t}u_{1}\|_{L^{2}(\Gamma)}  + \|\partial_{\bf n}\phi_{t} \|_{L^{2}(\Gamma)}  + \|\partial_{x}\phi_{t}\|_{L^{2}(\Gamma)}\right) \nonumber \\
& \le \frac{1}{2} \| \Delta\omega  \|^2_{L^2}  + \frac{\beta}{2} \|  \partial_{t}u_{1} \|_{L^{2}(\Gamma)}^{2} + \frac{1}{6} \|\partial_{\bf n}\phi_{t} \|_{L^{2}(\Gamma)}^{2} + \epsilon  \| {\bf u}  \|^2_{H^3}  + \frac{1}{4}\|\partial_{x}\phi_{t}\|_{L^{2}(\Gamma)}^{2} +  C\mathcal{E}(t)^{4},
\end{align}
where we have used the following fact
\begin{align*}
&  \omega_{t}   = \pm \left[\beta \partial_{t}u_{1}- \left(\partial_{\bf n}\phi_{t} +  \gamma_{fs}^{(2)}(\phi) \phi_{t}\right)\partial_{x}\phi -  \left(\partial_{\bf n}\phi + \gamma_{fs}'(\phi) \right)\partial_{x}\phi_{t}\right],  ~~{\rm on}~~ \{y=\pm 1 \}.
\end{align*}

\ \\[-2em]

{\bf Step 4. Closure of the estimates.}

Collecting $(\ref{C-1})$,  $(\ref{C-yb})$, $(\ref{C-2})$,  $(\ref{u-H33})$ and $(\ref{C-3})$ together, using Lemma $\ref{uuu}$ and  Korn's inequality, and choosing $\epsilon$ sufficiently small such that $C\epsilon<\frac{1}{12}$, we finally arrive at
\begin{align}\label{e-1}
 \sup\limits_{0\le t \le T} &\left( \| \nabla \phi_{t}\|_{L^{2}}^{2}  +  \|{\bf u}_t  \|^2_{L^2} + \| {\bf u}   \|^2_{H^2}  \right) + \int_{0}^{T} \mathcal{D}_{1}(t) {\rm d}t \nonumber \\
 & \hspace{8em} \le  C \left(\| {\bf u}_{0}\|_{H^{2}}, \| \phi_{0}\|_{H^{3}}, \| \mu_{0}\|_{H^{3}}\right) \int_{0}^{T} \mathcal{E}(t)^{4} {\rm d}t,
\end{align}
where $\mathcal{D}_{1}(t)$ is defined as in $(\ref{D-1})$.
Note that
\begin{align}\label{00}
\langle \phi_t   \rangle   & = \frac{1}{|\Omega|}\int_\Omega \phi_t   \mathrm{d}x{\rm d}y  = \frac{1}{|\Omega|}\int_\Omega \Delta\mu   \mathrm{d}x{\rm d}y  - \frac{1}{|\Omega|}\int_\Omega {\bf u} \cdot \nabla \phi   \mathrm{d}x {\rm d}y
\nonumber \\
& = - \frac{1}{|\Omega|}\int_{\Gamma} \partial_{\bf n} \mu {\rm d}x{\rm d}y  + \frac{1}{|\Omega|}\int_\Omega {\rm div} {\bf u} \,  \phi   \mathrm{d}x{\rm d}y -\frac{1}{|\Omega|}\int_{\Gamma} {\bf u}\cdot {\bf n} \phi   \mathrm{d}x = 0.
\end{align}
Hence,  using Poincar${\rm \acute{e}}$ inequality, we have
\begin{align}\label{dijie}
 \| \phi_{t}\|_{L^{2}}^{2}   = \left\| \phi_{t} -  \langle \phi_{t}   \rangle  \right\|_{L^{2}}^{2}  & \le C  \| \nabla\phi_{t}\|_{L^{2}}^{2} .
 \end{align}
From $(\ref{e-1})$, Lemma $\ref{energy}$ and $(\ref{dijie})$, the conclusion follows.
\end{proof}

Next, in order to get higher order estimates of $\phi$, we also need higher order estimates of $\mu  $.

\begin{lemma}
Under the assumptions of Theorem $\ref{TH}$, then for $0<T<1$,  there holds
\begin{align}\label{1-mu-2-phi-t}
&  \|\mu   \|^2_{H^2}  + \int_{0}^{T}  \left( \|\mu_{t}  \|^2_{L^2} +    \|\mu   \|^2_{H^3}  \right) \mathrm{d}t \le C \int_{0}^{T} \mathcal{E}(t)^{4} {\rm d}t + C,
\end{align}
where $C$ is a constant that depends on $\Omega$, $\beta$, $ \hat{\mathcal{E}}(0)$ and the initial value  $\left(\| {\bf u}_{0}\|_{H^{2}}, \| \phi_{0}\|_{H^{3}}, \| \mu_{0}\|_{H^{3}}\right)$.
\end{lemma}

\begin{proof}[\bf {Proof}]
First, from $(\ref{111111})_{4}$ and $(\ref{BC})_{4}$, we derive that
\begin{align}\label{p-mu}
\left| \langle \mu  _t    \rangle  \right| & = \left|   - \langle \Delta\phi_{t}   \rangle + \langle f_{t}  \rangle  \right| = \left| -  \frac{1}{|\Omega|}\int_{\Gamma} \partial_{\bf n}\phi_{t} {\rm d}x + \langle f_{t} \rangle  \right|
\nonumber \\
& \le C(|\Omega|,|\Gamma|) \Big( \|\partial_{\bf n}\phi_{t}  \|_{L^2(\Gamma)} + \|f_{t}\|_{L^2} \Big),
\end{align}
which implies
\begin{align}\label{mu-t-L2-L2}
\int_0^T  \|\mu  _{t}  \|^{2}_{L^2} \mathrm{d}t & \le C(\Omega) \int_0^T\left( \| \nabla \mu  _{t}   \|^{2}_{L^2} +  \left| \langle \mu  _t    \rangle  \right|^{2} \right) \mathrm{d}t \le C  \int_{0}^{T} \mathcal{E}(t)^{4} {\rm d}t + C.
\end{align}
The direct calculations yield that
\begin{align}
\dfrac{1}{2} \dfrac{\rm d}{{\rm d}t}   \| \mu     \|_{L^2}^2  & =  \int_{\Omega}  \mu\mu_{t}   \mathrm{d}x
\le  \| \mu  \|_{L^2}^2 +  \| \mu_{t}  \|^2_{L^2} \le C\mathcal{E}(t) + C  \int_{0}^{T} \mathcal{E}(t)^{4} {\rm d}t.
\end{align}
Integrating over $(0,T)$, we can deduce that
\begin{align}\label{mu-0}
\| \mu  \|_{L^2}^2
\le C (T) \int_{0}^{T} \mathcal{E}(t)^{4} {\rm d}t.
\end{align}
Finally, applying the elliptic estimates for $(\ref{111111})_{3}$ with boundary condition $(\ref{BC})_{2}$ gives that
\begin{align}\label{mu-H2}
\| \mu     \|^2_{H^2}  \le C \left( \| \mu     \|^2_{L^2} + \| \phi  _t \|^2_{L^2} + \| {\bf u} \|^2_{L^\infty}\| \nabla \phi   \|^2_{L^2} \right) \le C  \int_{0}^{T} \mathcal{E}(t)^{4} {\rm d}t + C,
\end{align}
and
\begin{align}\label{mu-H3}
\int_{0}^{T }\| \mu     \|^2_{H^3} {\rm d}t  & \le C \int_{0}^{T } \left( \| \mu  \|^2_{L^2} + \|  \phi  _t \|^2_{H^1}  + \| {\bf u} \|^2_{L^\infty}\| \nabla \phi   \|^2_{L^2} \right)  {\rm d}t
\nonumber \\
& ~\,\,\, + C \int_{0}^{T } \left( \| \nabla {\bf u} \|^2_{L^4}\| \nabla \phi   \|^2_{L^4} + \| {\bf u} \|^2_{L^\infty}\| \nabla^{2} \phi   \|^2_{L^2} \right) {\rm d}t
\nonumber \\
& \le C   \int_{0}^{T} \mathcal{E}(t)^{4} {\rm d}t + C,
\end{align}
where we have used $(\ref{mu-0})$, $(\ref{D-1})$. Combining  with $(\ref{mu-t-L2-L2})$ and $(\ref{mu-H2})$ leads $(\ref{1-mu-2-phi-t})$.
\end{proof}

\begin{lemma} \label{phi-H33}
Under the assumptions of Theorem $\ref{TH}$, then for $0<T<1$,  there holds
\begin{align}
& \| \nabla\phi \|^2_{H^2} + \| \nabla^{3}\mu\|_{L^2}^{2} + \int_{0}^{T} \left( \| \partial_{\bf n}\phi \|_{H^{2}(\Gamma)}^{2} +  \| \nabla^{4}\mu \|_{L^2}^{2} +  \| \nabla^{2}\phi_{t} \|_{L^2}^{2} \right)   {\rm d}t   \le        C  \int_{0}^{T} \mathcal{E}(t)^{8} {\rm d}t + C,
\end{align}
where $C$ is a constant that depends on $\Omega$, $\beta$,  $ \hat{\mathcal{E}}(0)$ and the initial value  $\left(\| {\bf u}_{0}\|_{H^{2}}, \| \phi_{0}\|_{H^{3}}, \| \mu_{0}\|_{H^{3}}\right)$.
\end{lemma}

\begin{proof}[\bf{Proof.}]
We complete the proof by several steps.

{\bf Step 1. Estimates for $\nabla \phi$ involving tangential derivatives.}

 For $k=1,2$, operating $\nabla\partial_{x}^{k}$, $\partial_{x}^{k}$ and $\partial_{x}^{k}$  to $(\ref{111111})_{3}$, $(\ref{111111})_{4}$ and $(\ref{N-BC})$ respectively, we get
\begin{align}\label{2-x}
\begin{cases}
\nabla\partial_{x}^{k}\Delta\mu  = \nabla\partial_{x}^{k}\phi_{t} + [\nabla\partial_{x}^{k}, {\bf u} \cdot \nabla] \phi+{\bf u} \cdot \nabla^{2} \partial_{x}^{k} \phi , & {\rm in}~~ \Omega,
\\
\partial_{x}^{k}\Delta\phi =  -\partial_{x}^{k}\mu  + \partial_{x}^{k}f, & {\rm in}~~\Omega,
\\
\partial_{x}^{k}\Delta\mu = - \partial_{x}^{k}\partial_{\bf n}\phi -  \partial_{x}^{k}\gamma_{fs}'(\phi) ,  & {\rm on}~~ \Gamma,
\end{cases}
\end{align}
where $[A,B] = AB-BA$ is commutators. We now consider
\begin{align}\label{E-phi-3}
-\langle \partial_{x}^{k}\Delta\mu,  \partial_{x}^{k}  \Delta\phi \rangle & =  \langle \nabla \partial_{x}^{k} \Delta\mu, \partial_{x}^{k} \nabla\phi \rangle  -  \int_{\Gamma} \partial_{x}^{k} \Delta\mu \partial_{x}^{k}\partial_{\bf n}\phi {\rm d}x
 =: J_{1} + J_{2}.
\end{align}
So as to obtain the information of $J_{1}$, putting $(\ref{2-x})_{1}$ into $J_{1}$, using integration by parts, and noting that
\begin{align*}
\int_{\Omega}  {\bf u} \cdot \nabla^{2} \partial_{x}^{k} \phi \cdot \nabla \partial_{x}^{k}  \phi  {\rm d}x{\rm d}y
 & =  \frac{1}{2}\int_{\Omega} {\bf u} \cdot \nabla |\nabla\partial_{x}^{k} \phi|^{2}  {\rm d}x{\rm d}y
\\
& = -\frac{1}{2}\int_{\Omega}  {\rm div}{\bf u}  |\nabla\partial_{x}^{k} \phi|^{2}  {\rm d}x{\rm d}y
+ \frac{1}{2}\int_{\Gamma} {\bf u} \cdot {\bf n} |\nabla\partial_{x}^{k} \phi|^{2}  {\rm d}x
=0,
\end{align*}
we have
\begin{align}\label{E-2-0}
J_{1} & = \int_{\Omega} \left(  \nabla\partial_{x}^{k}\phi_{t} + [\nabla\partial_{x}^{k}, {\bf u} \cdot \nabla] \phi  \right) \cdot \nabla\partial_{x}^{k}\phi  {\rm d}x{\rm d}y
\nonumber \\
 & =   \frac{1}{2} \frac{\rm d}{{\rm d}t} \| \nabla\partial_{x}^{k} \phi\|_{L^{2}}^{2} +\int_{\Omega} [\nabla\partial_{x}^{k}, {\bf u} \cdot \nabla] \phi  \cdot \nabla \partial_{x}^{k}\phi  {\rm d}x{\rm d}y.
\end{align}
We next put $(\ref{2-x})_{3}$ into $J_{2}$ and obtain
\begin{align}\label{E-2-1}
J_{2}  & =  \int_{\Gamma} \left( \partial_{x}^{k}\partial_{\bf n}\phi + \partial_{x}^{k}\gamma_{fs}'(\phi)  \right) \partial_{x}^{k} \partial_{\bf n}\phi {\rm d}x
 =  \| \partial_{x}^{k} \partial_{\bf n}\phi \|_{L^{2}(\Gamma)}^{2}  +  \int_{\Gamma}\partial_{x}^{k}\gamma_{fs}'(\phi)  \partial_{x}^{k}\partial_{\bf n}\phi  {\rm d}x.
\end{align}
On the other hand, it holds that
\begin{align}\label{E-2-2}
&-\langle \partial_{x}^{k} \Delta\mu, \partial_{x}^{k} \Delta\phi \rangle
\nonumber \\
& = \int_{\Omega}  \partial_{x}^{k} \Delta\mu   \left(  \partial_{x}^{k} \mu  - \partial_{x}^{k}f \right)  {\rm d}x{\rm d}y
=- \| \nabla \partial_{x}^{k} \mu \|_{L^{2}}^{2} + \int_{\Omega} \partial_{x}^{k}  \nabla\mu  \cdot \nabla \partial_{x}^{k} f {\rm d}x{\rm d}y,
\end{align}
where we have used $(\ref{2-x})_{2}$, integration by parts and $(\ref{BC})_{2}$.

Substituting $(\ref{E-2-0})$, $(\ref{E-2-1})$ and $(\ref{E-2-2})$ into $(\ref{E-phi-3})$, for $k=1,2$, we arrive at
\begin{align}\label{e-3}
& \frac{1}{2} \frac{\rm d}{{\rm d}t} \| \nabla\partial_{x}^{k} \phi\|_{L^{2}}^{2} +\| \nabla \partial_{x}^{k} \mu \|_{L^{2}}^{2} + \| \partial_{x}^{k} \partial_{\bf n}\phi \|_{L^{2}(\Gamma)}^{2}
\nonumber\\
& =-\int_{\Omega} [\nabla\partial_{x}^{k}, {\bf u} \cdot \nabla] \phi \cdot \nabla \partial_{x}^{k}\phi  {\rm d}x{\rm d}y
 +\int_{\Omega}  \partial_{x}^{k}\nabla\mu \cdot \nabla \partial_{x}^{k} f {\rm d}x{\rm d}y  -   \int_{\Gamma}   \partial_{x}^{k}\gamma_{fs}'(\phi)  \partial_{x}^{k}\partial_{\bf n}\phi  {\rm d}x
\nonumber \\
&=: \sum\limits_{i=1}^{3}M_{i}^{k}.
\end{align}
Next, for $k=1$, we estimate each term on the right hand side of $(\ref{e-3})$ as follows. Direct calculations show that
\begin{align*}
M_{1}^{1} & = - \int_{\Omega}  [\nabla\partial_{x}, {\bf u} \cdot \nabla] \phi \cdot\nabla \partial_{x}\phi  {\rm d}x{\rm d}y   = -\int_{\Omega}  \left( \nabla\partial_{x} ( {\bf u} \cdot \nabla\phi )- {\bf u} \cdot \nabla^{2}\partial_{x}\phi \right) \cdot \nabla \partial_{x}\phi  {\rm d}x{\rm d}y
 \\
 & = -\int_{\Omega}   \left( \nabla\partial_{x}  {\bf u} \cdot \nabla\phi + \partial_{x} {\bf u} \cdot \nabla^{2}\phi + \nabla  {\bf u} \cdot \nabla  \partial_{x}  \phi  \right) \cdot \nabla \partial_{x}\phi  {\rm d}x{\rm d}y
 \\
 & \le C \left( \| \nabla^{2}{\bf u} \|_{L^2} \| \nabla\phi \|_{L^{\infty}} + \| \nabla  {\bf u} \|_{L^4} \| \nabla^{2}\phi \|_{L^{4}} \right) \| \nabla^{2}\phi \|_{L^{2}}
  \\
 & \le C \left( \| \nabla^{2}{\bf u} \|_{L^2} \| \nabla\phi \|_{H^{2}} + \| \nabla  {\bf u} \|_{H^1} \| \nabla^{2}\phi \|_{H^{1}} \right) \| \nabla^{2}\phi \|_{L^{2}}
 \\
 & \le C \| {\bf u} \|_{H^2}^{2} + C\| \phi \|_{H^{3}} ^{4} \le C\mathcal{E}(t)^{2},
\\[1em]
M_{2}^{1} &= \int_{\Omega}   \nabla\partial_{x}  \mu \cdot \left[ (3\phi^{2}-1)\nabla\partial_{x}\phi  + 6 \phi \partial_{x} \phi\nabla\phi\right]  {\rm d}x{\rm d}y
 \\
 & \le \frac{1}{2} \| \nabla\partial_{x}  \mu \|_{L^2}^{2} + C \left(  \|   \phi \|_{L^\infty}^{4} + 1 \right)   \| \nabla  \partial_{x}\phi \|_{L^2}^{2} + C  \|   \phi \|_{L^\infty}^{2}   \|  \partial_{x} \phi \|_{L^4}^{2}    \|  \nabla\phi \|_{L^4}^{2}
  \\
 & \le \frac{1}{2} \| \nabla\partial_{x}  \mu \|_{L^2}^{2} + C \left(  \|   \phi \|_{H^2}^{4} + 1 \right)   \| \phi \|_{H^2}^{2} + C\|  \phi \|_{H^3}^{6} \le \frac{1}{2} \| \nabla\partial_{x}  \mu \|_{L^2}^{2} + C\mathcal{E}(t)^{3}.
 \end{align*}
It follows from Lemma $(\ref{Trace})$ that
\begin{align*}
M_{3}^{1} &= -\int_{\Gamma} \partial_{x} \gamma_{fs}'(\phi)  \partial_{x}\partial_{\bf n}\phi  {\rm d}x  = -\int_{\Gamma} \gamma_{fs}^{(2)}(\phi) \partial_{x} \phi   \partial_{x}\partial_{\bf n}\phi  {\rm d}x
\\
& \le  \frac{1}{2} \|   \partial_{x}\partial_{\bf n}\phi \|_{L^{2}(\Gamma)}^{2} +C \|  \partial_{x}\phi \|_{L^{2}(\Gamma)}^{2} \le   \frac{1}{2} \|  \partial_{x}\partial_{\bf n}\phi \|_{L^{2}(\Gamma)}^{2} + C \| \phi \|_{H^{2}}^{2}
\\
& \le    \frac{1}{2} \| \partial_{x}\partial_{\bf n}\phi \|_{L^{2}(\Gamma)}^{2}  + C\mathcal{E}(t).
 \end{align*}
And then, for $k=2$, we estimate each term on the right hand side of $(\ref{e-3})$ as follows. Direct calculations give that
\begin{align*}
M_{1}^{2} & = - \int_{\Omega}  [\nabla\partial_{x}^{2}, {\bf u} \cdot \nabla] \phi \cdot \nabla \partial_{x}^{2}\phi  {\rm d}x{\rm d}y  =  -\int_{\Omega} \left( \nabla\partial_{x}^{2} ( {\bf u} \cdot \nabla\phi )- {\bf u} \cdot \nabla^{2}\partial_{x}^{2}\phi \right) \cdot \nabla \partial_{x}^{2}\phi  {\rm d}x{\rm d}y
 \\
 & = - \int_{\Omega}  \left( \nabla\partial_{x}^{2}  {\bf u} \cdot \nabla\phi + 2\nabla\partial_{x}  {\bf u} \cdot \nabla \partial_{x}\phi + \partial_{x} ^{2}{\bf u}\cdot \nabla^{2}\phi   \right) \cdot \nabla \partial_{x}^{2}\phi  {\rm d}x{\rm d}y
 \\
 & ~\,\,\,+ \int_{\Omega}  \left(  \partial_{x} {\bf u}\cdot \nabla^{2}\partial_{x}\phi + \nabla  {\bf u} \cdot \nabla  \partial_{x}^{2}  \phi  \right) \cdot \nabla \partial_{x}^{2}\phi  {\rm d}x{\rm d}y
 \\
 & \le C \left( \| \nabla^{3}{\bf u} \|_{L^2} \| \nabla\phi \|_{L^{\infty}} + \| \nabla^{2}  {\bf u} \|_{L^4} \| \nabla^{2}\phi \|_{L^{4}}  +  \| \nabla{\bf u} \|_{L^\infty} \| \nabla^{3}\phi \|_{L^{2}} \right) \| \nabla^{3}\phi \|_{L^{2}}
  \\
 & \le C  \| \nabla {\bf u} \|_{H^2} \| \nabla\phi \|_{H^{2}}  \| \nabla^{3}\phi \|_{L^{2}}
  \le C \| {\bf u} \|_{H^3}^{2} + C\| \phi \|_{H^{3}} ^{4} \le  C \| {\bf u} \|_{H^3}^{2} + C\mathcal{E}(t)^{2},
\\[1em]
M_{2}^{2} &= \int_{\Omega}  \nabla\partial_{x}^{2} \mu \cdot \big[ (3\phi^{2}-1)\nabla\partial_{x}^{2}\phi  + 12 \phi \partial_{x} \phi \nabla\partial_{x}\phi + 6  \partial_{x} \phi \partial_{x} \phi \nabla\phi +  6 \phi \partial_{x}^{2} \phi\nabla\phi\big]  {\rm d}x{\rm d}y
 \\
 & \le \frac{1}{2} \| \nabla\partial_{x}^{2}  \mu \|_{L^2}^{2} + C \left(  \|   \phi \|_{L^\infty}^{4} + 1 \right)   \| \nabla  \partial_{x}^{2}\phi \|_{L^2}^{2} + C  \|   \phi \|_{L^\infty}^{2}   \|  \nabla \phi \|_{L^4}^{2}    \|  \nabla^2\phi \|_{L^4}^{2}
 \\
 &~\,\,\,+ C   \|   \nabla\phi \|_{L^\infty}^{4}\|   \nabla\phi \|_{L^2}^{2}
  \\
 & \le \frac{1}{2} \| \nabla\partial_{x}^{2}  \mu \|_{L^2}^{2} + C \left(  \|   \phi \|_{H^2}^{4} + 1 \right)   \| \phi \|_{H^3}^{2} + C \|  \phi \|_{H^3}^{6}
   \\
 & \le \frac{1}{2} \|\nabla\partial_{x}^{2}  \mu \|_{L^2}^{2} + C \mathcal{E}(t)^{3}.
 \end{align*}
Using Lemma $(\ref{Trace})$, one has
\begin{align*}
M_{3}^{2} &= -\int_{\Gamma} \partial_{x}^{2} \gamma_{fs}'(\phi)  \partial_{x}^{2}\partial_{\bf n}\phi  {\rm d}x
\\
& = - \int_{\Gamma}  \left(  \gamma_{fs}^{(2)}(\phi) \partial_{x}^{2} \phi  + \gamma_{fs}^{(2)}(\phi) \partial_{x} \phi\partial_{x} \phi \right)  \partial_{x}^{2}\partial_{\bf n}\phi  {\rm d}x
\\
& \le  \frac{1}{2} \|  \partial_{x}^{2}\partial_{\bf n}\phi \|_{L^{2}(\Gamma)}^{2} +C \|  \partial_{x}^{2}\phi \|_{L^{2}(\Gamma)}^{2}  + C \|  \partial_{x} \phi \|_{L^{4}(\Gamma)}^{4}
\\
& \le   \frac{1}{2} \|  \partial_{x}^{2}\partial_{\bf n}\phi \|_{L^{2}(\Gamma)}^{2} + C \| \phi \|_{H^{3}}^{4}
 \le    \frac{1}{2} \| \partial_{x}^{2}\partial_{\bf n}\phi \|_{L^{2}(\Gamma)}^{2}  + C\mathcal{E}(t)^{2}.
 \end{align*}
Finally, putting $M_{1}^{k}$-$M_{3}^{k}(k=1,2)$  into $(\ref{e-3})$ and combining with $(\ref{D-1})$, we obtain
\begin{align}\label{hz}
 &  \sup\limits_{0\le t \le T}\|  \nabla\partial_{x}^{k} \phi\|_{L^{2}}^{2}   + \int_{0}^{T} \left(  \|  \partial_{x}^{k} \nabla \mu \|_{L^{2}}^{2} +  \| \partial_{x}^{k} \partial_{\bf n}\phi \|_{L^{2}(\Gamma)}^{2}\right) {\rm d}t   \le C  \int_{0}^{T} \mathcal{E}(t)^{4} {\rm d}t + C, \hspace{0.2cm}{\rm for }~k=1,2.
 \end{align}

\ \\[-2em]

{\bf Step 2. Estimate for $\phi$ involving normal derivatives.}

 It follows from $(\ref{111111})_{4}$ that
\begin{align}\label{F-X}
 &\|  \partial_{y}^{2}\phi\|_{L^{2}}^{2} +  \|  \nabla\partial_{y}^{2} \phi\|_{L^{2}}^{2}
 \nonumber \\
 & \le  C \left( \| \partial_{x}^{2} \phi\|_{L^{2}}^{2}  + \|  \mu\|_{L^{2}}^{2} + \|  (\phi^{3}-\phi)\|_{L^{2}}^{2} \right)
 \nonumber \\
 &~\,\,\,+ C \left(  \|\nabla \partial_{x}^{2} \phi\|_{L^{2}}^{2}  + \| \nabla\mu\|_{L^{2}}^{2} + \| (3\phi^{2}-1)\nabla \phi\|_{L^{2}}^{2} \right)
 \nonumber \\
 &  \le  C \left( \| \nabla\partial_{x} \phi\|_{L^{2}}^{2} + \| \nabla\partial_{x}^{2} \phi\|_{L^{2}}^{2}\right)  + C\| \mu\|_{H^{1}}^{2} + C \| \phi \|_{H^{1}}^{3} + C + C\| \phi^{2}  \|_{L^4}^{2} \|\nabla \phi\|_{L^{4}}^{2}
  \nonumber \\
 &  \le    C \left( \| \nabla\partial_{x} \phi\|_{L^{2}}^{2} + \|  \nabla\partial_{x}^{2} \phi\|_{L^{2}}^{2}\right)  + C\| \mu\|_{H^{1}}^{2}   + C \| \phi \|_{H^{1}}^{3} + C + C\| \phi\|_{H^1}^{4}\|\nabla \phi\|_{L^{4}}^{2}
 \nonumber \\
 &\le \frac{1}{2} \| \partial_{y}^{2}\phi\|_{L^{2}}^{2} +   C  \int_{0}^{T} \mathcal{E}(t)^{4} {\rm d}t + C,
 \end{align}
where we have used  $(\ref{basic-E})$, $(\ref{D-1})$, $(\ref{1-mu-2-phi-t})$, $(\ref{hz})$ and the following fact
\begin{align*}
 \| \phi\|_{H^1}^{4} \|\nabla \phi\|_{L^{4}}^{2}  & \le   C\| \phi\|_{H^1}^{4}\|\nabla \phi\|_{L^{2}} \| \nabla \phi\|_{H^{1}}   \le C\|\phi\|_{H^{1}}^{10} + \epsilon\|\nabla^{2} \phi\|_{L^{2}}^{2}
 \\
 &\le  C\| \phi\|_{H^{1}}^{10} + C\epsilon \left( \nabla\partial_{x}\phi\|_{L^{2}}^{2}  +  \| \partial_{y}^{2} \phi\|_{L^{2}}^{2} \right).
 \end{align*}
Thus, we have
\begin{align}\label{N-B}
 \|  \nabla\phi\|_{H^2}^{2} &\le C \left(    \| \nabla\partial_{x} \phi\|_{L^{2}}^{2} +  \| \partial_{y}^{2}\phi\|_{L^{2}}^{2} \right) +  C \left(    \| \nabla\partial_{x}^{2} \phi\|_{L^{2}}^{2} +  \| \nabla \partial_{y}^{2} \phi\|_{L^{2}}^{2} \right)
 \nonumber \\
 & \le  C \int_{0}^{T} \mathcal{E}(t)^{4} {\rm d}t + C.
 \end{align}

{\bf Step 3. Estimates for $\Delta \mu$. }

Applying the elliptic estimates for $(\ref{111111})_{3}$ with boundary condition $(\ref{BC})_{2}$ and using $(\ref{D-1})$ and $(\ref{N-B})$, we have
\begin{align}
\| \nabla^{3}\mu\|_{L^2}^{2} &\le  C \left( \| \phi_{t}\|_{H^1}^{2}  + \| {\bf u}\cdot \nabla\phi  \|^2_{H^1} \right)
\nonumber \\
& \le C \left( \| \phi_{t}\|_{H^1}^{2}   +  \| \nabla{\bf u} \|^2_{L^2} \|  \nabla \phi  \|^2_{L^\infty}  +    \| {\bf u} \|^2_{L^\infty} \| \nabla^{2} \phi  \|^2_{L^2}  \right)
\nonumber \\
&\le  C \int_{0}^{T} \mathcal{E}(t)^{8} {\rm d}t + C.
\end{align}
Operating $\nabla\partial_{x}$, $\partial_{t}\partial_{x}$ and $\partial_{x}$  to $(\ref{111111})_{3}$, $(\ref{111111})_{4}$ and $(\ref{N-BC})$ respectively, we get
\begin{align}\label{2-x-t}
\begin{cases}
\nabla\partial_{x}\Delta\mu    =  \nabla\partial_{x}\phi_{t}  +  \nabla\partial_{x} {\bf u} \cdot  \nabla\phi + \partial_{x} {\bf u} \cdot \nabla^{2} \phi  +  \nabla {\bf u} \cdot \nabla \partial_{x} \phi + {\bf u} \cdot \nabla^{2} \partial_{x} \phi ,& {\rm in}~ \Omega,
\\
\partial_{x} \Delta\phi_{t} =  -\partial_{x} \mu_{t}  + (3\phi^{2}-1)\partial_{x} \phi_{t} +6\phi\phi_{t}\partial_{x} \phi, & {\rm in}~ \Omega,
\\
\partial_{x} \Delta\mu = - \partial_{x} \partial_{\bf n}\phi -   \gamma_{fs}^{(2)}(\phi) \partial_{x} \phi , &{\rm on}~ \Gamma.
\end{cases}
\end{align}
We now consider
\begin{align}\label{E-mu-4}
-\langle  \partial_{x} \Delta\mu,   \partial_{x}  \Delta\phi_{t}  \rangle & =   \langle \nabla \partial_{x}  \Delta\mu,  \partial_{x} \nabla\phi_{t}  \rangle-   \int_{\Gamma}   \partial_{x} \Delta\mu  \partial_{x} \partial_{\bf n}\phi_{t} {\rm d}x
 =: N_{1} + N_{2}.
\end{align}
Putting  $(\ref{2-x-t})_{1}$ into $N_{1}$, we know
\begin{align}\label{N-1}
 N_{1} & = \int_{\Omega}  \left(   \nabla\partial_{x}\phi_{t} + \nabla\partial_{x} {\bf u} \cdot  \nabla\phi + \partial_{x} {\bf u} \cdot \nabla \partial_{x} \phi +  \nabla {\bf u} \cdot \nabla \partial_{x} \phi  +  {\bf u} \cdot \nabla \partial_{x}^{2} \phi \right) \cdot \nabla\partial_{x}\phi_{t}  {\rm d}x{\rm d}y
 \nonumber \\
 & = \|\nabla\partial_{x}\phi_{t} \|_{L^{2}}^{2}  + \int_{\Omega} \left(    \nabla\partial_{x} {\bf u} \cdot  \nabla\phi + \partial_{x} {\bf u} \cdot \nabla \partial_{x} \phi  \right) \cdot \nabla\partial_{x}\phi_{t}  {\rm d}x{\rm d}y
 \nonumber \\
&~\,\,\, + \int_{\Omega}  \left(     \nabla {\bf u} \cdot \nabla \partial_{x} \phi  + {\bf u} \cdot \nabla \partial_{x}^{2} \phi \right) \cdot \nabla\partial_{x}\phi_{t} {\rm d}x{\rm d}y.
 \end{align}
We next put $(\ref{2-x-t})_{3}$ into $N_{2}$, use integration by parts and get
\begin{align}\label{N-3}
N_{2} & =   \int_{\Gamma}  \left( \partial_{x} \partial_{\bf n}\phi + \gamma_{fs}^{(2)}(\phi) \partial_{x} \phi \right) \partial_{x} \partial_{\bf n}\phi_{t} {\rm d}x
\nonumber \\
 & = \frac{1}{2}\frac{\rm d}{{\rm d}t} \| \partial_{x} \partial_{\bf n}\phi \|_{L^{2}(\Gamma)}^{2} + \int_{\Gamma}    \gamma_{fs}^{(2)}(\phi) \partial_{x} \phi \partial_{x} \partial_{\bf n}\phi_{t} {\rm d}x
 \nonumber \\
 & = \frac{1}{2}\frac{\rm d}{{\rm d}t} \|  \partial_{x} \partial_{\bf n}\phi \|_{L^{2}(\Gamma)}^{2} -  \int_{\Gamma} \left(  \gamma_{fs}^{(3)}(\phi) |\partial_{x} \phi|^{2} +\gamma_{fs}^{(2)}(\phi)\partial_x^2 \phi  \right)  \partial_{\bf n}\phi_{t} {\rm d}x.
\end{align}
On the other hand, it follows from integration by parts that
\begin{align}\label{E-mu}
-\langle \partial_{x} \Delta\mu,  \partial_{x}  \Delta\phi_{t} \rangle  & =  \int_{\Omega}  \partial_{x} \Delta\mu   \left(  \partial_{x} \mu_{t}  - (3\phi^{2}-1)\partial_{x} \phi_{t} -6\phi\phi_{t}\partial_{x} \phi \right)  {\rm d}x{\rm d}y
\nonumber \\
& =  - \frac{1}{2}\frac{\rm d}{{\rm d}t} \|  \partial_{x} \nabla \mu \|_{L^{2}}^{2} - \int_{\Omega}  \partial_{x} \Delta\mu  \left(    (3\phi^{2}-1)\partial_{x} \phi_{t} +6\phi\phi_{t}\partial_{x} \phi \right)  {\rm d}x{\rm d}y,
\end{align}
where we have used $(\ref{2-x-t})_{2}$ and $(\ref{BC})_{2}$.

Substituting $(\ref{N-1})$, $(\ref{N-3})$ and $(\ref{E-mu})$  into $(\ref{E-mu-4})$,  we have
\begin{align}\label{mu-high}
&\frac{1}{2}\frac{\rm d}{{\rm d}t} \left( \|  \partial_{x} \partial_{\bf n}\phi \|_{L^{2}(\Gamma)}^{2} + \| \partial_{x} \nabla \mu \|_{L^{2}}^{2} \right) + \| \nabla\partial_{x}\phi_{t} \|_{L^{2}}^{2}
\nonumber \\
&=  -\int_{\Omega} \left(    \nabla\partial_{x} {\bf u} \cdot  \nabla\phi + \partial_{x} {\bf u} \cdot \nabla \partial_{x} \phi + \nabla {\bf u} \cdot \nabla \partial_{x} \phi  + {\bf u} \cdot \nabla \partial_{x}^{2} \phi \right)\cdot\nabla\partial_{x}\phi_{t}  {\rm d}x{\rm d}y
\nonumber \\
& ~\,\,\,+  \int_{\Gamma} \left(   \gamma_{fs}^{(3)}(\phi) |\partial_{x} \phi|^{2} +  \gamma_{fs}^{(2)}(\phi)\partial_x^2 \phi  \right)  \partial_{\bf n}\phi_{t} {\rm d}x
\nonumber \\
& ~\,\,\,  - \int_{\Omega} \partial_{x} \Delta\mu  \left(    (3\phi^{2}-1)\partial_{x} \phi_{t} +6\phi\phi_{t}\partial_{x} \phi \right)  {\rm d}x{\rm d}y
\nonumber \\
& \le \frac{1}{2}  \| \nabla\partial_{x}\phi_{t} \|_{L^{2}}^{2} +  C\left(  \| \nabla^{2} {\bf u} \|_{L^{2}}^{2} \| \nabla \phi  \|_{L^{\infty}}^{2} + \| \nabla {\bf u} \|_{L^{4}}^{2} \| \nabla^{2} \phi  \|_{L^{4}}^{2} + \|   {\bf u} \|_{L^{\infty}}^{2} \| \nabla^{3} \phi  \|_{L^{2}}^{2}\right)
\nonumber \\
& ~\,\,\, + C\| \partial_{\bf n}\phi_{t} \|_{L^{2}(\Gamma)}^{2} + C\left(   \|  \partial_{x}\phi \|_{L^{4}(\Gamma)}^{4}   +  \| \partial_x^2\phi \|_{L^{2}(\Gamma)}^{2} \right)
\nonumber \\
& ~\,\,\, + C\| \nabla^{3}\mu\|_{L^{2}}^{2} + C \left( (\|\phi \|_{L^{\infty}}^{4} + 1) \|  \nabla\phi_{t}\|_{L^{2}}^{2} +\|\phi \|_{L^{\infty}}^{2} \| \phi_{t}\|_{L^{4}}^{2} \|  \nabla\phi \|_{L^{4}}^{2} \right)
\nonumber \\
& \le \frac{1}{2}  \| \nabla\partial_{x}\phi_{t}\|_{L^{2}}^{2} +  C\|  {\bf u} \|_{H^{2}}^{2} \|  \phi  \|_{H^{3}}^{2}  + C\| \partial_{\bf n}\phi_{t} \|_{L^{2}(\Gamma)}^{2} + C\|  \phi \|_{H^{3}}^{4} + C
\nonumber \\
& ~\,\,\, + C \| \mu \|_{H^{3}}^{2} + C \left[ (\| \phi \|_{H^{2}}^{4} + 1) \|  \nabla\phi_{t}\|_{L^{2}}^{2} + \|\phi \|_{H^{2}}^{2} \| \phi_{t}\|_{H^{1}}^{2}  \|  \nabla\phi \|_{H^{1}}^{2} \right]
\nonumber \\
& \le \frac{1}{2}  \|\nabla\partial_{x}\phi_{t} \|_{L^{2}}^{2} +  C\mathcal{E}(t)^{3} + C\| \partial_{\bf n}\phi_{t} \|_{L^{2}(\Gamma)}^{2} .
\end{align}
Using Gronwall inequality and $(\ref{D-1})$, we find
\begin{align}\label{mu-high-F}
\sup\limits_{0\le t\le T} \left( \|  \partial_{x} \partial_{\bf n}\phi \|_{L^{2}(\Gamma)}^{2} + \| \partial_{x} \nabla \mu \|_{L^{2}}^{2} \right) +   \int_{0}^{T}  \| \nabla\partial_{x}\phi_{t} \|_{L^{2}}^{2} {\rm d}t \le C  \int_{0}^{T} \mathcal{E}(t)^{4} {\rm d}t  + C.
\end{align}
Due to $(\ref{111111})_{4}$ and $(\ref{basic-E})$, one has
\begin{align}
\int_{0}^{T}  \|  \partial_{y}^{2} \phi_{t}\|_{L^{2}}^{2}   {\rm d}t  &\le  C \int_{0}^{T} \left( \|\partial_x^2 \phi_{t}\|_{L^{2}}^{2}  + \| \mu_{t}\|_{L^{2}}^{2} + \| (3\phi^{2}-1)\phi_{t}\|_{L^{2}}^{2} \right)   {\rm d}t
\nonumber \\
&\le  C \int_{0}^{T} \left( \| \nabla\partial_{x} \phi_{t}\|_{L^{2}}^{2}  + \|  \mu_{t}\|_{L^{2}}^{2} + (\| \phi\|_{L^{8}}^{4}+1) \|\phi_{t}\|_{L^{4}}^{2} \right)   {\rm d}t
\nonumber \\
&\le  C \int_{0}^{T} \left( \| \nabla\partial_{x} \phi_{t}\|_{L^{2}}^{2}  + \|  \mu_{t}\|_{L^{2}}^{2} + (\| \phi\|_{H^{1}}^{4}+1) \|\phi_{t}\|_{H^{1}}^{2} \right)   {\rm d}t
\nonumber \\
&\le  C \int_{0}^{T} \left( \| \nabla\partial_{x} \phi_{t}\|_{L^{2}}^{2} + \|  \mu_{t}\|_{L^{2}}^{2}\right)   {\rm d}t  + C \sup\limits_{0\le t\le T} \left( \| \phi\|_{H^{1}}^{4}+1 \right) \int_{0}^{T}   \|\phi_{t}\|_{H^{1}}^{2}    {\rm d}t
\nonumber \\
 & \le  C  \int_{0}^{T} \mathcal{E}(t)^{4} {\rm d}t  + C,
 \end{align}
which together with $(\ref{mu-high-F})$ implies that
\begin{align}\label{phi-t-2}
\int_{0}^{T}  \|  \nabla^{2} \phi_{t}\|_{L^{2}}^{2}   {\rm d}t  \le  C  \int_{0}^{T} \mathcal{E}(t)^{4} {\rm d}t  + C.
 \end{align}
Finally, it follows from $(\ref{111111})_{3}$ that
\begin{align}\label{mu-40}
\int_{0}^{T} \| \nabla^{2}\Delta \mu \|_{L^2}^{2}  {\rm d}t \le  C \int_{0}^{T} \mathcal{E}(t)^{4} {\rm d}t + C.
\end{align}
Combining $(\ref{hz})$,  $(\ref{mu-high-F})$, $(\ref{phi-t-2})$ and $(\ref{mu-40})$, one gets the desired estimates.
\end{proof}



\section{Construction of solutions to nonlinear NSCH problem}\label{sec:5}

Our goal of this section is to construct solution to the nonlinear NSCH problem $(\ref{111111})$, $(\ref{BC})$, $(\ref{initial})$ on a time interval $[0,T]$ for some $T>0$.

\subsection{The $\delta$-approximate problem}\label{subsec:4}

For any fixed $0<\delta < 1$, we consider $\delta$-approximate problem in channel domain $\Omega=\mathbb{T}\times(-1,1)$ as follows
\begin{align}\label{A-NSCH}
\begin{cases}
{\bf u}_t+({\bf u}\cdot\nabla){\bf u}+ \nabla p = {\rm div}\mathbb{S}({\bf u})+\mu\nabla \phi, &~~\text{in}~~\Omega,
\\
{\rm div}{\bf u}=0, &~~\text{in}~~\Omega,
\\
\phi_t - \delta \partial_{x}^{2} \phi +({\bf u}\cdot\nabla)\phi=\Delta\mu, &~~\text{in}~~\Omega,
\\
\mu=\delta\phi_{t}-\Delta\phi+f(\phi), &~~\text{in}~~\Omega,
\end{cases}
\end{align}
with the boundary condition
\begin{align}\label{A-BC}
\begin{cases}
u_{2}= 0, &~~\text{on}~~\Gamma,
\\
\partial_{\bf n}\mu =0, &~~\text{on}~~\Gamma,
\\
\beta u_1+ \partial_{\bf n} u_1= \big( \partial_{\bf n}\phi+\gamma_{fs}'(\phi)\big) \partial_{x} \phi, &~~\text{on}~~\Gamma,
\\
\phi_t - \delta  \partial_{x}^{2} \phi+u_1  \partial_{x} \phi = -\partial_{\bf n}\phi-\gamma_{fs}'(\phi), &~~\text{on}~~\Gamma,
\end{cases}
\end{align}
and the initial condition
\begin{equation} \label{A-IC}
({\bf u},\phi)\big|_{t=0}=({\bf u}_0,\phi_0),  ~~\text{in}~~\Omega.
\end{equation}

\begin{remark}\label{R-a}
By adding this artificial viscosity term $ -\delta \partial_{x}^{2}\phi $ in \eqref{A-BC}$_{4}$, inspired by \cite{C-W-X, G-1, M-Z}, we can obtain from \eqref{A-BC}$_{4}$ the high regularity of ${\bf tr}(\phi)$, which is used for fixed point argument. 
\end{remark}

\begin{remark}\label{R-add}
The aim we add the term $-\delta \partial_{x}^{2}\phi$ in \eqref{A-NSCH}$_{3}$ is to obtain the uniform-in-$\delta$ estimates since, to obtain these estimates,  we need this term to match the boundary condition \eqref{A-BC}$_{4}$. Otherwise,  \eqref{A-NSCH}$_{3}$ is
$$
\phi_t +({\bf u}\cdot\nabla)\phi=\Delta\mu,  ~\,\,\text{in}~~\Omega,
$$
which leads to, by taking trace, that
\begin{align}\label{d-sm}
\Delta \mu = \phi_t + u_1  \partial_{x} \phi  , ~\,\,\text{on}~ \Gamma.
\end{align}
It follows from \eqref{d-sm} and \eqref{A-BC}$_{4}$ that 
\begin{align}\label{ddd}
\Delta \mu =   \delta  \partial_{x}^{2} \phi -\partial_{\bf n}\phi-\gamma_{fs}'(\phi), ~\,\,\text{on}~ \Gamma.
\end{align}
Compared with \eqref{N-BC},  there is a high order term  $\delta  \partial_{x}^{2} \phi $ in the right hand side of \eqref{ddd}. Therefore, it is difficult for us to obtain $\delta$-independent estimates.
\end{remark}

\begin{remark}\label{R-O}
For general domain,  the term $-\delta\partial_{x}^{2}\phi$ in \eqref{A-BC}$_{4}$ should be replaced by $-\delta \Delta_{\boldsymbol \tau} \phi $, and the corresponding term $-\delta\partial_{x}^{2}\phi$ in \eqref{A-NSCH}$_{3}$ should be replaced by $-\delta \Delta_{\boldsymbol \tau} \phi $ too. The operator $ \Delta_{\boldsymbol \tau} $ is well-defined in the neighborhood near the boundary. However,  it is not well-defined in whole domain. Thus, we only consider channel domain here.
\end{remark}

It is worth noting that we do not  impose boundary conditions for $\phi$ and thus there is no boundary layer appearing as $\delta\to 0$. Formally, \eqref{A-NSCH}$_{4}$ is the parabolic equation of $\phi$ with the source term $\mu-f(\phi)$. However, we can not obtain the uniqueness of $\mu$ from the elliptic equation \eqref{A-NSCH}$_{3}$ with Neumann boundary condition \eqref{A-BC}$_{2}$. Thus, we now put $\phi_{t}$ of $(\ref{A-NSCH})_{3}$ into $(\ref{A-NSCH})_{4}$, then rewrite  \eqref{A-NSCH} as follows
\begin{align}\label{RA-NSCH}
\begin{cases}
{\bf u}_t+({\bf u}\cdot\nabla){\bf u}+ \nabla p = {\rm div}\mathbb{S}({\bf u})+\mu\nabla \phi, &~~\text{in}~~\Omega,
\\
{\rm div}{\bf u}=0, &~~\text{in}~~\Omega,
\\
\mu - \delta\Delta\mu = -\Delta\phi+f(\phi) +  \delta^{2}   \partial_{x}^{2} \phi- \delta {\bf u}\cdot \nabla \phi, &~~\text{in}~~\Omega,
\\
\delta\phi_{t}-\Delta\phi = \mu - f(\phi), &~~\text{in}~~\Omega,
\end{cases}
\end{align}

Note that the problem \eqref{A-NSCH}-\eqref{A-IC} and the problem \eqref{RA-NSCH}, \eqref{A-BC}, \eqref{A-IC} are equivalent.  Next, our goal  is to derive the solvability of the problem \eqref{RA-NSCH}, \eqref{A-BC}, \eqref{A-IC}.

We denote $\mu_0^\delta=\mu(x,y,0)$ as the initial value of $\mu$ for the $\delta$-approximate problem. Indeed, it follows from $(\ref{RA-NSCH})_3$ and $(\ref{A-BC})_2$ that $\mu_0^\delta$ is the unique solution to the following elliptic problem
\begin{align}\label{mu-000}
\begin{cases}
\mu^\delta_0-\delta \Delta \mu^\delta_0  = - \delta{\bf u}_0\cdot\nabla\phi_0  + \delta^{2}\partial_x^2 \phi_{0}+ \mu_0, ~&\text{in}~\Omega,
\\
\partial_{\bf n}\mu^\delta_0\big|_{\Gamma}=0, ~&\text{on}~\Gamma.
\end{cases}
\end{align}

The construction of the solution for the $\delta$-approximate problem relies on linearized problems with fixed point argument. The key point is the solvability of each toy model and the uniform $\delta$-independent estimates of the solutions, which allows us to pass to the limits as $\delta\to 0$ to produce the solution to the nonlinear NSCH
problem $(\ref{111111})$, $(\ref{BC})$, \eqref{initial}.

\subsection{Analysis for linearized toy models} The first linear problem is the following linear elliptic problem of $\mu$ with given $G_1$:
\begin{align}\label{L-mu}
\begin{cases}
\mu - \delta\Delta\mu = G_1, &~~\text{in}~~\Omega,
\\
\partial_{\bf n}\mu =0, &~~\text{on}~~\Gamma.
\end{cases}
\end{align}

\begin{proposition}\label{P-mu}
Given $G_1 \in L^{\infty}(0,T;H^{2})$, then the problem $(\ref{L-mu})$ admits a unique solution $\mu$ that satisfies
\begin{align}\label{L-mu-E}
  \| \mu \|_{H^{3}}^{2}  \le  C(\delta^{-1})   \| G_1 \|_{H^{1}}^{2}.
 \end{align}
Moreover, if $G_1 = \left(-\Delta +\partial_x^2\right)G_{11} + G_{12}$, where $ \nabla \partial_{t}G_{11} \in L^{\infty}(0,T; L^{2})$ and $ \partial_{t}G_{12} \in L^{\infty}(0,T;$
$L^{2})$, then
\begin{align}\label{L-mu-E2}
\| \mu_{t} \|_{H^{1}}^{2}  \le  C(\delta^{-1}) \| \nabla  \partial_{t} G_{11}  \|_{L^{2}}^{2} + C \| \partial_{t} G_{12}  \|_{L^{2}}^{2}.
\end{align}
\end{proposition}

\begin{proof}[\bf{Proof}]
The existence, uniqueness and $H^{3}$-regularity of the problem $(\ref{L-mu})$ is classical. We shall show that $(\ref{L-mu-E2})$ holds.

Differentiating $(\ref{L-mu})_{1}$ with respect to $t$, and then testing it by $\mu_{t}$, we find
\begin{align*}
&\| \mu_{t} \|_{L^{2}}^{2} + \delta \| \nabla\mu_{t} \|_{L^{2}}^{2}
\nonumber \\
& = \int_{\Omega} \left(-\Delta +\partial_x^2\right) \partial_{t}G_{11} \mu_{t}   {\rm d}x{\rm d}y + \int_{\Omega}  \partial_{t} G_{12} \mu_{t}   {\rm d}x{\rm d}y
\nonumber \\
& = \int_{\Omega} \nabla  \partial_{t} G_{11}  \cdot \nabla \mu_{t}   {\rm d}x{\rm d}y - \int_{\Omega} \partial_{x}  \partial_{t} G_{11}  \cdot \partial_{x} \mu_{t}   {\rm d}x{\rm d}y + \int_{\Omega}  \partial_{t} G_{12} \mu_{t}   {\rm d}x{\rm d}y
-  \int_{\Gamma} \partial_{\bf n}  \partial_{t} G_{11}  \mu_{t}   {\rm d}x
\nonumber \\
& \le \frac{\delta}{4} \| \mu_{t} \|_{H^{1}}^{2} + C(\delta^{-1}) \| \nabla  \partial_{t} G_{11}  \|_{L^{2}}^{2} + C \| \partial_{t} G_{12}  \|_{L^{2}}^{2}
+ C \| \partial_{\bf n}  \partial_{t} G_{11}  \|_{H^{-\frac{1}{2}}(\Gamma)} \|   \mu_{t}  \|_{H^{\frac{1}{2}}(\Gamma)}
\nonumber \\
& \le \frac{\delta}{4} \| \mu_{t} \|_{H^{1}}^{2}  + C(\delta^{-1}) \| \nabla  \partial_{t} G_{11}  \|_{L^{2}}^{2} + C \| \partial_{t} G_{12}  \|_{L^{2}}^{2}
+ C \| \partial_{\bf n}  \partial_{t} G_{11}  \|_{L^{2}} \|   \mu_{t}  \|_{H^{1}}
\nonumber \\
& \le \frac{\delta}{2} \| \mu_{t} \|_{H^{1}}^{2}  + C(\delta^{-1}) \| \nabla  \partial_{t} G_{11}  \|_{L^{2}}^{2} + C \| \partial_{t} G_{12}  \|_{L^{2}}^{2}.
\end{align*}
This completes the proof of Proposition $\ref{P-mu}$.
\end{proof}

The second linear problem is the linear parabolic problem of $\psi$ with given $G_2$:
\begin{align}\label{L-psi}
\begin{cases}
\psi_{t} -\delta \partial_{x}^{2} \psi = G_2, &~~\text{on}~\Gamma,
\\
\psi \big|_{t=0} = \psi_{0} := \phi_{0} \big|_{\Gamma}  , & ~~\text{on}~~\Gamma.
\end{cases}
\end{align}

\begin{proposition}
Given $G_2 \in L^{2}(0,T;H^{2}(\Gamma))$ and $ \partial_{t}G_2 \in L^{2}(0,T;L^{2}(\Gamma))$, and assume $\psi_{0}\in H^{3}(\Gamma)$, then the problem $(\ref{L-psi})$ admits a unique solution $\mu$ that satisfies
\begin{align}\label{L-psi-E}
&\sup\limits_{0\le t \le T} \left(  \| \psi_{t} \|_{H^{1}(\Gamma)}^{2} +  \| \psi \|_{H^{3}(\Gamma)}^{2} \right)
+ \int_{0}^{T}\left( \| \psi_{tt} \|_{L^{2}(\Gamma)}^{2} + \| \psi_{t} \|_{H^{2}(\Gamma)}^{2} +  \| \psi \|_{H^{4}(\Gamma)}^{2} \right) {\rm d}t
\nonumber \\
&\le  C(\delta^{-1})\left[ \| \psi_{0} \|_{H^{3}(\Gamma)}^{2} +  \| G_{2}(0) \|_{H^{1}(\Gamma)}^{2}  +   \int_{0}^{T} \mathcal{D}_{2}(t) {\rm d}t\right],
\end{align}
where
\begin{align*}
\mathcal{D}_{2}(t) := \| \partial_{t} G_2 \|_{L^{2}(\Gamma)}^{2}  +  \| G_2 \|_{H^{2}(\Gamma)}^{2} .
\end{align*}
\end{proposition}

\begin{proof}[\bf{Proof}]
The existence, uniqueness and regularity of the problem $(\ref{L-psi})$ can be obtained via the classical theory for parabolic equations and we omit the details here.
\end{proof}

The third linear problem is the linear parabolic problem of $\phi$ with given $G_{3}$ and $g_{3}$:
\begin{align}\label{L-phi}
\begin{cases}
\delta\phi_{t} -\Delta \phi = G_{3}, &~~\text{in}~~\Omega,
\\
\phi =g_{3}, &~~\text{on}~~\Gamma,
\\
\phi\big|_{t=0} = \phi_{0}, & ~~\text{in}~~\Omega.
\end{cases}
\end{align}

\begin{proposition}\label{P-phi}
Given $G_{3} \in  L^{\infty}(0,T;H^{1})\cap L^{2} (0,T;H^{2}) $, $\partial_{t}G_{3} \in  L^{2}(0,T;L^{2}) $,
 $g_{3} \in L^{2}(0,T;$ $H^{\frac{7}{2}}(\Gamma))$ $\partial_{t}g_{3} \in L^{2}(0,T;H^{\frac{3}{2}}(\Gamma)) $ and $ \partial_{t}^{2} g_{3} \in L^{2}(0,T;H^{-\frac{1}{2}}(\Gamma)) $, and assume that $\phi_{0} \in H^{3}$ and $ G_{0}^{3} \in H^{1}$. Then
\begin{align}\label{L-phi-E}
& \sup\limits_{0\le t \le T} \left( \| \phi \|_{H^{3}}^{2} + \| \phi_{t} \|_{H^{1}}^{2} \right) + \int_{0}^{T} \left( \| \phi_{t} \|_{H^{2}}^{2} + \| \phi_{tt} \|_{L^{2}}^{2} + \| \phi \|_{H^{4}}^{2} \right)  {\rm d}t
\nonumber \\
& \le  C (\delta^{-1}) \bigg[ \| \phi_{0}\|_{H^{3}}^{2} +  \|  G_{3}(0)\|_{H^{1}}^{2} +  \| G_{3} \|_{H^{1}}^{2} + \int_{0}^{T}  \mathcal{D}_{3}(t)  {\rm d}t \bigg].
\end{align}
where
\begin{align*}
\mathcal{D}_{3}(t) := \| \partial_{t}G_{3}\|_{L^{2}}^{2} + \| G_{3}\|_{H^{2}}^{2} +  \| \partial_{t}^{2} g_{3} \|_{H^{-\frac12}(\Gamma)}^{2}  + \| \partial_{t}  g_{3} \|_{H^{\frac32}(\Gamma)}^{2} +  \|   g_{3} \|_{H^{\frac72}(\Gamma)}^{2}.
\end{align*}
\end{proposition}

\begin{proof}[\bf{Proof}]
The existence and uniqueness of the problem $(\ref{L-phi})$ is classical, we next consider the regularity of solution to  the problem $(\ref{L-phi})$.

Operating $ \partial_{t} $ to $(\ref{L-phi})_{1}$,  we get
\begin{align*}
\delta \phi_{tt} - \Delta \phi_{t} = \partial_{t}G_{3}  .
\end{align*}
Testing it by $\phi_{t} - \Delta \phi_{t}$ and integrating over $\Omega$ by parts, we find
\begin{align*}
&\frac{\delta}{2}\frac{\rm d}{{\rm d}t} \| \phi_{t}\|_{H^{1}}^{2}  +  \| \nabla\phi_{t}\|_{L^{2}}^{2} +  \| \Delta \phi_{t}\|_{L^{2}}^{2}
\\
& =    \int_{\Omega}  \partial_{t}G_{3}  \phi_{t} {\rm d}x{\rm d}y - \int_{\Omega}  \partial_{t}G_{3}  \Delta \phi_{t} {\rm d}x{\rm d}y  +  \int_{\Gamma}  \partial_{\bf n} \phi_{t}  \partial_{t} g_{3}{\rm d}x + \delta \int_{\Gamma}  \partial_{\bf n} \phi_{t}  \partial_{t}^{2} g_{3}{\rm d}x
\\
& \le \frac{1}{4} \| \Delta \phi_{t}\|_{L^{2}}^{2} + \| \phi_{t}\|_{H^{1}}^{2} + C \| \partial_{t}G_{3}\|_{L^{2}}^{2} + \epsilon \delta \| \partial_{\bf n} \phi_{t} \|_{H^{\frac12}(\Gamma)}^{2} + C\| \partial_{t} g_{3} \|_{H^{-\frac12}(\Gamma)}^{2}  + C \| \partial_{t}^{2} g_{3} \|_{H^{-\frac12}(\Gamma)}^{2}
\\
& \le \frac{1}{2} \| \Delta \phi_{t}\|_{L^{2}}^{2} + C \| \phi_{t}\|_{
H^{1}}^{2}  + C \| \partial_{t}g_{3}\|_{L^{2}}^{2} + C\| \partial_{t} g_{3} \|_{H^{\frac32}(\Gamma)}^{2} + C \| \partial_{t}^{2} g_{3} \|_{H^{-\frac12}(\Gamma)}^{2},
\end{align*}
which together with Gronwall inequality implies that
\begin{align*}
&  \sup\limits_{0\le t \le T}  \delta \| \phi_{t}\|_{H^{1}}^{2}  +  \int_{0}^{T} \left(\| \nabla\phi_{t}\|_{L^{2}}^{2} + \| \Delta \phi_{t}\|_{L^{2}}^{2} \right) {\rm d}t
\nonumber\\
&\le  C\left[   \|  \phi_0\|_{H^{3}}^{2} + \|  G_{3}(0)\|_{H^{1}}^{2} +    \int_{0}^{T} \left( \| \partial_{t}G_{3}\|_{L^{2}}^{2} + \| \partial_{t} g_{3} \|_{H^{\frac32}(\Gamma)}^{2}+  \| \partial_{t}^{2} g_{3} \|_{H^{-\frac12}(\Gamma)}^{2} \right) {\rm d}t \right].
\end{align*}

Testing $(\ref{L-phi})_{1}$ by $\phi$ and integrating over $\Omega$ by parts, for $\epsilon$ small enough, we find 
\begin{align*}
\frac{\delta}{2}\frac{\rm d}{{\rm d}t} \| \phi\|_{L^{2}}^{2}  +  \| \nabla\phi\|_{L^{2}}^{2} & =    \int_{\Omega}  G_3  \phi {\rm d}x{\rm d}y +  \int_{\Gamma}  \partial_{\bf n} \phi   g_3{\rm d}x
\\
& \le \delta \|  \phi \|_{L^{2}}^{2} + C(\delta^{-1}) \| G_{3}\|_{L^{2}}^{2} + \epsilon  \| \partial_{\bf n} \phi \|_{H^{-\frac12}(\Gamma)}^{2} + C\| g_{3}\|_{H^{\frac12}(\Gamma)}^{2}
\\
& \le  \delta \|  \phi \|_{L^{2}}^{2} + C(\delta^{-1}) \| G_{3}\|_{L^{2}}^{2} + \frac{1}{2} \| \nabla \phi \|_{L^{2}}^{2} + C\| g_{3} \|_{H^{\frac12}(\Gamma)}^{2},
\end{align*}
which implies
\begin{align*}
 \sup\limits_{0\le t \le T} \| \phi\|_{L^{2}}^{2}  + \int_{0}^{T} \| \nabla\phi\|_{L^{2}}^{2} {\rm d}t  \le  C(\delta^{-1},T) \left[  \| \phi_{0}\|_{L^{2}}^{2} + \int_{0}^{T} \left( \| G_{3}\|_{L^{2}}^{2} + \| g_{3} \|_{H^{\frac12}(\Gamma)}^{2} \right) {\rm d}t \right].
\end{align*}

Operating $  \partial_{x}^{2} $ to $(\ref{L-phi})_{1}$,  we get
\begin{align*}
\delta   \partial_{x}^{2}  \phi_{t} - \Delta  \partial_{x}^{2} \phi  = \partial_{x}^{2} G_{3}  .
\end{align*}
Testing it by $ \partial_{x}^{2} \phi_{t} $ and integrating over $\Omega$ by parts, we find
\begin{align}\label{A-phi-3}
& \frac{1}{2}\frac{\rm d}{{\rm d}t} \| \nabla \partial_{x}^{2} \phi\|_{L^{2}}^{2} +  \delta \| \partial_{x}^{2} \phi_{t}\|_{L^{2}}^{2}
\nonumber\\
& =  \int_{\Omega}\partial_{x}^{2} G_{3}\partial_{x}^{2} \phi_{t} {\rm d}x{\rm d}y  + \int_{\Gamma}  \partial_{\bf n} \partial_{x}^{2} \phi \,  \partial_{t} \partial_{x}^{2}  g_{3}{\rm d}x
\nonumber\\
& \le \| \partial_{x}^{2} \phi_{t}\|_{L^{2}}^{2} +  \| G_{3}\|_{H^{2}}^{2} + C  \| \partial_{\bf n} \partial_{x}^{2}\phi \|_{H^{\frac12}(\Gamma)}^{2} + C \| \partial_{t} \partial_{x}^{2} g_{3} \|_{H^{-\frac12}(\Gamma)}^{2}
\nonumber\\
& \le  \|  \phi_{t}\|_{H^{2}}^{2}  +  \| G_{3}\|_{H^{2}}^{2} + C \| \phi \|_{H^{4}}^{2}  + C \| \partial_{t}  g_{3} \|_{H^{\frac32}(\Gamma)}^{2}.
\end{align}
Applying elliptic regularity theory for $(\ref{L-phi})$, we have
\begin{align}\label{A-phi-4}
 \| \phi_{t} \|_{H^{2}}^{2}   +   \| \phi \|_{H^{4}}^{2}  & \le  C  \left( \| \phi_{t} \|_{H^{2}}^{2}  + \| G_{3}\|_{H^{2}}^{2}  +  \|   g_{3} \|_{H^{\frac72}(\Gamma)}^{2}\right)
\nonumber\\
&\le  C  \left( \| \Delta \phi_{t} \|_{L^{2}}^{2}  + \| G_{3}\|_{H^{2}}^{2} + \| \partial_{t}  g_{3} \|_{H^{\frac32}(\Gamma)}^{2} +  \|   g_{3} \|_{H^{\frac72}(\Gamma)}^{2}\right) .
\end{align}
Putting $(\ref{A-phi-4})$ into $(\ref{A-phi-3})$ and using Gronwall inequality, we infer
\begin{align}\label{A-E1}
\sup\limits_{0\le t \le T} \| \nabla \partial_{x}^{2} \phi\|_{L^{2}}^{2}  &\le   C\left[   \|  \phi_0\|_{H^{3}}^{2} + \|  G_{3}(0)\|_{H^{1}}^{2} +    \int_{0}^{T} \left( \| \partial_{t}G_{3}\|_{L^{2}}^{2} +  \| \partial_{t}^{2} g_{3} \|_{H^{-\frac12}(\Gamma)}^{2}\right) {\rm d}t \right]
 \nonumber\\
 &~\,\,\,+ C\int_{0}^{T} \left(\| G_{3}\|_{H^{2}}^{2} + \| \partial_{t}  g_{3} \|_{H^{\frac32}(\Gamma)}^{2} +  \|   g_{3} \|_{H^{\frac72}(\Gamma)}^{2} \right) {\rm d}t .
\end{align}
It follows from $(\ref{L-phi})_{1}$ that
\begin{align}\label{A-E2}
 \| \nabla \partial_{y}^{2} \phi\|_{L^{2}}^{2}  &\le C \| \nabla \partial_{x}^{2} \phi\|_{L^{2}}^{2} + C\| \nabla\phi_{t}\|_{L^{2}}^{2}+ C\| G_{3}\|_{H^{1}}^{2}.
\end{align}
Together with $(\ref{A-E1})$ and $(\ref{A-E2})$, we obtain estimates of  $ \| \nabla^{3} \phi\|_{L^{2}}^{2} $. With the estimates of $ \| \nabla^{3} \phi\|_{L^{2}}^{2} $ and $ \|  \phi\|_{L^{2}}^{2} $ at hand, we can get the estimates of $ \| \nabla^{2} \phi\|_{L^{2}}^{2} $.
The proof of Proposition $\ref{P-phi}$ is completed.
\end{proof}

The last linear problem is the linear parabolic problem of ${\bf u}$ with given $G_{4}$ and $g_{4}$:
\begin{align}\label{L-u}
\begin{cases}
{\bf u}_t -  {\rm div}\mathbb{S}({\bf u}) + \nabla p = G_{4}, &~~\text{in}~~\Omega,
\\
{\rm div}{\bf u}=0, &~~\text{in}~~\Omega,
\\
u_{2} = 0, &~~\text{on}~~\Gamma,
\\
\beta u_1+ \partial_{\bf n} u_1= g_{4}, &~~\text{on}~~\Gamma,
\\
{\bf u}\big|_{t=0} = {\bf u}_{0}, & ~~\text{in}~~\Omega.
\end{cases}
\end{align}

\begin{proposition}\label{P-u}
Given $G_{4} \in L^{2} (0,T;H^{1}) $,  $ \partial_{t}G_{4} \in  L^{2}(0,T;L^{2})$ and $g_{4} \in   L^{2}(0,T;H^{\frac32}(\Gamma)) $ and $\partial_{t}g_{4} \in L^{2}(0,T;H^{\frac12}(\Gamma))$, and assume that ${\bf u}_{0} \in  H^{2}$. Then the problem $(\ref{L-u})$ has a unique solution ${\bf u}$ and
\begin{align}\label{L-u-E}
& \sup\limits_{0\le t \le T} \left( \| {\bf u}\|_{H^{2}}^{2}  +  \| {\bf u}_{t}\|_{L^{2}} ^{2} \right) + \int_{0}^{T} \left(  \| {\bf u}\|_{H^{3}}^{2}  +  \| {\bf u}_{t}\|_{H^{1}} ^{2}   \right) {\rm d}t
\nonumber \\
& \le  C \left[ \| {\bf u}_{0}\|_{H^{2}}^{2}  + \|G_{4}(0)\|^2_{L^2} +  \int_{0}^{T} \mathcal{D}_{4}(t) {\rm d}t \right],
\end{align}
where
\begin{align*}
\mathcal{D}_{4}(t) := \varepsilon \| \partial_{t} G_4\|_{L^{2}}^{2} +  \|  \partial_{t} g_4 \|_{H^{-\frac12}(\Gamma)}^{2} + \|G_{4}  \|^2_{H^1}  + \| g_{4} \|_{H^{\frac32}(\Gamma)}^{2},
\end{align*}
in which $ \varepsilon$ is a sufficiently small constant to be determined below.
\end{proposition}

\begin{proof}[\bf{Proof}]
The existence and uniqueness of the problem $(\ref{L-u})$ can be derived directly from the result in \cite{N-T} (see Theorem 2.3 therein), we only focus on the regularity of solution.

Testing $(\ref{L-u})_{1}$ by ${\bf u}$ and integrating over $\Omega$ by parts, we find
\begin{align*}
 \frac{1}{2}\frac{\rm d}{{\rm d}t} \| {\bf u}\|_{L^{2}}^{2} +   \| \mathbb{S}({\bf u})\|_{L^{2}}^{2} & + \beta \| u_1\|_{L^{2}(\Gamma)}^{2}
 =  \int_{\Omega}  G^4 {\bf u} {\rm d}x{\rm d}y  + \int_{\Gamma}  u_1  g^4{\rm d}x
\nonumber\\
& \le \frac{\beta}{2} \| u_1\|_{L^{2}(\Gamma)}^{2}
 + \| {\bf u}\|_{L^{2}}^{2} +  \| G_4\|_{L^{2}}^{2} +  C(\beta^{-1})\|  g_4 \|_{L^{2}(\Gamma)}^{2},
\end{align*}
which implies that
\begin{align}\label{A-E-u-0}
\sup\limits_{0\le t \le T} \| {\bf u}\|_{L^{2}}^{2}  & +   \int_{0}^{T} \left(\| \mathbb{S}({\bf u})\|_{L^{2}}^{2} + \beta \| u_1\|_{L^{2}(\Gamma)}^{2} \right) {\rm d}t
\nonumber\\
& \le C(T)\left[ \| {\bf u}_{0}\|_{L^{2}}^{2} +     \int_{0}^{T} \left( \| G_4\|_{L^{2}}^{2} +  \|  g_4 \|_{L^{2}(\Gamma)}^{2} \right) {\rm d}t  \right].
\end{align}

Operating $\partial_{t}$ to $(\ref{L-u})_{1}$ yields that
\begin{align*}
 {\bf u}_{tt} -  {\rm div}\mathbb{S}({\bf u}_{t}) + \nabla p_{t} = \partial_{t}G_{4}.
\end{align*}
Testing it by ${\bf u}_{t}$ and integrating over $\Omega$ by parts, we get
\begin{align}\label{A-u-0}
 &\frac{1}{2}\frac{\rm d}{{\rm d}t} \| {\bf u}_{t}\|_{L^{2}}^{2} +   \| \mathbb{S}({\bf u}_{t})\|_{L^{2}}^{2} + \beta \| \partial_{t}u_1\|_{L^{2}(\Gamma)}^{2}
\nonumber \\
 & =  \int_{\Omega}  \partial_{t}G_4 {\bf u}_{t} {\rm d}x{\rm d}y  + \int_{\Gamma}  \partial_{t} u_1  \partial_{t}g_4{\rm d}x
\nonumber\\
& \le \epsilon \| \partial_{t} u_1\|_{H^{\frac12}(\Gamma)}^{2}
 + C \| {\bf u}_{t}\|_{L^{2}}^{2} +  \varepsilon\| \partial_{t} G_4\|_{L^{2}}^{2} + C \|  \partial_{t} g_4 \|_{H^{-\frac12}(\Gamma)}^{2}
 \nonumber\\
& \le \frac{1}{2} \| \mathbb{S}({\bf u}_{t})\|_{L^{2}}^{2}
 + C\| {\bf u}_{t}\|_{L^{2}}^{2} + \varepsilon  \| \partial_{t} G_4\|_{L^{2}}^{2} + C \|  \partial_{t} g_4\|_{H^{-\frac12}(\Gamma)}^{2}.
\end{align}
 By integration by parts, we have
\begin{align}
\|{\bf u}_t(t) \|^2_{L^2} & = \int_\Omega {\rm div}\mathbb{S}({\bf u})\cdot{\bf u}_t \mathrm{d}x + \int_\Omega G_{4}\cdot{\bf u}_t \mathrm{d}x
\nonumber \\
& \le \dfrac{1}{2}\|{\bf u}_t(t) \|^2_{L^2} + C \left ( \|\Delta{\bf u}(t) \|^2_{L^2} + \|G_4(t) \|^2_{L^2} \right). \nonumber
\end{align}
Taking $t\to 0$, we get
\begin{align}
\|{\bf u}_t(0) \|^2_{L^2} \le   C \left ( \|\Delta{\bf u}(0) \|^2_{L^2} + \|G_{4}(0) \|^2_{L^2} \right) . \nonumber
\end{align}
Using the Gronwall inequality, we have
\begin{align}\label{A-E-u-t}
\sup\limits_{0\le t \le T} \| {\bf u}_{t}\|_{L^{2}}^{2}  & +   \int_{0}^{T} \left(\| \mathbb{S}({\bf u}_{t})\|_{L^{2}}^{2} + \beta \| \partial_{t}u_1\|_{L^{2}(\Gamma)}^{2} \right) {\rm d}t
\nonumber\\
& \le C(T)\left[ \| {\bf u}_{0}\|_{H^{2}}^{2}  + \|G_4(0) \|^2_{L^2} +     \int_{0}^{T} \left(   \varepsilon \| \partial_{t} G_4\|_{L^{2}}^{2} +  \|  \partial_{t} g_4 \|_{H^{-\frac12}(\Gamma)}^{2}\right) {\rm d}t  \right].
\end{align}

Let $\omega = \nabla \times {\bf u}$. Operating ${\nabla \times}$ to $(\ref{L-u})_{1}$ and using  Lemma $\ref{lemma-curl-u}$, we have
\begin{align}\label{A-curl-u}
\begin{cases}
 \omega_t - \Delta\omega  =  \nabla \times G_{4}  , & {\rm in}~~\Omega,
\\
 \omega   =\pm\big(\beta u_{1} - g_{4}\big), & {\rm on}~~ \{ y =\pm 1 \}.
\end{cases}
\end{align}
Testing $(\ref{A-curl-u})_1$ by $ \omega - \Delta\omega $ and integrating over $\Omega$ by parts, we have
\begin{align}\label{A-E-u-2}
& \dfrac{1}{2} \dfrac{\rm d}{{\rm d}t} \| \omega  \|^2_{H^1} + \| \nabla\omega  \|^2_{L^2} + \| \Delta\omega  \|^2_{L^2}
\nonumber \\
& = \int_{\Omega}  \nabla \times G_{4} (\omega - \Delta\omega)  {\rm d}x{\rm d}y + \int_{\Gamma} \partial_{\bf n} \omega \left( \omega + \omega_{t} \right) {\rm d}x
\nonumber \\
& \le \frac{1}{4} \| \Delta\omega  \|^2_{L^2}  + \| \omega  \|^2_{L^2}  +  \|G_{4}  \|^2_{H^1}
+ \epsilon \| \partial_{\bf n} \omega \|_{H^{\frac12}(\Gamma)} ^{2} + C \| \omega \|_{H^{-\frac12}(\Gamma)}^{2}+ C\| \omega_{t} \|_{H^{-\frac12}(\Gamma)}^{2}
\nonumber \\
& \le \frac{1}{4} \| \Delta\omega  \|^2_{L^2}  + C\| \omega  \|^2_{L^2}  +  \|G_{4}  \|^2_{H^1}
+ C\epsilon \| \nabla\omega \|_{H^{1}}^{2} + C\| \omega_{t} \|_{L^{2}}^{2}
\nonumber \\
& \le \frac{1}{2} \| \Delta\omega  \|^2_{L^2}  + C\| \omega  \|^2_{L^2}  +  \|G_{4}  \|^2_{H^1}
+ C \| u_1 \|_{H^{\frac32}(\Gamma)}^{2} + \| g_{4} \|_{H^{\frac32}(\Gamma)}^{2} + C\| \mathbb{S}({\bf u}_{t}) \|_{L^{2}}^{2},
\end{align}
where we have used Korn's inequality and the following fact
\begin{align*}
 \| \nabla\omega \|_{H^{1}}^{2} \le C \| \Delta\omega  \|^2_{L^2}  + C \| u_1 \|_{H^{\frac32}(\Gamma)}^{2} + \| g_{4} \|_{H^{\frac32}(\Gamma)}^{2}.
 \end{align*}
Applying  Lemma $\ref{Trace}$ and Lemma $\ref{uuu}$, we infer
\begin{align}\label{A-E-u-b}
 \| u_1 \|_{H^{\frac32}(\Gamma)}^{2}  \le C \| {\bf u} \|_{H^{2}}^{2} \le C\| \omega  \|^2_{H^1} + C\| {\bf u}  \|^2_{L^2}.
\end{align}
Putting $(\ref{A-E-u-b})$ into $(\ref{A-E-u-2})$, using Gronwall inequality, $(\ref{A-E-u-0})$ and $(\ref{A-E-u-t})$,  and noting that $\| \omega_{0}  \|^2_{H^1} \le C \| {\bf u}_{0}  \|^2_{H^2}$, we have
\begin{align*}
&\sup\limits_{0\le t \le T} \| \omega  \|^2_{H^1} +  \int_{0}^{T} \left( \| \nabla\omega  \|^2_{L^2} + \| \Delta\omega  \|^2_{L^2} \right) {\rm d}t
\nonumber \\
& \le C \left[ \| {\bf u}_{0}\|_{H^{2}}^{2}  + \|G_4(0) \|^2_{L^2} +  \int_{0}^{T} \left( \varepsilon\| \partial_{t} G_4\|_{L^{2}}^{2} +  \|  \partial_{t} g_4 \|_{H^{-\frac12}(\Gamma)}^{2} + \|G_{4}  \|^2_{H^1}  + \| g_{4} \|_{H^{\frac32}(\Gamma)}^{2} \right) {\rm d}t \right].
\end{align*}
Applying elliptic estimates and Lemma $\ref{uuu}$, we obtain the conclusion.
\end{proof}

\subsection{Solvability of the problem \eqref{RA-NSCH}, \eqref{A-BC}, \eqref{A-IC} } \label{S-A}
We employ a fixed point argument in order to produce a solution to the problem  \eqref{RA-NSCH}, \eqref{A-BC}, \eqref{A-IC}.
\begin{theorem}\label{d-TH}
Suppose that $\Omega=\mathbb{T}\times (-1,1)$, $({\bf u}_0,\phi_{0})\in \mathbb{D}$, $ {\bf tr}(\phi_0) \in H^3(\Gamma)$ and $\mu_{0}^{\delta} \in H^{1}$. Then there is a positive constant $T^{*}$, which may depend on the initial data and on $\delta$, such that there exists a unique local strong solution $({\bf u}, \phi)$ to the problem \eqref{RA-NSCH}, \eqref{A-BC}, \eqref{A-IC}, satisfying
\begin{align}
&{\bf u} \in L^{\infty}(0,T^{*};H^{2}) \cap L^{2}(0,T^{*};H^{3}), ~\,\,\, {\bf u}_{t} \in L^{\infty}(0,T^{*};L^{2}) \cap L^{2}(0,T^{*};H^{1}),
\nonumber \\
& \phi \in L^{\infty}(0,T^{*};H^{3}) \cap L^{2}(0,T^{*};H^{4}),~\,\,\, \phi_{t} \in L^{\infty}(0,T^{*};H^{1}) \cap L^{2}(0,T^{*};H^{2}). \nonumber
\end{align}
\end{theorem}

\begin{proof}[\bf{Proof}] For $0<T<1$ and $M>1$, we define the metric space in which to work:
\begin{align}
\mathbb{X}(M,T) = &\Big\{({\bf U},\Phi) \Big| {\bf U}=(U_{1},U_{2}),\Phi \text{ satisfy that}\,   ({\bf U},\Phi)\big|_{t=0} = ({\bf u}_{0}, \phi_{0}),\nonumber\\
 &\ {\rm div}\  {\bf U} = 0 \ \text{in}\ \Omega, \ {\bf U}\cdot {\bf n} = 0\ \text{on}\ \Gamma, \, \text{and}
\, \mathcal{E}({\bf U},\Phi) + \int_{0}^{T} \mathcal{D}({\bf U},\Phi)  {\rm d}t \le M \Big\},
\end{align}
where
\begin{align*}
\mathcal{E}({\bf U},\Phi)&:= \| \Phi \|_{H^{3}}^{2} + \| \Phi_{t} \|_{H^{1}}^{2}+ \| {\bf U} \|_{H^{2}}^{2}  +  \| {\bf U}_{t}\|_{L^{2}} ^{2},
\\
\mathcal{D}({\bf U},\Phi)&:=   \| \Phi_{t} \|_{H^{2}}^{2} + \| \Phi_{tt} \|_{L^{2}}^{2} + \| \Phi \|_{H^{4}}^{2} +\| {\bf U} \|_{H^{3}}^{2}  +  \| {\bf U}_{t}\|_{H^{1}} ^{2} .
\end{align*}
Moreover, let us introduce
 $$\| ({\bf U},\Phi) \|_{\mathbb{X}} :=  \mathcal{E}({\bf U},\Phi) + \int_{0}^{T} \mathcal{D}({\bf U},\Phi)  {\rm d}t .$$ It is easy to see that $\mathbb{X}(M,T)$ is a Banach space. To perform the fixed point argument, we define a mapping
 $$\mathcal{M}: \mathbb{X}(M,T) \to \mathbb{X}(M,T) \ \text{ as }\ \mathcal{M}({\bf U},\Phi) = ({\bf u}, \phi),$$ where ${\bf u}$ and $\phi$ are determined as follows.

 For given $ ({\bf U},\Phi) \in \mathbb{X}(M,T) $, we consider the linearized system of the problem \eqref{RA-NSCH}, \eqref{A-BC}, \eqref{A-IC} as follows
\begin{align}\label{L-RA-NSCH}
\begin{cases}
{\bf u}_t- {\rm div}\mathbb{S}({\bf u})+ \nabla p =- {\bf U}\cdot\nabla {\bf U}+ \mu \nabla \Phi, &~~\text{in}~~\Omega,
\\
{\rm div}{\bf u}=0, &~~\text{in}~~\Omega,
\\
\mu - \delta\Delta\mu = -\Delta\Phi+f(\Phi) +  \delta^{2} \partial_x^2 \Phi- \delta {\bf U}\cdot \nabla \Phi, &~~\text{in}~~\Omega,
\\
\delta\phi_{t}-\Delta\phi = \mu - f(\Phi), &~~\text{in}~~\Omega,
\end{cases}
\end{align}
with the boundary condition
\begin{align}\label{L-RA-BC}
\begin{cases}
u_{2}= 0, &~~\text{on}~~\Gamma,
\\
\partial_{\bf n}\mu =0, &~~\text{on}~~\Gamma,
\\
\beta u_1+\partial_{\bf n} u_1= \big( \partial_{\bf n}\Phi+\gamma_{fs}'(\Phi)\big) \partial_{x} \Phi, &~~\text{on}~~\Gamma,
\\
\psi_t - \delta\partial_x^2 \psi = - U_{1}  \partial_{x} \Phi  -\partial_{\bf n}\Phi-\gamma_{fs}'(\Phi), &~~\text{on}~~\Gamma,
\\
\phi\big|_{\Gamma} = \psi, &~~\text{on}~~ \Gamma,
\end{cases}
\end{align}
and the initial condition
\begin{equation} \label{L-RA-IC}
({\bf u},\phi)\big|_{t=0}=({\bf u}_0,\phi_0),  ~~\text{in}~~\Omega, ~\,\,\, \psi \big|_{t=0} = \psi_{0} := \phi_{0}\big |_{\Gamma} , ~~ \text{on}~~\Gamma.
\end{equation}
The direct calculations yield that
\begin{align}\label{L-f-T}
\| f(\Phi) \|_{H^{1}} & = \left\| \int_{0}^{t} f_{t}(\Phi(x,s)) {\rm d}s +f(\Phi(x,0)) \right\|_{H^{1}}
 \le \int_{0}^{t} \| (3\Phi^{2}-1) \Phi_{t}\|_{H^{1}}   {\rm d}s + C\|  \phi_{0} \|_{H^{2}}^{3}
\nonumber \\
& \le \int_{0}^{t} \left( \| (3\Phi^{2}-1)\|_{L^{\infty}}  \| \Phi_{t}\|_{H^{1}} +  6\| \Phi\|_{L^{\infty}} \| \nabla\Phi\|_{L^{\infty}} \| \Phi_{t}\|_{L^{2}} \right) {\rm d}s + C\|  \phi_{0} \|_{H^{2}}^{3}
\nonumber \\
& \le CM^{\frac32}T + C\|  \phi_{0} \|_{H^{2}}^{3}.
\end{align}

Given $ ({\bf U},\Phi) \in \mathbb{X}(M,T) $, we define $\mu$ as the solution to problem  $(\ref{L-mu})$ with $G_1 = -\Delta\Phi+f(\Phi) +  \delta^{2}\partial_x^2 \Phi- \delta {\bf U}\cdot \nabla \Phi $. Besides, note that
$$
G_1=\left( -\Delta +\partial_x^2\right) \Phi + \left[ f(\Phi)- \delta {\bf U}\cdot \nabla \Phi \right], \ \text{i.e.}, \  G_{11} = \Phi,  G_{12}=f(\Phi)- \delta {\bf U}\cdot \nabla \Phi,
$$
we infer from $(\ref{L-mu-E})$ and $(\ref{L-mu-E2})$ that
\begin{align}\label{L-mu-T}
\| \mu_{t} \|_{H^{1}}^{2} +  \| \mu \|_{H^{3}}^{2}  \le  C(\delta^{-1}) \left(  \| \nabla\partial_{t} G_{11} \|_{L^{2}}^{2} + \| \partial_{t} G_{12} \|_{L^{2}}^{2} +  \| G_1 \|_{H^{1}}^{2} \right)\le CM^{3},
\end{align}
which implies that
\begin{align}\label{L-mu-TT}
\| \mu \|_{H^{1}} = \left\| \int_{0}^{t} \mu_{t}(x,s) {\rm d}x + \mu_{0}^{\delta} \right\|_{H^{1}} \le \int_{0}^{t} \|\mu_{t}(x,s) \|_{H^{1}} {\rm d}s + \|  \mu^{\delta}_{0} \|_{H^{1}} \le CM^{\frac32}T + \|  \mu^{\delta}_{0} \|_{H^{1}},
\end{align}
in which $\mu^{\delta}_{0}$ is defined as in $(\ref{mu-000})$.

With this, we define $\psi$ as the solution to problem  $(\ref{L-psi})$ with $G_2 = - U_{1} \partial_{x} \Phi  -\partial_{\bf n}\Phi-\gamma_{fs}'(\Phi) $ and the initial data $\psi_{0}$. Indeed, we find that
\begin{align*}
\| \partial_{t} G_2 \|_{L^{2}(\Gamma)}^{2}  &\le C \| \partial_{t} U_{1} \|_{L^{2}(\Gamma)}^{2} \| \partial_{x} \Phi \|_{L^{\infty}(\Gamma)}^{2}  +  C \|  U_{1} \|_{L^{\infty}(\Gamma)}^{2}  \| \partial_{x} \Phi_{t} \|_{L^{2}(\Gamma)}^{2}
\\
&~\,\,\,+  C \| \partial_{\bf n}\Phi_{t}  \|_{L^{2}(\Gamma)}^{2}  + C \| \Phi_{t}  \|_{L^{2}(\Gamma)}^{2}
\\
& \le C \| {\bf U}_{t} \|_{H^{\frac12}}^{2} \| \partial_{x} \Phi \|_{H^{1}(\Gamma)}^{2} + C( \|  U_{1} \|_{H^{1}(\Gamma)}^{2} +1 ) \| \nabla \Phi_{t} \|_{H^{\frac12}}^{2} + C \|  \Phi_{t} \|_{H^{\frac12}}^{2}
\\
& \le C \| {\bf U}_{t} \|_{L^{2}} \| {\bf U}_{t} \|_{H^{1}} \| \Phi \|_{H^{3}}^{2} + C( \|  {\bf U} \|_{H^{2}}^{2} +1 ) \| \nabla \Phi_{t} \|_{L^{2}} \| \nabla \Phi_{t} \|_{H^{1}}  +  C\| \Phi_{t} \|_{H^{1}}^{2}
\\
& \le C M^{\frac{3}{2}} ( \| {\bf U}_{t} \|_{H^{1}} +  \| \nabla \Phi_{t} \|_{H^{1}} ) +  C M,
\end{align*}
and
\begin{align*}
 \| G_2 \|_{H^{2}(\Gamma)}^{2} & \le C \| \partial_{x}^{2} U_{1} \|_{L^{2}(\Gamma)}^{2} \| \partial_{x} \Phi \|_{L^{\infty}(\Gamma)}^{2}  +  C\| \partial_{x} U_{1} \|_{L^{4}(\Gamma)}^{2} \| \partial_{x}^{2} \Phi \|_{L^{4}(\Gamma)}^{2}
 \\
& ~\,\,\, +  \|  U_{1} \|_{L^{\infty}(\Gamma)}^{2} \| \partial_{x}^{3} \Phi \|_{L^{2}(\Gamma)}^{2}  +  \| \partial_{\bf n} \Phi \|_{H^{2}(\Gamma)}^{2} + \| \Phi \|_{H^{2}(\Gamma)}^{2}
 \\
 & \le C  \| {\bf U} \|_{H^{\frac52}}^{2} \| \partial_{x} \Phi \|_{H^{1}(\Gamma)}^{2}+  C\| \partial_{x} U_{1} \|_{H^{\frac12}(\Gamma)}^{2}  \| \partial_{x}^{2} \Phi \|_{H^{\frac12}(\Gamma)}^{2}
  \\
& ~\,\,\, + C \|  U_{1} \|_{H^{1}(\Gamma)}^{2} \| \partial_{x}^{3} \Phi \|_{H^{\frac12}}^{2}  +  C\| \partial_{\bf n} \Phi \|_{H^{\frac52}}^{2} + C\| \Phi \|_{H^{\frac52}}^{2}
 \\
 & \le C  \| {\bf U} \|_{H^{2}} \| {\bf U} \|_{H^{3}}  \| \Phi \|_{H^{3}}^{2}+  C\|  {\bf U}\|_{H^{2}}^{2}  \| \Phi \|_{H^{3}}^{2}
  \\
& ~\,\,\, +  C\|  {\bf U} \|_{H^{2}}^{2} \| \Phi \|_{H^{3}}  \| \Phi \|_{H^{4}}  +   C\| \Phi \|_{H^{3}}  \| \Phi \|_{H^{4}} + C\| \Phi \|_{H^{3}}^{2}
\\
& \le C  M^{\frac{3}{2}} ( \| {\bf U} \|_{H^{3}} + \| \Phi \|_{H^{4}} )+  CM^{2},
\end{align*}
which yield that
\begin{align}\label{D-11}
\int_{0}^{T} \mathcal{D}_{2}(t) {\rm d}t =  \int_{0}^{T} \left( \| \partial_{t} G_2 \|_{L^{2}(\Gamma)}^{2}  +  \| G_2 \|_{H^{2}(\Gamma)}^{2} \right) {\rm d}t \le CM^{2}T^{\frac{1}{2}} + CMT.
\end{align}
Meanwhile, it follows from $(\ref{L-psi-E})$ that
\begin{align}\label{L-psi-T}
&\sup\limits_{0\le t \le T} \left(  \| \psi_{t} \|_{H^{1}(\Gamma)}^{2} +  \| \psi \|_{H^{3}(\Gamma)}^{2} \right)
+ \int_{0}^{T}\left( \| \psi_{tt} \|_{L^{2}(\Gamma)}^{2} + \| \psi_{t} \|_{H^{2}(\Gamma)}^{2} +  \| \psi \|_{H^{4}(\Gamma)}^{2} \right) {\rm d}t
\nonumber \\
&\le  C(\delta^{-1})\left( \| \psi_{0} \|_{H^{3}(\Gamma)}^{2} +  \| {\bf u}_{0} \|_{H^{2}}^{2}\|\phi_{0} \|_{H^{3}}^{2} +  M^{2}T^{\frac{1}{2}} + MT \right),
\end{align}
where \eqref{D-11} has been used.


With $\mu$ and $\psi$ at hand, we define $\phi$ as the solution to problem $(\ref{L-phi})$ with $G_{3} =
\mu - f(\Phi) $, $g_{3}=\psi$ and the initial data $\phi_{0}$, where $\mu$ and $\psi$ are the above functions constructed before.
Thanks to $(\ref{L-mu-T})$ and  $(\ref{L-psi-T})$, one has
\begin{align*}
& \| \partial_{t}G_{3}\|_{L^{2}}^{2}   \le   \| \mu_{t} \|_{L^{2}}^{2} + \| (3\Phi^{2} - 1)\Phi_{t} \|_{L^{2}}^{2}  \le  CM^{3},
\\[1em]
&\| G_{3}\|_{H^{2}}^{2} \le  \| \mu \|_{H^{2}}^{2} + \| \Phi^{3} - \Phi \|_{H^{2}}^{2}  \le CM^{3},
\end{align*}
and
\begin{align*}
& \| \partial_{t}^{2}  g_{3} \|_{H^{-\frac12}(\Gamma)}^{2} +  \| \partial_{t}  g_{3} \|_{H^{\frac32}(\Gamma)}^{2} + \|   g_{3} \|_{H^{\frac72}(\Gamma)}^{2}  \le  C\| \psi_{tt} \|_{L^{2}(\Gamma)}^{2} +C \| \psi_{t} \|_{H^{2}(\Gamma)}^{2} +C \| \psi \|_{H^{4}(\Gamma)}^{2} .
\end{align*}
Thus, the estimate $(\ref{L-phi-E})$ together with the above estimates $(\ref{L-f-T})$, $(\ref{L-mu-TT})$ and $(\ref{L-psi-T})$ implies
\begin{align}\label{L-phi-T}
& \sup\limits_{0\le t \le T} \left( \| \phi \|_{H^{3}}^{2} + \| \phi_{t} \|_{H^{1}}^{2} \right) + \int_{0}^{T} \left( \| \phi_{t} \|_{H^{2}}^{2} + \| \phi_{tt} \|_{L^{2}}^{2} + \| \phi \|_{H^{4}}^{2} \right)  {\rm d}t
\nonumber \\
& \le  C (\delta^{-1}) \bigg[ \| \phi_{0}\|_{H^{3}}^{6} +  \|  \mu_{0}^{\delta} \|_{H^{1}}^{2} +\| \psi_{0} \|_{H^{3}(\Gamma)}^{2} +  \| {\bf u}_{0} \|_{H^{2}}^{2}\|\phi_{0} \|_{H^{3}}^{2} +  M^{2}T^{\frac{1}{2}} +  M^{3}T  \bigg].
\end{align}

Finally, we define ${\bf u}$ as the solution to the problem $(\ref{L-u})$ with $G_{4} = - {\bf U} \cdot \nabla {\bf U} + \mu \nabla\Phi$, $g_{4}= ( \partial_{\bf n}\Phi+\gamma_{fs}'(\Phi) ) \partial_{x} \Phi$ and the initial data ${\bf u}_{0}$. Moreover, choosing $\varepsilon = M^{-2}$, we get from the estimate $(\ref{L-u-E})$ that
\begin{align}\label{L-u-T}
& \sup\limits_{0\le t \le T} \left( \| {\bf u}\|_{H^{2}}^{2}  +  \| {\bf u}_{t}\|_{L^{2}} ^{2} \right) + \int_{0}^{T} \left(  \| {\bf u}\|_{H^{3}}^{2}  +  \| {\bf u}_{t}\|_{H^{1}} ^{2}   \right) {\rm d}t
\nonumber \\
& \le  C \left[ \| {\bf u}_{0}\|_{H^{2}}^{4}  + \|\mu_{0}^{\delta} \|^2_{L^2} \|\phi_{0} \|^2_{H^3} + M^{4}T + 1\right],
\end{align}
where we have used the following fact
\begin{align*}
\varepsilon \| \partial_{t} G_4\|_{L^{2}}^{2} & \le  \varepsilon \| {\bf U}_{t} \|_{L^{2}}^{2} \| \nabla {\bf U} \|_{L^{\infty}}^{2} + \varepsilon \| {\bf U} \|_{L^{\infty}}^{2} \| \nabla {\bf U}_{t} \|_{L^{2}}^{2} + \varepsilon   \| \mu_{t} \|_{L^{2}}^{2} \| \nabla \Phi \|_{L^{\infty}}^{2} + \varepsilon  \| \mu \|_{L^{\infty}}^{2} \| \nabla \Phi_{t} \|_{L^{2}}^{2}
\\
& \le C \varepsilon M \left(  \| {\bf U} \|_{H^{3}}^{2} + \| \nabla {\bf U}_{t} \|_{L^{2}}^{2} + M^{3} \right),
\\[1em]
\|G_{4}  \|^2_{H^1}  &\le C \| {\bf U} \|_{L^{\infty}}^{2} \| \nabla {\bf U} \|_{L^{2}}^{2} + C \| \nabla {\bf U}\|_{L^{4}}^{2} \| \nabla {\bf U}\|_{L^{4}}^{2} + C \| {\bf U} \|_{L^{\infty}}^{2} \| \nabla^{2} {\bf U} \|_{L^{2}}^{2}
\\
&  ~\,\,\, + C \| \mu \|_{L^{4}}^{2} \| \nabla \Phi \|_{L^{4}}^{2} + C \| \mu \|_{L^{\infty}}^{2} \| \nabla ^{2}\Phi \|_{L^{2}}^{2} \le CM^{4},
\\[1em]
\| g_{4} \|_{H^{\frac32}(\Gamma)}^{2} &\le \| \partial_{\bf n}\Phi + \gamma_{fs}'(\Phi) \|_{H^{\frac32}(\Gamma)}^{2} \|\partial_{x} \Phi \|_{H^{\frac32}(\Gamma)}^{2} \le  C\left( \|  \partial_{\bf n}\Phi \|_{H^{2}}^{2}+1\right) \|\partial_{x} \Phi \|_{H^{2}}^{2} \le CM^{2},
\end{align*}
and
\begin{align*}
 \| \partial_{t} g_4 \|_{H^{-\frac12}(\Gamma)}^{2}  & \le C  \|  \partial_{\bf n}\Phi_{t}+ \Phi_{t} \|_{L^{2}}^{2} \|\partial_{x} \Phi \|_{L^{\infty}}^{2} + C \| \partial_{\bf n}\Phi+ \gamma_{fs}'(\Phi)\|_{L^{\infty}}^{2} \|\partial_{x} \Phi_{t} \|_{L^{2}}^{2}  \le CM^{2}.
\end{align*}
Therefore, we take $M$ as follows
\begin{align*}
M := C\big( \| \phi_{0}\|_{H^{3}}^{6} +  \|  \mu_{0}^{\delta} \|_{H^{1}}^{2} +\| \psi_{0} \|_{H^{3}(\Gamma)}^{2} +  \| {\bf u}_{0} \|_{H^{2}}^{2}\|\phi_{0} \|_{H^{3}}^{2} +   \| {\bf u}_{0}\|_{H^{2}}^{4}  + \|\mu_{0}^{\delta} \|^2_{L^2} \|\phi_{0} \|^2_{H^3} +4 \big),
\end{align*}
and choose $T$ is sufficiently small such that $ M^{4}T<1 $. Then $({\bf u}, \phi) \in \mathbb{X}(M,T)$, which implies that the mapping $\mathcal{M}: \mathbb{X}(M,T) \to \mathbb{X}(M,T)$ is well-defined.

Next we show that the mapping $\mathcal{M}$ has a fixed point in the space $\mathbb{X}(M,T)$ by proving the contraction. To this end, let $({\bf U}_{i},\Phi_{i}) \in \mathbb{X}(M,T)$, $ {\bf U}_{i} = (U_{1}^{i}, U_{2}^{i})$ and $({\bf u}_{i},\phi_{i})= \mathcal{M}({\bf U}_{i},\Phi_{i}) \in \mathbb{X}(M,T), i=1,2$. Define
\begin{align*}
 \Phi &: =\Phi_{1}-\Phi_{2},\ {\bf U}:={\bf U}_{1}-{\bf U}_{2}=(U_{1},U_{2})=(U_{1}^{1}-U_{1}^{2},U_{2}^{1}-U_{2}^{2}),
 \\
f(\Phi)&:=f(\Phi_{1})-f(\Phi_{2}) = \Phi \left( \Phi_{1}^{2} + \Phi_{1}\Phi_{2} + \Phi_{2}^{2} -1\right) ,
\\
\gamma_{fs}'(\Phi)& := \gamma_{fs}'(\Phi_{1}) - \gamma_{fs}'(\Phi_{2}) = \gamma_{fs}^{(2)}(\theta\Phi_{1} + (1-\theta)\Phi_{2}) \Phi,
 \end{align*}
 where $\theta\in(0,1)$.
 The estimate $(\ref{L-f-T})$ becomes that
\begin{align}\label{L-f-T-C}
\| f(\Phi) \|_{H^{1}} & = \left\| \int_{0}^{t} f_{t}(\Phi(x,s)) {\rm d}s \right\|_{H^{1}}
 \nonumber \\
& \le \int_{0}^{t}\left(  \| \left( \Phi_{1}^{2} + \Phi_{1}\Phi_{2} + \Phi_{2}^{2} -1\right) \Phi_{t}\|_{H^{1}}+ \| \left( \Phi_{1}^{2} + \Phi_{1}\Phi_{2} + \Phi_{2}^{2} -1\right)_{t} \Phi\|_{H^{1}}  \right)  {\rm d}s
 \nonumber \\
& \le CMT \| ({\bf U},\Phi) \|_{\mathbb{X}}^{\frac12}.
\end{align}
We can find $ \mu = \mu_{1}-\mu_{2}$ that solves $(\ref{L-mu})$ with
\begin{align*}
G_1=\left( -\Delta +\partial_x^2\right) \Phi + \left[ f(\Phi)- \delta {\bf U} \cdot \nabla \Phi_{1} - \delta {\bf U}_{2} \cdot \nabla \Phi \right],
\end{align*}
i.e.,
\begin{align*}
G_{11} = \Phi,  G_{12}=f(\Phi)- \delta {\bf U} \cdot \nabla \Phi_{1} - \delta {\bf U}_{2} \cdot \nabla \Phi,\ G_1\big|_{t=0} = 0.
\end{align*}
One can solve $(\ref{mu-000})$ with $ u_0=\phi_0=\mu_0=0$ and obtain $\mu_{0}^{\delta} = 0$. Then $(\ref{L-mu-T})$ and $(\ref{L-mu-TT})$ are reduced to
\begin{align}\label{L-mu-T-C}
\| \mu_{t} \|_{H^{1}}^{2} +  \| \mu \|_{H^{3}}^{2}  \le  C(\delta^{-1}) \left(  \| \nabla\partial_{t} G_{11} \|_{L^{2}}^{2} + \| \partial_{t} G_{12} \|_{L^{2}}^{2} +  \| G_1 \|_{H^{1}}^{2}\right)\le CM^{2}  \| ({\bf U},\Phi) \|_{\mathbb{X}},
\end{align}
and
\begin{align}\label{L-mu-TT-C}
\| \mu \|_{H^{1}} = \left\| \int_{0}^{t} \mu_{t}(x,s) {\rm d}s  \right\|_{H^{1}} \le \int_{0}^{t} \|\mu_{t}(x,s) \|_{H^{1}} {\rm d}s  \le CMT \| ({\bf U},\Phi) \|_{\mathbb{X}}^{\frac12},
\end{align}
respectively.

We find $\psi=\psi_{1}-\psi_{2}$ which solves $(\ref{L-psi})$ with 
$$G_2 = - U_{1} \partial_{x} \Phi_{1} -  U_{1}^{2}\partial_{x} \Phi  -\partial_{\bf n}\Phi-\gamma_{fs}'(\Phi), $$ and the initial data $\psi_{0} = 0$. Meanwhile, the estimate $(\ref{L-psi-T})$ becomes into
\begin{align}\label{L-psi-T-C}
&\sup\limits_{0\le t \le T} \left(  \| \psi_{t} \|_{H^{1}(\Gamma)}^{2} +  \| \psi \|_{H^{3}(\Gamma)}^{2} \right)
+ \int_{0}^{T}\left( \| \psi_{tt} \|_{L^{2}(\Gamma)}^{2} + \| \psi_{t} \|_{H^{2}(\Gamma)}^{2} +  \| \psi \|_{H^{4}(\Gamma)}^{2} \right) {\rm d}t
\nonumber \\
&\le  C(\delta^{-1})  ( MT^{\frac{1}{2}} + T ) \| ({\bf U},\Phi) \|_{\mathbb{X}}.
\end{align}

We will find $\phi=\phi_{1}-\phi_{2}$ which solves $(\ref{L-phi})$ with $G_{3} = \mu - f(\Phi) $, $g_{3}=\psi$ and the initial data $\phi_{0} = 0$. The estimate $(\ref{L-phi-T})$ is reduced to
\begin{align}\label{L-phi-T-C}
& \sup\limits_{0\le t \le T} \left( \| \phi \|_{H^{3}}^{2} + \| \phi_{t} \|_{H^{1}}^{2} \right) + \int_{0}^{T} \left( \| \phi_{t} \|_{H^{2}}^{2} + \| \phi_{tt} \|_{L^{2}}^{2} + \| \phi \|_{H^{4}}^{2} \right)  {\rm d}t
\nonumber \\
& \le  C (\delta^{-1})   ( MT^{\frac{1}{2}} +  M^{2}T )\| (U,\Phi) \|_{\mathbb{X}} .
\end{align}

Finally, we find ${\bf u}$ which solves  $(\ref{L-u})$ with $G_{4} = - {\bf U} \cdot \nabla {\bf U}_{1}- {\bf U}_{2} \cdot \nabla {\bf U} + \mu \nabla\Phi_{1} + \mu_{2} \nabla\Phi$, $g_{4}= ( \partial_{\bf n}\Phi+\gamma_{fs}'(\Phi) ) \partial_{x} \Phi_{1} +  ( \partial_{\bf n}\Phi_{2}+\gamma_{fs}'(\Phi_{2}) ) \partial_{x} \Phi$ and the initial data ${\bf u}_{0}=0$. Moreover, choosing $ C\varepsilon <  \frac{1}{4M^{2}}$, we find that the estimates $(\ref{L-u-T})$ becomes into
\begin{align}\label{L-u-T-C}
& \sup\limits_{0\le t \le T} \left( \| {\bf u}\|_{H^{2}}^{2}  +  \| {\bf u}_{t}\|_{L^{2}} ^{2} \right) + \int_{0}^{T} \left(  \| {\bf u}\|_{H^{3}}^{2}  +  \| {\bf u}_{t}\|_{H^{1}} ^{2}   \right) {\rm d}t
 \le  \left( C M^{3}T  + 1/4 \right) \| ({\bf U},\Phi) \|_{\mathbb{X}}.
\end{align}
Combining \eqref{L-phi-T-C} and \eqref{L-u-T-C}, and choosing $T$ small enough such that $CM^{3}T<1/12$, we have
\begin{align*}
  \| ({\bf u},\phi) \|_{\mathbb{X}} \le \frac{1}{2}\| ({\bf U},\Phi) \|_{\mathbb{X}},
\end{align*}
which implies the mapping $\mathcal{M}: \mathbb{X}(M,T) \to \mathbb{X}(M,T)$ is a contraction. Therefore, the mapping $\mathcal{M}$ admits a unique fixed point in $\mathbb{X}(M,T)$, i.e., the problem $(\ref{RA-NSCH})$, $(\ref{A-BC})$, $(\ref{A-IC})$ admits a unique solution $({\bf u}, \phi) \in \mathbb{X}(M,T)$.
The proof is completed.
\end{proof}

\subsection{Uniform $\delta$-independent estimates of $(\ref{A-NSCH})$-$(\ref{A-IC})$}\label{subs-d}
For each $\delta>0$, by subsection $\ref{S-A}$, there exists $T_{\delta}=T_{\delta}({\bf u}_{0},\phi_{0},\psi_{0},\mu_{0}^{\delta})>0$ such that there is a unique solution $({\bf u}, \phi)=({\bf u}^{\delta}, \phi^{\delta})\in \mathbb{X}(M,T)$ to  the problem $(\ref{A-NSCH})$-$(\ref{A-IC})$ on $[0,T_{\delta}]$. For notational simplifications, we will not explicitly write the dependence of the solution on $\delta$. The purpose of this subsection is to obtain the $\delta$-independent estimates of the solutions to $(\ref{A-NSCH})$-$(\ref{A-IC})$, which allow us to pass the limit as $\delta\to 0$.

We now define the energy functional
\begin{align}\label{tE-delta}
\mathcal{E}^{\delta}(t) & = 1 + \| {\bf u}\|_{H^{2}}^{2} + \| {\bf u}_{t}\|_{L^{2}}^{2} + \| \phi\|_{H^{3}}^{2} + \| \phi_{t}\|_{H^{1}}^{2} +  \| \mu \|_{H^{3}}^{2}
\nonumber \\
& ~\,\,\,+ \delta^{2}\| \partial_{x}\phi_{t} \|_{L^{2}}^{2} + \delta\|\partial_{x}^{2}\phi\|_{L^{2}(\Gamma)}^{2}+  \delta^{2} \| \partial_{x}^{2}\phi \|_{L^{2}}^{2}+ \delta \| \nabla\partial_{x}^{2}\phi \|_{L^{2}}^{2} + \delta \|\partial_{x} \Delta\mu  \|_{L^{2}}^{2} ,
\end{align}
and dissipative functional
\begin{align}
\mathcal{D}^{\delta}(t) & = \| {\bf u}\|_{H^{3}}^{2} + \| {\bf u}_{t}\|_{H^{1}}^{2} + \| \Delta\phi\|_{H^{2}}^{2} + \| \phi_{t}\|_{H^{2}}^{2} +  \| \mu \|_{H^{4}}^{2} +   \| \mu_{t} \|_{H^{1}}^{2} +   \| \partial_{\bf n}\phi \|_{H^{2}(\Gamma)}^{2}
\nonumber \\
& ~\,\,\,  +   \delta  \| \phi_{tt} \|_{L^{2}}^{2}  +     \delta  \| \nabla\partial_{x}^{3} \phi \|_{L^{2}}^{2}.
\end{align}
We will show that $\mathcal{E}^{\delta}(t)$ remains bounded on a time interval independent of $\delta$, which is stated as the following theorem.
\begin{theorem}\label{P-E-d}
Under the assumptions of Theorem $\ref{TH}$. Then there exists a time $T_{0}$ independent of $\delta$ such that
\begin{align}\label{d-i-E}
\sup\limits_{0\le t\le T_{0}}\mathcal{E}^{\delta}(t) + \int_{0}^{T} \mathcal{D}^{\delta}(t)  {\rm d}t  \le  2C_{0},
\end{align}
where $C_{0}$ is a constant that depend on $\Omega$, $\beta$, $ \hat{\mathcal{E}}(0)$ and the initial value  $(\| {\bf u}_{0}\|_{H^{2}}, \| \phi_{0}\|_{H^{3}}, \| \partial_{x}^{2}\phi_{0}\|_{H^{2}},$
$\| \mu_{0}\|_{H^{3}})$, but not on $\delta$.
\end{theorem}
To prove Theorem \ref{P-E-d}, we need the following Proposition $\ref{P-mu}$.

\begin{proposition}\label{P-mu}
Under the assumptions of Theorem $\ref{TH}$. Then
\begin{align}\label{mu-0-h3}
   \| \mu^\delta_0 \|^2_{H^3}
\le  C \| {\bf u}_0\|^2_{H^2}\|\phi_0\|^2_{H^3} + C \|\partial_x^2\phi_0 \|_{H^{2}}^2 +C\| \mu_0\|^2_{H^3},
\end{align}
where $C$  is independent of $\delta$.
\end{proposition}

\begin{proof}[\bf Proof]
Let $\mu_0 \in H^1$, $\Delta\mu_0 \in H^1$, $\partial_x^2\phi_0\in H^{2}$ and $({\bf u}_0,\phi_{0})\in \mathbb{D}$, we now consider the following elliptic equation
\begin{align}\label{mu-00}
\begin{cases}
\mu^\delta_0-\delta \Delta \mu^\delta_0  = - \delta{\bf u}_0\cdot\nabla\phi_0  + \delta^{2}\partial_x^2 \phi_{0}+ \mu_0, ~&\text{in}~\Omega,
\\
\partial_{\bf n}\mu^\delta_0\big|_{\Gamma}=0, ~&\text{on}~\Gamma.
\end{cases}
\end{align}

First, multiplying $(\ref{mu-00})_1$ by $\mu^\delta_0$ and $-\Delta \mu^\delta_0$, respectively, and taking the sum of the results, we get
\begin{align}\label{aaaa}
\| \mu^\delta_0 \|^2_{H^1} + \delta \| \Delta\mu^\delta_0 \|^2_{L^2} \le \frac{\delta}{2}  \| \Delta\mu^\delta_0 \|^2_{L^2} +  \frac{1}{2} \| \mu^\delta_0 \|^2_{H^1}+ \delta  \| {\bf u}_0\cdot\nabla\phi_0 \|^2_{L^2} + \delta \|\partial_{x}^{2}\phi_0\|^2_{L^2} + \|\mu_0\|^2_{H^1}.
\end{align}

Next, it follows from $(\ref{mu-00})_1$ that
\begin{align*}
\nabla\mu^\delta_0-\delta \nabla\Delta \mu^\delta_0  =- \delta\nabla({\bf u}_0\cdot\nabla\phi_0) + \delta^{2}\nabla \partial_{x}^{2}\phi_0+ \nabla\mu_0,~\text{in}~\Omega.
\end{align*}
Multiplying it by $ -\nabla\Delta\mu^\delta_0$ and integrating the result over $\Omega$ by parts, we have
\begin{align}\label{4.3}
\| \Delta \mu^\delta_0 \|^2_{L^2}  + \delta\| \nabla\Delta \mu^\delta_0 \|^2_{L^2} \le  \delta  \| \nabla( {\bf u}_0\cdot\nabla\phi_0) \|^2_{L^2} +\delta \|\nabla\partial_{x}^{2}\phi_0\|^2_{L^2} + \|\Delta\mu_0\|^2_{L^2}.
\end{align}

And then, operating $\nabla \partial_{x}$  and $\partial_{x}^{2}$ to $(\ref{mu-00})_1$ respectively,  we have 
\begin{align}\label{mm}
\begin{cases}
\nabla \partial_{x}\mu^\delta_0-\delta \Delta \nabla \partial_{x}\mu^\delta_0  = -\delta\nabla \partial_{x}({\bf u}_0\cdot\nabla\phi_0) + \delta^{2}\nabla \partial_x^3\phi_0  + \nabla \partial_{x}\mu_0, ~&\text{in}~\Omega,
\\
\partial_{x}^{2}\mu^\delta_0-\delta \Delta  \partial_{x}^{2}\mu^\delta_0  = -\delta\partial_{x}^{2}({\bf u}_0\cdot\nabla\phi_0) + \delta^{2} \partial_x^4\phi_0 +  \partial_{x}^{2}\mu_0,~&\text{in}~\Omega.
\end{cases}
\end{align}
Multiplying \eqref{mm}$_{1}$ and \eqref{mm}$_{2}$ by $ -\Delta\nabla\partial_{x}\mu^\delta_0$ and $ -\Delta\partial_{x}^{2}\mu^\delta_0$, respectively, and integrating the result over $\Omega$ by parts, one has
\begin{align}\label{qiexiang}
& \| \Delta \partial_{x} \mu^\delta_0 \|^2_{L^2} + \| \nabla \partial_{x}^{2} \mu^\delta_0 \|^2_{L^2}  + \delta\| \nabla\Delta \partial_{x}\mu^\delta_0 \|^2_{L^2}  + \delta\| \Delta \partial_{x}^{2}\mu^\delta_0 \|^2_{L^2}
\nonumber \\
& \le  \delta  \| \nabla\partial_{x}( {\bf u}_0\cdot\nabla\phi_0) \|^2_{L^2} + \delta \|\partial_x^2\phi_0 \|_{H^{2}}^{2} + \|\nabla^{3}\mu_0\|^2_{L^2}.
\end{align}

Moreover, operating $ \partial_{\bf n}$ to $(\ref{mu-00})_1$ yields
\begin{align*}
\partial_{\bf n}\mu^\delta_0-\delta \Delta \partial_{\bf n}\mu^\delta_0  = -\delta\partial_{\bf n}({\bf u}_0\cdot\nabla\phi_0) + \delta^{2}\partial_{\bf n}\partial_{x}^{2}\phi_0 + \partial_{\bf n}\mu_0, ~\,\text{in}~\Omega,
\end{align*}
which together with $(\ref{mu-00})_2$ and  $ \partial_{\bf n} \mu_0 \big|_{\Gamma}=0$ implies
\begin{align}\label{bjbj}
-\Delta \partial_{\bf n}\mu^\delta_0 = -\partial_{\bf n}({\bf u}_0\cdot\nabla\phi_0) + \delta\partial_{\bf n}\partial_{x}^{2}\phi_0,~~\text{on}~~ \Gamma.
\end{align}

Finally, operating $ \nabla\partial_{y}$ to $(\ref{mu-00})_1$ holds
\begin{align*}
\nabla\partial_{y}\mu^\delta_0-\delta \nabla \Delta \partial_{y}\mu^\delta_0  = -\delta\nabla\partial_{y}({\bf u}_0\cdot\nabla\phi_0) + \delta^{2}\nabla\partial_{y}\partial_{x}^{2}\phi_0 + \nabla\partial_{y}\mu_0,~&\text{in}~\Omega.
\end{align*}
Multiplying it by $ -\nabla\Delta\partial_{y}\mu^\delta_0$ and integrating the result over $\Omega$ by parts, one achieves from $(\ref{bjbj})$ and Lemma $\ref{Trace}$ that
\begin{align}\label{4.88}
&\| \Delta \partial_{y}\mu^\delta_0 \|^2_{L^2}  + \delta\| \nabla\Delta \partial_{y}\mu^\delta_0 \|^2_{L^2} \nonumber\\
&\le  \dfrac{\delta}{2}\| \nabla\Delta \partial_{y}\mu^\delta_0 \|^2_{L^2} + \dfrac{1}{4}\| \Delta \partial_{y}\mu^\delta_0 \|^2_{L^2} + \int_{\Gamma} \partial_{\bf n}\partial_{y}\mu^\delta_0\, \Delta \partial_{y}\mu^\delta_0 {\rm d}x + \delta  \| \nabla\partial_{y}( {\bf u}_0\cdot\nabla\phi_0) \|^2_{L^2}
\nonumber\\
&~\,\,\,+ \delta\| \nabla\partial_{y} \partial_{x}^{2}\phi_0 \|_{L^{2}}^{2}+\|\Delta \partial_{y}\mu_0\|^2_{L^2} + \int_{\Gamma} \partial_{\bf n}\partial_{y}\mu_0\, \Delta \partial_{y}\mu^\delta_0 {\rm d}x
\nonumber\\
& =  \dfrac{\delta}{2}\| \nabla\Delta \partial_{y}\mu^\delta_0 \|^2_{L^2} + \dfrac{1}{4}\| \Delta \partial_{y}\mu^\delta_0 \|^2_{L^2} + \int_{\Gamma} \partial_{\bf n}\partial_{y}\mu^\delta_0\,  \left[-\partial_{\bf n}({\bf u}_0\cdot\nabla\phi_0) + \delta\partial_{\bf n}\partial_{x}^{2}\phi_0  \right] {\rm d}x
\nonumber\\
&~\,\,\,+ \delta  \| \nabla\partial_{y}( {\bf u}_0\cdot\nabla\phi_0) \|^2_{L^2}+ \delta\| \nabla\partial_{y}\partial_{x}^{2} \phi_0 \|_{L^{2}}^{2}+\|\Delta \partial_{y}\mu_0\|^2_{L^2}
\nonumber\\
&~\,\,\, - \int_{\Gamma} \partial_{\bf n}\partial_{y}\mu_0\,  \left[-\partial_{\bf n}({\bf u}_0\cdot\nabla\phi_0) + \delta\partial_{\bf n} \partial_{x}^{2} \phi_0  \right] {\rm d}x
\nonumber\\
& \le \dfrac{\delta}{2}\| \nabla\Delta \partial_{y}\mu^\delta_0 \|^2_{L^2} + \dfrac{1}{4}\| \Delta \partial_{y}\mu^\delta_0 \|^2_{L^2} +\dfrac{1}{4}\|  \partial_{y}^{2} \mu^\delta_0 \|^2_{H^1} + C \|\partial_{y}({\bf u}_0\cdot\nabla\phi_0) \|^2_{H^1}
\nonumber\\
&~\,\,\,+ C\delta \|\partial_x^2\phi_0 \|_{H^{2}} +C\| \mu_0\|^2_{H^3},
\end{align}
which combined with  $(\ref{4.3})$, $(\ref{qiexiang})$ and the following fact
\begin{align*}
\| \partial_{y}^{2} \mu^\delta_0 \|^2_{H^1} \le \|\Delta\mu^\delta_0 \|^2_{H^1} + \| \partial_{x}^{2} \mu^\delta_0 \|^2_{H^1},
\end{align*}
gives
\begin{align*}
\| \Delta \mu^\delta_0 \|^2_{H^1}  + \delta\| \Delta \mu^\delta_0 \|^2_{H^2}
\le C \| {\bf u}_0\|^2_{H^2} \|\phi_0\|^2_{H^3} + C \|\partial_x^2\phi_0 \|_{H^{2}} +C\| \mu_0\|^2_{H^3}.
\end{align*}
This together with $(\ref{aaaa})$ and elliptic estimates implies $(\ref{mu-0-h3})$ holds.
\end{proof}



We are on a position to prove Theorem $\ref{P-E-d}$.
\begin{proof}[\bf Proof of Theorem ${\bf \ref{P-E-d}}$]
Indeed, compared with the case of $\delta=0$, there are no differences in the equation for velocity. Therefore, we only need to consider the estimates for $\phi$ and $\mu$.

The proof will be completed by several steps.

{\bf Step 1. Behavior near boundaries.}

It follows from $(\ref{A-NSCH})_{3}$ that
\begin{align*}
\phi_{t} - \delta\partial_x^2 \phi + u_1  \partial_{x}\phi  + u_{2}\partial_{y}\phi = \Delta\mu, ~~~~\text{in}~~\Omega.
\end{align*}
We now take trace for the above equation to give
\begin{align*}
\phi_{t} - \delta\partial_x^2 \phi + u_1  \partial_{x}\phi = \Delta\mu, ~~~~\text{on}~~\Gamma ,
\end{align*}
since $u_{2} = 0$ on $\Gamma $. Thus, we find from the above equation and $(\ref{A-BC})_{4}$ that
\begin{align}\label{NN-BC}
\Delta\mu = - \partial_{\bf n}\phi - \gamma_{fs}'(\phi), ~~~~\text{on}~~\Gamma.
\end{align}

{\bf Step 2.  Basic energy.}

Multiplying $(\ref{A-NSCH})_3$ by $\mu$ and integrating by parts, we get
\begin{align}\label{A-E-2}
&\frac{1}{2}\frac{{\rm d}}{{\rm d}t} \left[\| \nabla \phi \|_{L^2}^{2}+ \frac{1}{2}\| \phi^2-1 \|_{L^2}^2 + \delta^{2} \| \partial_{x}\phi \|_{L^2}^{2} \right]+\| \nabla \mu \|_{L^2}^{2} + \delta \| \phi_{t}\|_{L^{2}}^{2} + \delta\| \nabla\partial_{x} \phi \|_{L^2}^{2}
\nonumber\\
&= -\int_{\Omega} \mu {\bf u} \cdot \nabla \phi {\rm d}x{\rm d}y + \delta \int_{\Omega} (1-3\phi^{2}) |\partial_{x}\phi|^{2} {\rm d}x{\rm d}y+ \int_{\Gamma} \left(\phi_{t} - \delta \partial_{x}^{2}\phi\right)\partial_{\bf n}\phi{\rm d}x.
\end{align}
Note that
\begin{align*}
 &\int_{\Gamma}  \left(\phi_{t} - \delta \partial_{x}^{2}\phi\right) \partial_{\bf n}\phi{\rm d}x 
\nonumber\\ 
 & =  \int_{\Gamma}  \left(\phi_{t} - \delta \partial_{x}^{2}\phi\right) \left( \mathcal{L}(\phi)  - \gamma'_{fs}(\phi) \right) {\rm d}x
 \nonumber\\
 & = - \frac{{\rm d}}{{\rm d}t}   \int_{\Gamma} \gamma_{fs}(\phi)  {\rm d}x - \delta \int_{\Gamma}     \gamma^{(2)}_{fs}(\phi) |\partial_{x} \phi |^{2} {\rm d}x+ \int_{\Gamma}  \left(\phi_{t} - \delta \partial_{x}^{2}\phi\right)  \mathcal{L}(\phi)  {\rm d}x
 \nonumber\\
 & = - \frac{{\rm d}}{{\rm d}t}   \int_{\Gamma} \gamma_{fs}(\phi)  {\rm d}x - \delta \int_{\Gamma}     \gamma^{(2)}_{fs}(\phi) |\partial_{x} \phi |^{2} {\rm d}x + \int_{\Gamma} \left(-\mathcal{L}(\phi) -u_1\partial_{x}\phi  \right)  \mathcal{L}(\phi)  {\rm d}x
 \nonumber\\
 & = - \frac{{\rm d}}{{\rm d}t}   \int_{\Gamma} \gamma_{fs}(\phi)  {\rm d}x - \delta \int_{\Gamma}     \gamma^{(2)}_{fs}(\phi) |\partial_{x} \phi |^{2} {\rm d}x -    \|    \mathcal{L}(\phi)\|_{L^2(\Gamma)}^{2} - \int_{\Gamma}    \mathcal{L}(\phi)  u_1\partial_{x}\phi{\rm d}x.
\end{align*}
Hence, we find that $(\ref{A-E-2})$ becomes 
\begin{align*}
&\frac{1}{2}\frac{{\rm d}}{{\rm d}t} \left[ \| \nabla \phi \|_{L^2}^{2}+ \frac{1}{2}\| \phi^2-1 \|_{L^2}^2 + \int_{\Gamma} \gamma_{fs}(\phi)  {\rm d}x+ \delta^{2} \| \partial_{x}\phi \|_{L^2}^{2}\right]   \nonumber \\
&\hspace{2em} +\| \nabla \mu \|_{L^2}^{2} + \delta \| \phi_{t}\|_{L^{2}}^{2}
 + \delta   \| \nabla\partial_{x}\phi\|_{L^2}^{2}  +   \|    \mathcal{L}(\phi)\|_{L^2(\Gamma)}^{2}
 \nonumber \\
& \le  -\int_{\Omega} \mu {\bf u} \cdot \nabla \phi {\rm d}x{\rm d}y  - \int_{\Gamma}    \mathcal{L}(\phi)  u_1\partial_{x}\phi{\rm d}x + \delta \| \partial_{x}\phi\|_{L^2}^{2} + C\delta \| \partial_{x} \phi\|_{L^2(\Gamma)}^{2}
 \nonumber \\
& \le  -\int_{\Omega} \mu {\bf u} \cdot \nabla \phi {\rm d}x{\rm d}y  - \int_{\Gamma}    \mathcal{L}(\phi)  u_1\partial_{x}\phi{\rm d}x + C \| \nabla\phi\|_{L^2}^{2} + C\delta \| \partial_{x} \phi\|_{L^2}\| \partial_{x} \phi\|_{H^1}
 \nonumber \\
& \le  -\int_{\Omega} \mu {\bf u} \cdot \nabla \phi {\rm d}x{\rm d}y  - \int_{\Gamma}    \mathcal{L}(\phi)  u_1\partial_{x}\phi{\rm d}x + C \| \nabla\phi\|_{L^2}^{2} + \frac{ \delta}{2}   \| \nabla\partial_{x}\phi\|_{L^2}^{2}.
\end{align*}

This together with $(\ref{E-1})$ and Gornwall inequality gives
\begin{align*}
\sup\limits_{0\le t\le T} \left( \hat{\mathcal{E}}_{1}(t) + \delta^{2} \| \partial_{x}\phi \|_{L^2}^{2} \right )+\int_{0}^{T} \left( \hat{\mathcal{D}}(t) + \delta \| \phi_{t}\|_{L^{2}}^{2} + \delta   \| \nabla\partial_{x}\phi\|_{L^2}^{2}  \right) {\rm d}t \le C\hat{\mathcal{E}}_{1}(0).
\end{align*}
It follows from $(\ref{A-NSCH})_{3}$ that
\begin{align}\label{000}
 \frac{\rm d}{{\rm d}t} \langle \phi   \rangle & = \frac{1}{|\Omega|}\int_\Omega \phi_t   \mathrm{d}x {\rm d}y  = \frac{1}{|\Omega|} \int_\Omega \Delta\mu   \mathrm{d}x{\rm d}y  - \frac{1}{|\Omega|}\int_\Omega {\bf u} \cdot \nabla \phi   \mathrm{d}x {\rm d}y + \delta\frac{1}{|\Omega|}\int_\Omega \partial_x^2 \phi   \mathrm{d}x{\rm d}y
\nonumber \\
& = - \frac{1}{|\Omega|}\int_{\Gamma} \partial_{\bf n} \mu {\rm d}x  + \frac{1}{|\Omega|}\int_\Omega {\rm div} {\bf u} \,  \phi   \mathrm{d}x {\rm d}y -\frac{1}{|\Omega|}\int_{\Gamma}  u_{2} \, \phi {\rm d}x+ 0 = 0.
\end{align}
Thus, using Poincar${\rm \acute{e}}$ inequality, we have
\begin{align}\label{dijie-11}
 \| \phi \|_{L^{2}}^{2}  \le    \left\| \phi -  \langle \phi   \rangle  \right\|_{L^{2}}^{2} + \left\|   \langle \phi   \rangle  \right\|_{L^{2}}^{2} \le    C \left(  \| \nabla\phi\|_{L^{2}}^{2} + \left|\langle \phi_{0} \rangle \right|^{2}  \right) \le C\hat{\mathcal{E}}(0) .
 \end{align}
Finally, we arrive at
\begin{align}\label{d-JL-1}
\sup\limits_{0\le t\le T} \hat{\mathcal{E}}(t) +\int_{0}^{T} \left( \hat{\mathcal{D}}(t) + \delta \| \phi_{t}\|_{L^{2}}^{2} + \delta   \| \nabla\partial_{x}\phi\|_{L^2}^{2}  \right) {\rm d}t \le C\hat{\mathcal{E}}(0).
\end{align}
\ \\[-0.5em]

{\bf Step 3. Estimates of  ${\boldsymbol \| \phi_{t} \|_{L^{\infty}(0,T;H^{1})} }$.}

Differentiating $(\ref{A-NSCH})_{3,4}$ and $(\ref{NN-BC})$ with respect to $t$ respectively, we get
\begin{align}\label{d-1-t}
\begin{cases}
\Delta\mu_{t}  = \phi_{tt} - \delta \partial_x^2 \phi_{t}+ {\bf u}_{t} \cdot \nabla\phi  + {\bf u} \cdot \nabla\phi_{t}  , & {\rm in}~~ \Omega,
\\
 \Delta\phi_{t} = \delta\phi_{tt} -\mu_{t}  +f_{t}, & {\rm in}~~ \Omega,
\\
\Delta\mu_{t} = - \partial_{\bf n}\phi_{t} - \gamma_{fs}^{(2)}(\phi)\phi_{t},  & {\rm on}~~ \Gamma.
\end{cases}
\end{align}
We now consider
\begin{align}\label{d-E-1-1}
-\langle \Delta\mu_{t},  \Delta\phi_{t}\rangle =   \langle  \nabla\Delta\mu_{t} ,  \nabla\phi_{t}  \rangle -   \int_{\Gamma}  \Delta\mu_{t}   \partial_{\bf n}\phi_{t}  {\rm d}x =: I_{1}^{\delta} + I_{2}^{\delta}.
\end{align}
Recall that $ I_{2}^{\delta} = I_{2}$, where $ I_{2}$ is defined as in $(\ref{E-1-1})$, we find
\begin{align}\label{d-I-1}
I_{1}^{\delta} = I_{1} - \delta  \int_{\Omega}  \nabla \partial_{x}^{2} \phi_{t} \cdot  \nabla\phi_{t}  {\rm d}x{\rm d}y =  I_{1} +  \delta \| \nabla \partial_{x} \phi_{t}\|_{L^{2}}^{2},
\end{align}
where $ I_{1}$ is defined as in $(\ref{E-1-0})$. Putting $(\ref{d-1-t})_{2}$ into $(\ref{d-E-1-1})$, we get 
\begin{align}\label{d-E-1-2}
-\langle \Delta\mu_{t},  \Delta\phi_{t}  \rangle & = \int_{\Omega} \Delta\mu_{t}  \left( -\delta \phi_{tt} +\mu_{t}  -f_{t} \right)  {\rm d}x{\rm d}y
\nonumber \\
& =   -\| \nabla \mu_{t} \|_{L^{2}}^{2} +\int_{\Omega} \nabla\mu_{t} \cdot \nabla f_{t}  {\rm d}x{\rm d}y - \delta \int_{\Omega} \Delta\mu_{t}   \phi_{tt}   {\rm d}x{\rm d}y.
\end{align}
Compared with $(\ref{E-1-2})$,  we only need to consider the last term on right hand side of $(\ref{d-E-1-2})$. Indeed, we find that
\begin{align}\label{Diff-t}
-\delta \int_{\Omega} \Delta\mu_{t}   \phi_{tt}   {\rm d}x{\rm d}y & = -\delta \int_{\Omega} \left(  \phi_{tt} - \delta \partial_x^2 \phi_{t}+ {\bf u}_{t} \cdot \nabla\phi  + {\bf u} \cdot \nabla\phi_{t} \right) \phi_{tt}   {\rm d}x{\rm d}y
\nonumber \\
& = -\delta \| \phi_{tt} \|_{L^{2}}^{2} - \frac{\delta^{2}}{2}\frac{\rm d}{{\rm d}t} \| \partial_{x}\phi_{t} \|_{L^{2}}^{2}- \delta \int_{\Omega} \left(   {\bf u}_{t} \cdot \nabla\phi  + {\bf u} \cdot \nabla\phi_{t} \right) \phi_{tt}   {\rm d}x{\rm d}y.
\end{align}
Thus, substituting $(\ref{d-I-1})$, $(\ref{d-E-1-2})$ and $(\ref{Diff-t})$ into $(\ref{d-E-1-1})$,  and following the similar arguments to $(\ref{C-1})$, we have
\begin{align}\label{d-C-1}
& \frac{1}{2} \frac{\rm d}{{\rm d}t} \left(  \| \nabla \phi_{t}\|_{L^{2}}^{2} +  \delta^{2}\| \partial_{x}\phi_{t} \|_{L^{2}}^{2} \right) +  \| \nabla \mu_{t} \|_{L^{2}}^{2} + \| \partial_{\bf n}\phi_{t}\|_{L^{2}(\Gamma)}^{2} + \delta \| \nabla \partial_{x} \phi_{t}\|_{L^{2}}^{2} + \delta \| \phi_{tt} \|_{L^{2}}^{2}
\nonumber\\
& \le \epsilon \left(  2 \| \nabla {\bf u}_{t} \|_{L^{2}}^{2} + \| \nabla {\bf u} \|_{H^{2}}^{2} \right)  +  \frac{1}{4} \| \nabla \mu_{t} \|_{L^{2}}^{2} +  \frac{1}{6} \| \partial_{\bf n}\phi_{t}\|_{L^{2}(\Gamma)}^{2}  + C\mathcal{E}(t)^{3}
\nonumber\\
&~\,\,\,-\delta \int_{\Omega} \left(   {\bf u}_{t} \cdot \nabla\phi  + {\bf u} \cdot \nabla\phi_{t} \right) \phi_{tt}   {\rm d}x{\rm d}y
\nonumber\\
&  \le \epsilon \left(  2 \| \nabla {\bf u}_{t} \|_{L^{2}}^{2} + \| \nabla {\bf u} \|_{H^{2}}^{2} \right)  +  \frac{1}{4} \| \nabla \mu_{t} \|_{L^{2}}^{2} +  \frac{1}{6} \| \partial_{\bf n}\phi_{t}\|_{L^{2}(\Gamma)}^{2}  + C\mathcal{E}(t)^{3}
\nonumber\\
& ~\,\,\, +  \frac{\delta}{2}\| \phi_{tt} \|_{L^{2}}^{2} + \delta \left( \|  {\bf u}_{t} \|_{L^{2}}^{2} \|  \nabla\phi\|_{L^{\infty}}^{2} + \|  {\bf u} \|_{L^{\infty}}^{2} \|  \nabla\phi_{t}\|_{L^{2}}^{2} \right)
\nonumber\\
&  \le \epsilon \left(  2 \| \nabla {\bf u}_{t} \|_{L^{2}}^{2} + \| \nabla {\bf u} \|_{H^{2}}^{2} \right)  +  \frac{1}{4} \| \nabla \mu_{t} \|_{L^{2}}^{2} +  \frac{1}{6} \| \partial_{\bf n}\phi_{t}\|_{L^{2}(\Gamma)}^{2}  + C\mathcal{E}(t)^{3} +  \frac{\delta}{2}\| \phi_{tt} \|_{L^{2}}^{2} .
\end{align}
Operating $\partial_{x}$ to $(\ref{A-BC})_{4}$, and testing the result by $\partial_{x}\phi_{t}$, then we find that $(\ref{C-yb})$ is replaced by the following inequality
\begin{align}\label{d-C-yb}
\frac{\delta}{2} \frac{\rm d}{{\rm d}t}  \|\partial_{x}^{2}\phi\|_{L^{2}(\Gamma)}^{2}  + \|\partial_{x}\phi_{t}\|_{L^{2}(\Gamma)}^{2}  \le  C\mathcal{E}(t)^{2}.
\end{align}
Putting $(\ref{d-C-1})$,  $(\ref{d-C-yb})$, $(\ref{C-2})$,  $(\ref{u-H33})$ and $(\ref{C-3})$ together, using Lemma $\ref{uuu}$ and  Korn's inequality, and choosing $\epsilon$ sufficiently small such that $C\epsilon<\frac{1}{12}$, we finally arrive at
\begin{align}\label{d-e-1}
 &\sup\limits_{0\le t \le T} \left( \| \nabla \phi_{t}\|_{L^{2}}^{2}  +  \|{\bf u}_t  \|^2_{L^2} + \| {\bf u}   \|^2_{H^2}  +   \delta^{2}\| \partial_{x}\phi_{t} \|_{L^{2}}^{2} + \delta\|\partial_{x}^{2}\phi\|_{L^{2}(\Gamma)}^{2} \right)
 \nonumber \\
& ~\,\,\,+ \int_{0}^{T}\left( \mathcal{D}_{1}(t) + \delta \| \nabla \partial_{x} \phi_{t}\|_{L^{2}}^{2} + \delta \| \phi_{tt} \|_{L^{2}}^{2}
 \right) {\rm d}t
\nonumber \\
&\le C \left(\| {\bf u}_{0}\|_{H^{2}}, \| \phi_{0}\|_{H^{3}},  \| \psi_{0}\|_{H^{3}} , \| \mu_{0}^{\delta}\|_{H^{3}}\right)\int_{0}^{T} \mathcal{E}(t)^{4} {\rm d}t.
\end{align}
Noting $(\ref{000})$ and using Poincar${\rm \acute{e}}$ inequality, we have
\begin{align}\label{d-dijie}
 \| \phi_{t}\|_{L^{2}}^{2}   = \left\| \phi_{t} -  \langle \phi_{t}   \rangle  \right\|_{L^{2}}^{2}  & \le C  \| \nabla\phi_{t}\|_{L^{2}}^{2} .
 \end{align}
\ \\[-2em]

{\bf Step 4. Estimates of $\|\mu\|_{L^{\infty}(0,T;L^{2})\cap L^{2}(0,T;H^{3})}$.}

Note $\langle \phi_{t} \rangle = 0$ (see $(\ref{000})$). Thus, we find that $(\ref{p-mu})$-$(\ref{mu-0})$ also hold. Applying the elliptic estimates for $(\ref{A-NSCH})_{3}$ with boundary condition $(\ref{A-BC})_{2}$ gives that
\begin{align}\label{mu-H3}
\int_{0}^{T }\| \mu     \|^2_{H^3} {\rm d}t  & \le C \int_{0}^{T } \left( \| \mu  \|^2_{L^2} + \|  \phi  _t \|^2_{H^1}  + \|\partial_x^2\phi \|^2_{H^1} + \| {\bf u} \|^2_{L^\infty}\| \nabla \phi   \|^2_{L^2} \right)  {\rm d}t
\nonumber \\
& ~\,\,\, + C \int_{0}^{T } \left( \| \nabla {\bf u} \|^2_{L^4}\| \nabla \phi   \|^2_{L^4} + \| {\bf u} \|^2_{L^\infty}\| \nabla^{2} \phi   \|^2_{L^2} \right) {\rm d}t
\nonumber \\
& \le C   \int_{0}^{T} \mathcal{E}(t)^{4} {\rm d}t + C,
\end{align}
where we have used $(\ref{mu-0})$.
\ \\[-1em]

{\bf Step 5. Estimates of $\| \phi \|_{L^\infty(0,T;H^{3})}$ and $\| \mu\|_{L^\infty(0,T;H^{3})}$.}

For $k=1,2$, operating $\nabla\partial_{x}^{k}$, $\partial_{x}^{k}$ and $\partial_{x}^{k}$  to $(\ref{A-NSCH})_{3}$, $(\ref{A-NSCH})_{4}$ and $(\ref{NN-BC})$ respectively, we get
\begin{align}\label{d-2-x}
\begin{cases}
\nabla\partial_{x}^{k}\Delta\mu  = \nabla\partial_{x}^{k}\phi_{t} - \delta\nabla\partial_{x}^{k}\partial_x^2 \phi + [\nabla\partial_{x}^{k}, {\bf u} \cdot \nabla] \phi+{\bf u} \cdot \nabla^{2} \partial_{x}^{k} \phi , & {\rm in}~~ \Omega,
\\
\partial_{x}^{k}\Delta\phi =  \delta \partial_{x}^{k} \phi_{t}-\partial_{x}^{k}\mu  + \partial_{x}^{k}f, & {\rm in}~~\Omega,
\\
\partial_{x}^{k}\Delta\mu = - \partial_{x}^{k}\partial_{\bf n}\phi -  \partial_{x}^{k}\gamma_{fs}'(\phi) ,  & {\rm on}~~ \Gamma,
\end{cases}
\end{align}
where $[A,B] = AB-BA$. We now consider
\begin{align}\label{d-E-phi-3}
-\langle \partial_{x}^{k}\Delta\mu , \partial_{x}^{k}  \Delta\phi \rangle & = \langle \nabla \partial_{x}^{k} \Delta\mu,  \partial_{x}^{k} \nabla\phi \rangle  -  \int_{\Gamma} \partial_{x}^{k} \Delta\mu \partial_{x}^{k}\partial_{\bf n}\phi {\rm d}x
 =: J_{1}^{\delta} + J_{2}^{\delta}.
\end{align}
Noting $J_{2}^{\delta} = J_{2}$, where $J_{2}$ is defined as in $(\ref{E-2-1})$, we find
\begin{align}\label{d-E-2-0}
J_{1}^{\delta} & = J_{1} - \delta \int_{\Omega} \nabla\partial_{x}^{k}\partial_x^2 \phi \cdot \nabla\partial_{x}^{k}\phi  {\rm d}x{\rm d}y  = J_{1} + \delta \| \nabla\partial_{x}^{k+1}  \phi\|_{L^{2}},
\end{align}
where $J_{1}$ is defined as in $(\ref{E-2-0})$. We put $(\ref{d-2-x})_{2}$ into $(\ref{d-E-phi-3})$ and get
\begin{align}\label{d-E-2-2}
&-\langle  \partial_{x}^{k} \Delta\mu,  \partial_{x}^{k} \Delta\phi \rangle
 = \int_{\Omega}  \partial_{x}^{k} \Delta\mu   \left( -  \delta \partial_{x}^{k} \phi_{t} + \partial_{x}^{k} \mu  -  \partial_{x}^{k}f \right)  {\rm d}x{\rm d}y
\nonumber \\
&= - \| \nabla \partial_{x}^{k} \mu \|_{L^{2}}^{2} + \int_{\Omega} \partial_{x}^{k}\nabla\mu  \cdot \nabla \partial_{x}^{k} f {\rm d}x{\rm d}y -  \delta\int_{\Omega}  \partial_{x}^{k} \Delta\mu  \partial_{x}^{k} \phi_{t} {\rm d}x{\rm d}y.
\end{align}
Compared with $(\ref{E-2-2})$,  we only need to consider the last term on right hand side of $(\ref{d-E-2-2})$. It holds
\begin{align}\label{Diff-x}
 -\delta\int_{\Omega}  \partial_{x}^{k} \Delta\mu  \partial_{x}^{k} \phi_{t} {\rm d}x{\rm d}y & =  -\delta\int_{\Omega} \left( \partial_{x}^{k} \phi_{t} -\delta \partial_{x}^{k}\partial_x^2 \phi + \partial_{x}^{k} ({\bf u} \cdot \nabla \phi)\right) \partial_{x}^{k} \phi_{t} {\rm d}x{\rm d}y
 \nonumber\\
 & = -\delta \| \partial_{x}^{k} \phi_{t} \|_{L^{2}}^{2} - \frac{\delta^{2}}{2}\frac{\rm d}{{\rm d}t} \| \partial_{x}^{k+1} \phi \|_{L^{2}}^{2} -  \delta\int_{\Omega}  \partial_{x}^{k} ({\bf u} \cdot \nabla \phi) \partial_{x}^{k} \phi_{t} {\rm d}x{\rm d}y.
\end{align}
Thus, substituting $(\ref{d-E-2-0})$, $(\ref{d-E-2-2})$ and $(\ref{Diff-x})$ into $(\ref{d-E-phi-3})$,  and repeating the arguments in $(\ref{e-3})$-$(\ref{hz})$, we have, for $k=1,2$,
\begin{align}\label{d-hz}
 &\frac{1}{2}\frac{\rm d}{{\rm d}t} \left( \|  \nabla\partial_{x}^{k} \phi\|_{L^{2}}^{2} +  \delta^{2}\| \partial_{x}^{k+1} \phi \|_{L^{2}}^{2} \right)
   \nonumber\\
& +  \frac{1}{2} \|  \partial_{x}^{k} \nabla \mu \|_{L^{2}}^{2} +  \frac{1}{2}\| \partial_{x}^{k} \partial_{\bf n}\phi \|_{L^{2}(\Gamma)}^{2} + \delta \| \nabla\partial_{x}^{k+1}  \phi\|_{L^{2}}^{2} + \delta \| \partial_{x}^{k} \phi_{t} \|_{L^{2}}^{2}
 \nonumber\\
& \le  C\| {\bf u} \|_{H^3}^{2} + C \mathcal{E}(t)^{3} -   \delta\int_{\Omega}  \partial_{x}^{k} ({\bf u} \cdot \nabla \phi) \partial_{x}^{k} \phi_{t} {\rm d}x{\rm d}y
 \nonumber\\
& \le  C\| {\bf u} \|_{H^3}^{2} + C \mathcal{E}(t)^{3} + \frac{\delta}{2} \| \partial_{x}^{k} \phi_{t} \|_{L^{2}}^{2} +  C\| \partial_{x}^{k} {\bf u} \|_{L^{2}}^{2}  \| \nabla\phi \|_{L^{\infty}}^{2} + C\|  {\bf u} \|_{L^{\infty}}^{2}  \| \partial_{x}^{k}\nabla\phi \|_{L^{2}}^{2}
\nonumber\\
& \le  C\| {\bf u} \|_{H^3}^{2} + C \mathcal{E}(t)^{3} + \frac{\delta}{2} \| \partial_{x}^{k} \phi_{t} \|_{L^{2}}^{2}.
 \end{align}
Using Gronwall inequality and $(\ref{d-e-1})$, we obtain,  for $k=1,2$,
\begin{align}\label{d-phi-33}
& \sup\limits_{0\le t \le T} \left( \|  \nabla\partial_{x}^{k} \phi\|_{L^{2}}^{2} +  \delta^{2}\| \partial_{x}^{k+1} \phi \|_{L^{2}}^{2} \right)
   \nonumber\\
& +  \int_{0}^{T} \left( \|  \partial_{x}^{k} \nabla \mu \|_{L^{2}}^{2} +  \| \partial_{x}^{k} \partial_{\bf n}\phi \|_{L^{2}(\Gamma)}^{2} + \delta \| \nabla\partial_{x}^{k+1}  \phi\|_{L^{2}} + \delta \| \partial_{x}^{k} \phi_{t} \|_{L^{2}}^{2} \right) {\rm d}t
 \nonumber\\
& \le   C   \int_{0}^{T}  \mathcal{E}(t)^{4} {\rm d}t + C.
 \end{align}
 We will derive the estimate for $\phi$ involving normal derivatives. It follows from $(\ref{A-NSCH})_{4}$ that
\begin{align}\label{d-F-X-1}
 \|  \partial_{y}^{2} \phi\|_{L^{2}}^{2}
 & \le  C \left( \| \partial_{x}^{2}\phi\|_{L^{2}}^{2} + \delta^{2}  \|  \phi_{t}\|_{L^{2}}^{2}  + \|  \mu\|_{L^{2}}^{2} + \|  (\phi^{3}-\phi)\|_{L^{2}}^{2} \right)
 \nonumber \\
 &  \le  C \| \nabla\partial_{x} \phi\|_{L^{2}}^{2} +   C\|  \phi_{t}\|_{L^{2}}^{2}  + C\| \mu\|_{L^{2}}^{2} + C \| \phi \|_{H^{1}}^{3}
 \nonumber \\
 &\le    C  \int_{0}^{T} \mathcal{E}(t)^{4} {\rm d}t + C,
 \end{align}
where we have used $(\ref{d-JL-1} )$, $(\ref{d-e-1})$, $(\ref{d-dijie})$, $(\ref{mu-0})$ and $(\ref{d-phi-33})$. Thus,  $(\ref{d-JL-1} )$ together with $(\ref{d-phi-33})$ and $(\ref{d-F-X-1})$ gives
\begin{align}\label{d-phi-2}
 \| \phi\|_{H^{2}}^{2} \le    C  \int_{0}^{T} \mathcal{E}(t)^{4} {\rm d}t + C.
 \end{align}
Applying the elliptic estimates for $(\ref{A-NSCH})_{3}$ with boundary condition $(\ref{A-BC})_{2}$ gives that
\begin{align}\label{d-mu-H3}
\| \mu     \|^2_{H^3}   & \le C \left( \| \mu  \|^2_{L^2} + \|  \phi  _t \|^2_{H^1}  + \|\partial_x^2\phi \|^2_{H^1} + \| {\bf u} \|^2_{L^\infty}\| \nabla \phi   \|^2_{L^2} \right)
\nonumber \\
& ~\,\,\,+ C \left( \| \nabla{\bf u} \|^2_{L^4} \|  \nabla \phi  \|^2_{L^4}  +    \| {\bf u} \|^2_{L^\infty} \| \nabla^{2} \phi  \|^2_{L^2}  \right)
 \nonumber \\
& \le C   \int_{0}^{T} \mathcal{E}(t)^{8} {\rm d}t + C,
\end{align}
where we have used $(\ref{mu-0})$, $(\ref{d-e-1})$, $(\ref{d-dijie})$, $(\ref{d-phi-33})$ and $(\ref{d-phi-2})$. It follows from $(\ref{A-NSCH})_{4}$ that
\begin{align}\label{d-F-X-2}
 \| \nabla \partial_{y}^{2}\phi\|_{L^{2}}^{2}
 & \le  C \left( \|\nabla\partial_{x}^{2} \phi\|_{L^{2}}^{2} + \delta^{2}  \| \nabla \phi_{t}\|_{L^{2}}^{2}  + \|  \nabla \mu\|_{L^{2}}^{2} + \| \nabla (\phi^{3}-\phi)\|_{L^{2}}^{2} \right)
 \nonumber \\
 &  \le  C \| \nabla\partial_{x}^{2} \phi\|_{L^{2}}^{2} +   C\|  \nabla\phi_{t}\|_{L^{2}}^{2}  + C\| \nabla\mu\|_{L^{2}}^{2} + C \| \phi \|_{H^{2}}^{3}
 \nonumber \\
 &\le    C  \int_{0}^{T} \mathcal{E}(t)^{8} {\rm d}t + C,
 \end{align}
where we have used $(\ref{d-e-1})$, $(\ref{d-dijie})$, $(\ref{d-phi-33})$ and $(\ref{d-phi-2})$. Hence,  $(\ref{d-JL-1} )$ together with $(\ref{d-phi-33})$, $(\ref{d-F-X-2})$ and $(\ref{d-phi-2})$ gives
\begin{align}\label{d-phi-3}
 \| \phi\|_{H^{3}}^{2} \le    C  \int_{0}^{T} \mathcal{E}(t)^{8} {\rm d}t + C.
 \end{align}
\ \\[-1em]

{\bf Step 6. Estimates of $\| \Delta\mu \|_{L^{2}(0,T;H^{2})}$.}

Operating $\nabla\partial_{x}$, $\partial_{t}\partial_{x}$ and $\partial_{x}$  to $(\ref{A-NSCH})_{3}$, $(\ref{A-NSCH})_{4}$ and $(\ref{NN-BC})$ respectively, we get
\begin{align}\label{d-2-x-t}
\begin{cases}
  \nabla\partial_{x}\Delta\mu    = \nabla\partial_{x}\phi_{t} - \delta \nabla\partial_{x}\partial_x^2\phi  + \nabla\partial_{x} {\bf u} \cdot  \nabla\phi +\partial_{x} {\bf u} \cdot \nabla^{2} \phi
\\
 \hspace{1.6cm}+  \nabla {\bf u} \cdot \nabla \partial_{x} \phi + {\bf u} \cdot \nabla^{2} \partial_{x} \phi ,& {\rm in}~ \Omega,
\\
\partial_{x} \Delta\phi_{t} =  \delta\partial_{x}\phi_{tt}-\partial_{x} \mu_{t}  + (3\phi^{2}-1)\partial_{x} \phi_{t} +6\phi\phi_{t}\partial_{x} \phi, & {\rm in}~ \Omega,
\\
\partial_{x} \Delta\mu = - \partial_{x} \partial_{\bf n}\phi -   \gamma_{fs}^{(2)}(\phi) \partial_{x} \phi , &{\rm on}~ \Gamma.
\end{cases}
\end{align}
We now consider
\begin{align}\label{d-E-mu-4}
-\langle  \partial_{x} \Delta\mu,   \partial_{x}  \Delta\phi_{t}\rangle & =   \langle \nabla \partial_{x}  \Delta\mu  , \partial_{x} \nabla\phi_{t}\rangle-  \int_{\Gamma}   \partial_{x} \Delta\mu \partial_{x} \partial_{\bf n}\phi_{t} {\rm d}x
 =: N_{1}^{\delta} + N_{2}^{\delta}.
\end{align}
Recall that $N_{2}^{\delta} = N_{2}$, where $N_{2}$ is defined as in $(\ref{N-3})$, and then we find
\begin{align}\label{d-N-1}
 N_{1}^{\delta} & = N_{1} - \delta \int_{\Omega} \nabla\partial_{x}\partial_x^2\phi  \cdot  \nabla\partial_{x}\phi_{t}  {\rm d}x{\rm d}y =  N_{1} +\frac{\delta}{2}\frac{\rm d}{{\rm d}t}\| \nabla\partial_{x}^{2}\phi \|_{L^{2}}^{2},
 \end{align}
where $N_{1}$ is defined as in $(\ref{N-1})$. We put $(\ref{d-2-x-t})_{2}$ into the left side of $(\ref{d-E-mu-4})$ and get
\begin{align}\label{d-E-mu}
-\langle \partial_{x} \Delta\mu, \partial_{x}  \Delta\phi_{t}  \rangle
& =  \int_{\Omega}  \partial_{x} \Delta\mu   \left(- \delta\partial_{x}\phi_{tt}  +\partial_{x} \mu_{t}  - (3\phi^{2}-1)\partial_{x} \phi_{t} -6\phi\phi_{t}\partial_{x} \phi \right)  {\rm d}x{\rm d}y
\nonumber \\
& =  - \frac{1}{2}\frac{\rm d}{{\rm d}t} \|  \partial_{x} \nabla \mu \|_{L^{2}}^{2} -\int_{\Omega}  \partial_{x} \Delta\mu  \left(    (3\phi^{2}-1)\partial_{x} \phi_{t} +6\phi\phi_{t}\partial_{x} \phi \right)  {\rm d}x{\rm d}y 
\nonumber \\
&~\,\,\,- \delta\int_{\Omega}  \partial_{x} \Delta\mu    \partial_{x}\phi_{tt}   {\rm d}x{\rm d}y.
\end{align}
Compared with $(\ref{E-mu})$,  we only need to consider the last term on right side of $(\ref{d-E-mu})$. Indeed, it holds that
\begin{align}\label{dd-E-mu}
-\delta\int_{\Omega}  \partial_{x} \Delta\mu    \partial_{x}\phi_{tt}   {\rm d}x{\rm d}y
&=  -\delta\int_{\Omega}  \partial_{x} \Delta\mu   \left[\partial_{x} \Delta\mu_{t} + \delta \partial_{x}^{3}\phi_{t}  + \partial_{x} ({\bf u}\cdot \nabla \phi)_{t}\right]    {\rm d}x{\rm d}y
 \nonumber \\
 & = -\frac{\delta }{2}\frac{\rm d}{{\rm d}t}  \|\partial_{x} \Delta\mu  \|_{L^{2}}^{2} +\delta\int_{\Omega} \partial_x^2 \Delta\mu   \left[ \delta\partial_x^2\phi_{t}  +  ({\bf u}\cdot \nabla \phi)_{t}\right]    {\rm d}x{\rm d}y.
\end{align}
Thus, substituting $(\ref{d-N-1})$-$(\ref{dd-E-mu})$ into $(\ref{d-E-mu-4}
)$,  and following the similar arguments as in $(\ref{mu-high})$, we have
\begin{align}\label{d-mu-high}
&\frac{1}{2}\frac{\rm d}{{\rm d}t} \left( \|  \partial_{x} \partial_{\bf n}\phi \|_{L^{2}(\Gamma)}^{2} + \| \partial_{x} \nabla \mu \|_{L^{2}}^{2}  + \delta \| \nabla\partial_{x}^{2}\phi \|_{L^{2}}^{2} + \delta \|\partial_{x} \Delta\mu  \|_{L^{2}}^{2}\right) + \| \nabla\partial_{x}\phi_{t} \|_{L^{2}}^{2}
\nonumber \\
& \le \frac{1}{2}  \|\nabla\partial_{x}\phi_{t} \|_{L^{2}}^{2} +  C\mathcal{E}(t)^{3} + C\| \partial_{\bf n}\phi_{t} \|_{L^{2}(\Gamma)}^{2} + \delta\int_{\Omega} \partial_x^2 \Delta\mu   \left[ \delta\partial_x^2\phi_{t}  +  ({\bf u}\cdot \nabla \phi)_{t}\right]    {\rm d}x{\rm d}y
\nonumber \\
& \le \frac{1}{2}  \|\nabla\partial_{x}\phi_{t} \|_{L^{2}}^{2} +  C\mathcal{E}(t)^{3} + C\| \partial_{\bf n}\phi_{t} \|_{L^{2}(\Gamma)}^{2} + \| \delta\partial_x^2 \Delta\mu  \|_{L^{2}}^{2}  + \|  \delta\partial_x^2\phi_{t} \|_{L^{2}}^{2}
\nonumber \\
& ~\,\,\,+  \| {\bf u}_{t} \|_{L^{2}}^{2} \| \nabla \phi \|_{L^{\infty}}^{2} +  \| {\bf u}\|_{L^{\infty}}^{2} \| \nabla \phi_{t} \|_{L^{2}}^{2}
\nonumber \\
& \le \frac{1}{2}  \|\nabla\partial_{x}\phi_{t} \|_{L^{2}}^{2} +  C\mathcal{E}(t)^{3} + C\| \partial_{\bf n}\phi_{t} \|_{L^{2}(\Gamma)}^{2} + \| \delta\partial_x^2 \Delta\mu  \|_{L^{2}}^{2}  + \delta\| \partial_x^2\phi_{t} \|_{L^{2}}^{2}.
\end{align}
It follows from $(\ref{A-NSCH})_{3}$ that
\begin{align}\label{dd-mu-4}
\| \delta\partial_x^2 \Delta\mu  \|_{L^{2}}^{2} & \le C\left( \| \delta\partial_x^2 \phi_{t}  \|_{L^{2}}^{2} +  \| \delta^{2}\partial_x^{4} \phi  \|_{L^{2}}^{2} +   \| \delta\partial_x^2({\bf u}\cdot \nabla\phi)  \|_{L^{2}}^{2}  \right)
\nonumber \\
& \le C \left( \delta \| \partial_x^2 \phi_{t}  \|_{L^{2}}^{2} + \delta \| \nabla \partial_{x}^{3} \phi  \|_{L^{2}}^{2}  \right) + C\mathcal{E}(t)^{2}.
\end{align}
Putting $(\ref{dd-mu-4})$ into $(\ref{d-mu-high})$, and using Gronwall inequality, $(\ref{d-e-1})$ and $(\ref{d-phi-33})$ for $k=2$, we get
\begin{align}\label{d-phi-t-20}
& \sup\limits_{0\le t\le T} \left( \|  \partial_{x} \partial_{\bf n}\phi \|_{L^{2}(\Gamma)}^{2} + \| \partial_{x} \nabla \mu \|_{L^{2}}^{2}  + \delta \| \nabla\partial_{x}^{2}\phi \|_{L^{2}}^{2} + \delta \|\partial_{x} \Delta\mu  \|_{L^{2}}^{2}\right) + \int_{0}^{T}  \| \nabla\partial_{x}\phi_{t} \|_{L^{2}}^{2} {\rm d}t
\nonumber \\
& \le C  \int_{0}^{T} \mathcal{E}(t)^{8} {\rm d}t + C.
\end{align}
Due to $(\ref{A-NSCH})_{4}$, $(\ref{d-JL-1})$, $(\ref{d-e-1})$ and $(\ref{d-phi-t-2})$, one has
\begin{align*}
&\int_{0}^{T}  \|  \partial_{y}^{2} \phi_{t}\|_{L^{2}}^{2}   {\rm d}t
\nonumber \\
&\le  C \int_{0}^{T} \left( \|\partial_x^2 \phi_{t}\|_{L^{2}}^{2}  + \delta  \| \phi_{tt}\|_{L^{2}}^{2}  + \| \mu_{t}\|_{L^{2}}^{2} + \| (3\phi^{2}-\phi)\phi_{t}\|_{L^{2}}^{2} \right)   {\rm d}t
\nonumber \\
&\le  C \int_{0}^{T} \left( \| \nabla\partial_{x} \phi_{t}\|_{L^{2}}^{2} + \delta  \| \phi_{tt}\|_{L^{2}}^{2} + \|  \mu_{t}\|_{L^{2}}^{2} + (\| \phi\|_{L^{8}}^{4}+1) \|\phi_{t}\|_{L^{4}}^{2} \right)   {\rm d}t
\nonumber \\
&\le  C \int_{0}^{T} \left( \| \nabla\partial_{x} \phi_{t}\|_{L^{2}}^{2} +\delta  \| \phi_{tt}\|_{L^{2}}^{2} + \|  \mu_{t}\|_{L^{2}}^{2} + (\| \phi\|_{H^{1}}^{4}+1) \|\phi_{t}\|_{H^{1}}^{2} \right)   {\rm d}t
\nonumber \\
&\le  C \int_{0}^{T} \left( \| \nabla\partial_{x} \phi_{t}\|_{L^{2}}^{2} + \delta  \| \phi_{tt}\|_{L^{2}}^{2}+ \|  \mu_{t}\|_{L^{2}}^{2}\right)   {\rm d}t  + C \sup\limits_{0\le t\le T} \left( \| \phi\|_{H^{1}}^{4}+1 \right) \int_{0}^{T}   \|\phi_{t}\|_{H^{1}}^{2}    {\rm d}t
\nonumber \\
 & \le  C  \int_{0}^{T} \mathcal{E}(t)^{4} {\rm d}t  + C,
 \end{align*}
which together with $(\ref{d-phi-t-20})$ implies that
\begin{align}\label{d-phi-t-2}
\int_{0}^{T}  \|  \nabla^{2} \phi_{t}\|_{L^{2}}^{2}   {\rm d}t  \le  C  \int_{0}^{T} \mathcal{E}(t)^{4} {\rm d}t  + C.
 \end{align}
This combining with  $(\ref{A-NSCH})_{4}$, \eqref{d-mu-H3} and \eqref{d-phi-3} gives 
\begin{align}\label{d-phi-4}
\int_{0}^{T}  \|  \nabla^{2} \Delta \phi\|_{L^{2}}^{2}   {\rm d}t  \le  C  \int_{0}^{T} \mathcal{E}(t)^{8} {\rm d}t  + C.
 \end{align}
Finally, it follows from $(\ref{A-NSCH})_{4}$, $(\ref{d-phi-t-2})$ and  $(\ref{d-phi-33})$ for $k=2$ that
\begin{align}\label{d-phi-4-f}
& \int_{0}^{T} \delta\| \nabla\partial_{x}  \partial_{y}^{2} \phi\|_{L^2}^{2}  {\rm d}t
\nonumber \\
& \le C \int_{0}^{T} \left( \delta\| \nabla\partial_{x}  \partial_x^2 \phi\|_{L^2}^{2}  + \delta  \| \nabla^{2} \phi_{t}\|_{L^2}^{2}+  \| \nabla^{2} \mu\|_{L^2}^{2} +  \| \nabla^{2}(\phi^{3}-\phi)\|_{L^2}^{2} \right)  {\rm d}t
\nonumber \\
& \le C \int_{0}^{T} \left( \delta\| \nabla \partial_{x} ^{3} \phi\|_{L^2}^{2}  + \delta  \| \nabla^{2} \phi_{t}\|_{L^2}^{2}+  \| \nabla^{2} \mu\|_{L^2}^{2} +  \| (3\phi^{2}-1) \nabla^{2}\phi\|_{L^2}^{2} + \| 6\phi |\nabla\phi|^{2}\|_{L^2}^{2}  \right)  {\rm d}t
\nonumber \\
& \le  C \int_{0}^{T} \mathcal{E}(t)^{4} {\rm d}t + C.
\end{align}
Thus, this together with $(\ref{d-phi-33})$ for $k=2$ gives that
\begin{align}\label{d-phi-4-h}
&\int_{0}^{T} \delta \| \nabla^{2}   \partial_x^2  \phi\|_{L^2}^{2}  {\rm d}t
 \le C \int_{0}^{T}\delta \| \nabla   \partial_{x}^{3}  \phi\|_{L^2}^{2}  {\rm d}t  + C\int_{0}^{T} \delta\| \nabla    \partial_{x}   \partial_{ y}^{2}  \phi\|_{L^2}^{2}  {\rm d}t \le  C \int_{0}^{T} \mathcal{E}(t)^{4} {\rm d}t + C.
\end{align}
 It follows from $(\ref{A-NSCH})_{3}$, $(\ref{d-phi-t-2})$ and $(\ref{d-phi-4-h})$ that
\begin{align}\label{d-mu-40}
\int_{0}^{T} \| \nabla^{2}\Delta \mu \|_{L^2}^{2}  {\rm d}t \le  C \int_{0}^{T} \mathcal{E}(t)^{4} {\rm d}t + C.
\end{align}

{\bf Step 7. Conclusion.}

Putting the above estimates in Step 1-Step 6 together and using Proposition $\ref{P-mu}$, we find that
\begin{align}
\sup\limits_{0\le t\le T}\mathcal{E}^{\delta}(t) \le C_{0} + C_{1}T \left(\sup\limits_{0\le t\le T}\mathcal{E}^{\delta}(t) \right)^{8},
\end{align}
where $C_{0}$ and $C_{1}$ are constant that depend on $\Omega$, $\beta$, $ \hat{\mathcal{E}}(0)$ and the initial value  $\| {\bf u}_{0}\|_{H^{2}}, \| \phi_{0}\|_{H^{3}}$, $\| \partial_x^2\phi_{0}\|_{H^{2}}$ and $ \| \mu_{0}\|_{H^{3}}$, but not on $\delta$. Just as in Section 9 of \cite{C-S}, this provides us with a time of existence $T_{0}$ independent of $\delta$ and an estimate on $[0,T_{0}]$ independent of $\delta$ of the type:
\begin{align*}
\sup\limits_{0\le t\le T_{0}}\mathcal{E}^{\delta}(t)  \le  2C_{0}.
\end{align*}
The proof of Theorem $\ref{P-E-d}$ is completed.
\end{proof}

\section{Proof of Theorem \ref{TH}: Local well-posedness of $(\ref{111111})$, $(\ref{BC})$, $(\ref{initial})$} \label{sec:55}
In this section, we are ready to prove Theorem \ref{TH}.
\begin{proof}[\bf Proof of Theorem ${\bf \ref{TH}}$] It follows from Lemma \ref{Trace} that
\begin{align}\label{tr-phi-0}
\| {\bf tr}(\phi_0) \|_{H^3(\Gamma)}^2 \le C\|\phi_0\|^2_{H^3} + C\|\partial_x^2\phi_0 \|_{H^{2}}^2.
\end{align}
Thus, it is easy to verify that the conditions of Theorem $\ref{d-TH}$ are satisfied by using \eqref{tr-phi-0} and Proposition $\ref{P-mu}$. Now, for each $\delta>0$, we denote the solutions to the $\delta$-approximate problem $(\ref{A-NSCH})$-$(\ref{A-IC})$ by $({\bf u}^{\delta}, \phi^{\delta})$, which were constructed in Theorem $\ref{d-TH}$. With $\delta$-independent estimates $(\ref{d-i-E})$ at hand, it then follows from a standard procedure that as sequence of $\delta\to 0$, the solutions of $(\ref{A-NSCH})$-$(\ref{A-IC})$ approaches a limit which  is a strong solutions of $(\ref{111111})$, $(\ref{BC})$, \eqref{initial}. This shows the existence of solutions to $(\ref{111111})$, $(\ref{BC})$, \eqref{initial}.

It remains to prove the uniqueness. To do this, let  $ ({\bf u}_1,\phi_1) $ and $ ({\bf u}_2,\phi_2) $ be two solutions of the original problem $(\ref{111111})$, $(\ref{BC})$, \eqref{initial}, and let
$$
({\bf u},\phi)=({\bf u}_1-{\bf u}_2,\phi_1-\phi_2) = (({u}^{1}_1-{u}^{2}_1, {u}^{1}_2-{u}^{2}_2),\phi_1-\phi_2) = (({u}_1, {u}_2),\phi_1-\phi_2) .
$$
Then $({\bf u},\phi)$ and the corresponding $\mu, \nabla p$ satisfy the following equations
\begin{align}\label{OEq}
\begin{cases}
{\bf u}_t-{\rm div}\mathbb{S}({\bf u})+ \nabla p = -({\bf u}\cdot\nabla){\bf u}_1-({\bf u}_2\cdot\nabla){\bf u}  +\mu\nabla \phi_1 + \mu_2\nabla \phi, &~~\text{in}~~\Omega,
\\
{\rm div}{\bf u}=0, &~~\text{in}~~\Omega,
\\
\phi_ t + {\bf u} \cdot \nabla \phi_1  + {\bf u}_2 \cdot \nabla \phi = \Delta\mu   , &~~\text{in}~~\Omega,
\\
\mu=-\Delta\phi+\phi\left( \phi_1^2+\phi_1\phi_2+\phi_2^2 -1\right), &~~\text{in}~~\Omega,
\end{cases}
\end{align}
with boundary condition
\begin{align}\label{OBq}
\begin{cases}
u_{2}= 0, &~~\text{on}~~\Gamma,
\\
\partial_{\bf n}\mu =0, &~~\text{on}~~\Gamma,
\\
\beta u_1+\partial_{\bf n} u_1= g, &~~\text{on}~~\Gamma,
\\
\phi_t + u_{1} \partial_{x} \phi_1 + u^{2}_{1}\partial_{x} \phi  = -\partial_{\bf n}\phi - \gamma^{(2)}_{fs}(\theta\phi_{1}+(1-\theta)\phi_{2})\phi , &~~\text{on}~~\Gamma,
\end{cases}
\end{align}
where
\begin{align}
& g = \partial_{\bf n}\phi \partial_{x} \phi_1 -  \partial_{\bf n}\phi_2\partial_{x} \phi + \gamma^{(2)}_{fs}(\theta\phi_{1}+(1-\theta)\phi_{2})\phi \partial_{x} \phi_{1} + \gamma'_{fs}(\phi_{2})\partial_{x} \phi, ~\theta \in (0,1), \nonumber
\end{align}
and initial condition
\begin{equation}\label{chuzhi}
({\bf u},\phi)\big|_{t=0}=(0,0), ~~\text{\rm in}~\Omega.
\end{equation}

Applying \eqref{Product-2} with $r=s_{1}=\frac{1}{2} \le s_{2} = \frac{3}{2} $  and  $ s_{2} > r + \frac{N}{2} =  \frac{1}{2} + \frac{1}{2} =1$, one has
\begin{align}\label{h-r-1}
 \|    \partial_{\bf n}\phi_{2} \, u_1 \|_{H^{\frac12}(\Gamma)} \le C\|  \partial_{\bf n}\phi_{2}\|_{H^{\frac32}(\Gamma)} \| u_1 \|_{H^{\frac12}(\Gamma)}.
\end{align}
Multiplying $(\ref{OEq})_1$ by $ {\bf u}$, and integration by parts over $\Omega$, we obtain that
\begin{align}\label{weiyi-1}
&\dfrac{1}{2} \dfrac{\rm d}{{\rm d}t}  \|{\bf u} \|^2_{L^2}  + \dfrac{1}{2} \|\mathbb{S}({\bf u}) \|^2_{L^2}  + \beta\| u_1\|^2_{L^2(\Gamma)}
\nonumber \\
&= \int_\Omega  \left( -{\bf u}\cdot\nabla{\bf u}_1-{\bf u}_2\cdot\nabla{\bf u}  - \nabla \mu \phi_1 + \mu_2\nabla \phi \right) \cdot {\bf u} \,{\rm d}x{\rm d}y + \int_{\Gamma} g \cdot  u_1 {\rm d}x
\nonumber \\
& \le    \| \nabla {\bf u}_{1}\|_{L^{\infty}} \| {\bf u}\|_{L^{2}}^{2} +  \|  {\bf u}_{2}\|_{L^{\infty}} \| {\bf u}\|_{L^{2}} \| \nabla{\bf u}\|_{L^{2}} + \|  \nabla \mu \|_{L^{2}} \| \phi_{1}\|_{L^{\infty}} \| {\bf u}\|_{L^{2}} +  \|  \mu_{2} \|_{L^{\infty}} \| \nabla\phi\|_{L^{2}} \| {\bf u}\|_{L^{2}}
\nonumber \\
&~\,\,\,+\Big( \|  \partial_{\bf n}\phi \|_{L^{2}(\Gamma)} \|  \partial_{x} \phi_1  \|_{L^{\infty}(\Gamma)}  + \|\phi \|_{L^{2}(\Gamma)} \|  \partial_{x} \phi_1  \|_{L^{\infty}(\Gamma)} \Big)  \|   u_1 \|_{L^{2}(\Gamma)}
\nonumber \\
&~\,\,\, + \|  \partial_{x} \phi \|_{H^{-\frac12}(\Gamma)} \|    \partial_{\bf n}\phi_{2} \, u_1 \|_{H^{\frac12}(\Gamma)} + \| \partial_{x} \phi \|_{H^{-\frac12}(\Gamma)} \|    u_1 \|_{H^{\frac12}(\Gamma)}
\nonumber \\
& \le \epsilon   \| \nabla {\bf u} \|_{L^{2}} +  \frac{1}{4}  \| \nabla \mu \|_{L^{2}} +   \frac{1}{4}\| \partial_{\bf n} \phi \|_{L^{2}(\Gamma)}^{2} +  C \left( \| \nabla {\bf u}_{1}\|_{L^{\infty}} + \| {\bf u}_{2}\|_{H^{2}}^{2} +  \| \phi_{1}\|_{H^{2}}^{2} +  \| \mu_{2}\|_{H^{2}}^{2} \right)  \| {\bf u}\|_{L^{2}}^{2}
\nonumber \\
&~\,\,\,+ C\| \nabla\phi\|_{L^{2}}^{2} + C\|  \partial_{x} \phi_1  \|_{H^{1}(\Gamma)}^{2} \left( \| \phi \|_{L^{2}(\Gamma)}^{2} +  \|   u_1 \|_{L^{2}(\Gamma)}^{2} \right) +  C \|   u_1 \|_{L^{2}(\Gamma)}^{2} + \epsilon\|    u_1 \|_{H^{\frac12}(\Gamma)}^{2}
\nonumber \\
&~\,\,\, +C \left(  \|  \partial_{\bf n}\phi_{2} \|_{H^{\frac32}(\Gamma)}^{2} + 1  \right) \| \partial_{x} \phi \|_{H^{-\frac12}(\Gamma)}^{2}
\nonumber \\
& \le \epsilon   \| \nabla {\bf u} \|_{L^{2}} +  \frac{1}{4}  \| \nabla \mu \|_{L^{2}} +   \frac{1}{4}\| \partial_{\bf n} \phi \|_{L^{2}(\Gamma)}^{2} +  C \left( \| \nabla {\bf u}_{1}\|_{L^{\infty}} + 1 \right)  \| {\bf u}\|_{L^{2}}^{2}
\nonumber \\
&~\,\,\, + C\| \phi\|_{H^{1}}^{2}+ C \epsilon \|  {\bf u} \|_{H^{1}}^{2}  +C \| \partial_{x} \phi \|_{L^{2}}^{2}
\nonumber \\
& \le  C\epsilon   \| \nabla {\bf u} \|_{L^{2}} +  \frac{1}{4}  \| \nabla \mu \|_{L^{2}} +   \frac{1}{4}\| \partial_{\bf n} \phi \|_{L^{2}(\Gamma)}^{2} +  C \left( \| \nabla {\bf u}_{1}\|_{L^{\infty}} + 1 \right)  \| {\bf u}\|_{L^{2}}^{2}+ C\| \phi\|_{H^{1}}^{2},
\end{align}
where we have used \eqref{h-r-1}.

Next, we consider
\begin{align}\label{weiyi-2}
- \langle \Delta \mu, \Delta \phi \rangle = \langle \nabla\Delta \mu, \nabla\phi \rangle - \int_{\Gamma} \Delta \mu \partial_{\bf n} \phi {\rm d}x=: W_{1} + W_{2}.
\end{align}
It follows from $(\ref{OEq})_{3}$ and integration by parts that
\begin{align}\label{W-1}
W_{1} & = \int_\Omega \left[ \nabla\phi_ t + \nabla \left( {\bf u} \cdot \nabla \phi_1  + {\bf u}_2 \cdot \nabla \phi\right) \right] \cdot \nabla\phi{\rm d}x{\rm d}y
\nonumber \\
& = \frac{1}{2}\frac{\rm d}{{\rm d}t} \| \nabla \phi\|_{L^{2}}^{2} +\int_\Omega \left[ \nabla {\bf u} \cdot \nabla \phi_1 +  {\bf u} \cdot \nabla^{2} \phi_1   +  \nabla{\bf u}_2 \cdot \nabla \phi \right] \cdot \nabla\phi{\rm d}x{\rm d}y.
\end{align}
Note that $\Delta \mu = -\partial_{\bf n}\phi - \gamma^{(2)}_{fs}(\theta\phi_{1}+(1-\theta)\phi_{2})\phi$ on $\Gamma$, thus one has
\begin{align}\label{W-2}
W_{2} & = \int_{\Gamma} \left(  \partial_{\bf n} \phi  + \gamma^{(2)}_{fs}(\theta\phi_{1}+(1-\theta)\phi_{2})\phi  \right) \partial_{\bf n} \phi {\rm d}x
\nonumber \\
& =  \| \partial_{\bf n} \phi \|_{L^{2}(\Gamma)}^{2} +\int_{\Gamma}   \gamma^{(2)}_{fs}(\theta\phi_{1}+(1-\theta)\phi_{2})\phi  \partial_{\bf n} \phi {\rm d}x.
\end{align}
We obtain from $(\ref{OEq})_{4}$ that
\begin{align}\label{weiyi-3}
- \langle \Delta \mu, \Delta \phi \rangle &=  \int_\Omega \Delta \mu \left[ \mu-\phi\left( \phi_1^2+\phi_1\phi_2+\phi_2^2 -1\right) \right] {\rm d}x{\rm d}y
 \nonumber \\
 &= -\| \nabla \mu \|_{L^{2}} +\int_\Omega \nabla \mu \cdot \nabla \phi\left( \phi_1^2+\phi_1\phi_2+\phi_2^2 -1\right) {\rm d}x{\rm d}y
 \nonumber \\
 & ~\,\,\, +\int_\Omega \nabla \mu \cdot  \nabla\left( \phi_1^2+\phi_1\phi_2+\phi_2^2 -1\right) \phi  {\rm d}x{\rm d}y.
 \end{align}
Putting $(\ref{W-1})$-$(\ref{weiyi-3})$ into $(\ref{weiyi-2})$, we have
\begin{align}\label{weiyi-4}
& \frac{1}{2}\frac{\rm d}{{\rm d}t} \| \nabla \phi\|_{L^{2}}^{2} +   \| \nabla \mu \|_{L^{2}} + \| \partial_{\bf n} \phi \|_{L^{2}(\Gamma)}^{2}
\nonumber \\
& = -\int_\Omega \left[ \nabla {\bf u} \cdot \nabla \phi_1 +  {\bf u} \cdot \nabla^{2} \phi_1   +  \nabla{\bf u}_2 \cdot \nabla \phi \right] \cdot \nabla\phi{\rm d}x{\rm d}y  
\nonumber \\
& ~\,\,\,+  \int_\Omega \nabla \mu \cdot \nabla \phi\left( \phi_1^2+\phi_1\phi_2+\phi_2^2 -1\right) {\rm d}x{\rm d}y
\nonumber\\
&~\,\,\,+\int_\Omega \nabla \mu \cdot  \nabla\left( \phi_1^2+\phi_1\phi_2+\phi_2^2 -1\right) \phi  {\rm d}x{\rm d}y  -\int_{\Gamma}   \gamma^{(2)}_{fs}(\theta\phi_{1}+(1-\theta)\phi_{2})\phi  \partial_{\bf n} \phi {\rm d}x
\nonumber \\
& \le \| \nabla {\bf u} \|_{L^{2}} \| \nabla \phi_{1} \|_{L^{\infty}} \| \nabla \phi \|_{L^{2}}+\| {\bf u} \|_{L^{4}} \| \nabla^{2}\phi_{1} \|_{L^{4}} \| \nabla \phi \|_{L^{2}} 
\nonumber\\
&~\,\,\, + \| \nabla {\bf u}_{2} \|_{L^{\infty}} \| \nabla \phi \|_{L^{2}}^{2} + \| \nabla \mu \|_{L^{2}}  \| \nabla \phi \|_{L^{2}} \| \phi_1^2+\phi_1\phi_2+\phi_2^2 -1 \|_{L^{\infty}} 
\nonumber \\
&~\,\,\,+  \| \nabla \mu \|_{L^{2}}  \| \phi \|_{L^{4}} \| \nabla(\phi_1^2+\phi_1\phi_2+\phi_2^2 -1) \|_{L^{4}} + C\| \partial_{\bf n} \phi \|_{L^{2}(\Gamma)} \| \phi \|_{L^{2}(\Gamma)}
\nonumber \\
& \le \epsilon  \| \nabla {\bf u} \|_{L^{2}}^{2} + \frac{1}{4}\| \nabla \mu \|_{L^{2}}^{2} +  \frac{1}{4}\| \partial_{\bf n} \phi \|_{L^{2}(\Gamma)}^{2}  + C\| {\bf u} \|_{L^{2}}^{2}
\nonumber \\
& ~\,\,\, + C \| \phi \|_{L^{4}}^{2} \| \nabla(\phi_1^2+\phi_1\phi_2+\phi_2^2 -1) \|_{L^{4}}^{2}+C\| \phi \|_{L^{2}(\Gamma)}^{2}
\nonumber\\
&~\,\,\,+ C \left(  \| \nabla \phi_{1}\|_{L^{\infty}}^{2} + \| \nabla^{2}\phi_{1}\|_{L^{4}}^{2} +  \| \nabla {\bf u}_{2} \|_{L^{\infty}}  + \| \phi_1^2+\phi_1\phi_2+\phi_2^2 -1 \|_{L^{\infty}}^{2} \right) \| \nabla \phi \|_{L^{2}}^{2}
\nonumber \\
& \le \epsilon  \| \nabla {\bf u} \|_{L^{2}}^{2} + \frac{1}{4}\| \nabla \mu \|_{L^{2}}^{2} +  \frac{1}{4}\| \partial_{\bf n} \phi \|_{L^{2}(\Gamma)}^{2}  + C\| {\bf u} \|_{L^{2}}^{2} + C\left( \| \nabla {\bf u}_{2} \|_{L^{\infty}} + 1 \right) \| \phi \|_{H^{1}}^{2}.
\end{align}
Adding $(\ref{weiyi-1})$ to $(\ref{weiyi-4})$, we obtain
\begin{align}\label{weiyi-5}
& \frac{\rm d}{{\rm d}t} \left( \| {\bf u} \|_{L^{2}}^{2}  +  \| \nabla \phi\|_{L^{2}}^{2}  \right) +  \|\mathbb{S}({\bf u}) \|^2_{L^2}  + \beta\| u_1\|^2_{L^2(\Gamma)} +   \| \nabla \mu \|_{L^{2}} + \| \partial_{\bf n} \phi \|_{L^{2}(\Gamma)}^{2}
\nonumber\\
& \le C\epsilon  \| \nabla {\bf u} \|_{L^{2}}^{2} +   C \left( \| \nabla {\bf u}_{1}\|_{L^{\infty}} + \| \nabla {\bf u}_{2} \|_{L^{\infty}}  + 1 \right)  \left( \| {\bf u}\|_{L^{2}}^{2} + \| \phi \|_{H^{1}}^{2}\right).
\end{align}
Note that
\begin{align}
\frac{\rm d}{{\rm d}t}\langle \phi \rangle = \langle \phi_t \rangle & = \langle\Delta\mu  - {\bf u} \cdot \nabla \phi_1  - {\bf u}_2 \cdot \nabla \phi\rangle  = 0, ~\,\,\,{\rm i.e.,}~\,\langle \phi \rangle  = \langle\phi_0\rangle  = 0,
\nonumber
\end{align}
one gets follow from Poincar${\rm \acute{e}}$ inequality that
\begin{align}\label{weiyi-6}
\| \phi\|^2_{L^2} \le  C \| \nabla\phi \|^2_{L^2} + C\left| \langle \phi \rangle\right|^2 =  C \| \nabla\phi \|^2_{L^2}.
\end{align}
Therefore, we deduce from Gronwall inequality, $(\ref{chuzhi})$, $(\ref{weiyi-5})$ and $(\ref{weiyi-6})$ that
\begin{align}
 {\bf u} = \phi =0,~ {\rm a.e. ~in}~\Omega.\nonumber
\end{align}
The proof of Theorem $\ref{TH}$ is completed.
\end{proof}


\section*{Acknowledgments}
Ding's research is supported by the  Key Project of the National Natural Science Foundation of China (No.12131010), the National Natural Science Foundation of China (No.12271032). Li's work is supported by the National Natural Science Foundation of China (No.12371205).  Lin is supported by GuangDong Basic and Applied Basic Research Foundation (No.  2023A1515110500), and Young Teacher Research Fund of South China Normal University (No. 23KJ32).

\section*{Conflict of interest}
On behalf of all authors, the corresponding author states that there is no conflict of interest.




\end{document}